\documentclass[a4paper,12pt]{article}

\usepackage{terek_en}
\usepackage{fullpage}
\title{MAT6702 - Topics in Lorentz Geometry}
\author{Ivo Terek Couto}
\date{}
\newtheorem{problem}{Problem}
\usepackage{fancyhdr}
\fancypagestyle{plain}{}
\usepackage{makeidx}
\usepackage{nicefrac}
\makeindex

\definecolor{DarkGreen}{rgb}{0.0, 0.5, 0.0}
\newcommand\subb[2]{\ensuremath{ \vec{#1}_{\vec{#2}} }}
\newcommand\subg[2]{\ensuremath{ #1_{\vec{#2}} }}
\newcommand{\II}{\ensuremath{\textnormal{I\!I}}}
\newcommand\dd{\ensuremath{\rm d}}
\newcommand{\ptau}{\ensuremath{\textbf{\capricornus}}}
\begin{document}
\maketitle
\vspace*{\fill}
\begin{figure}[H]
  \centering
\includegraphics[width=.7\textwidth]{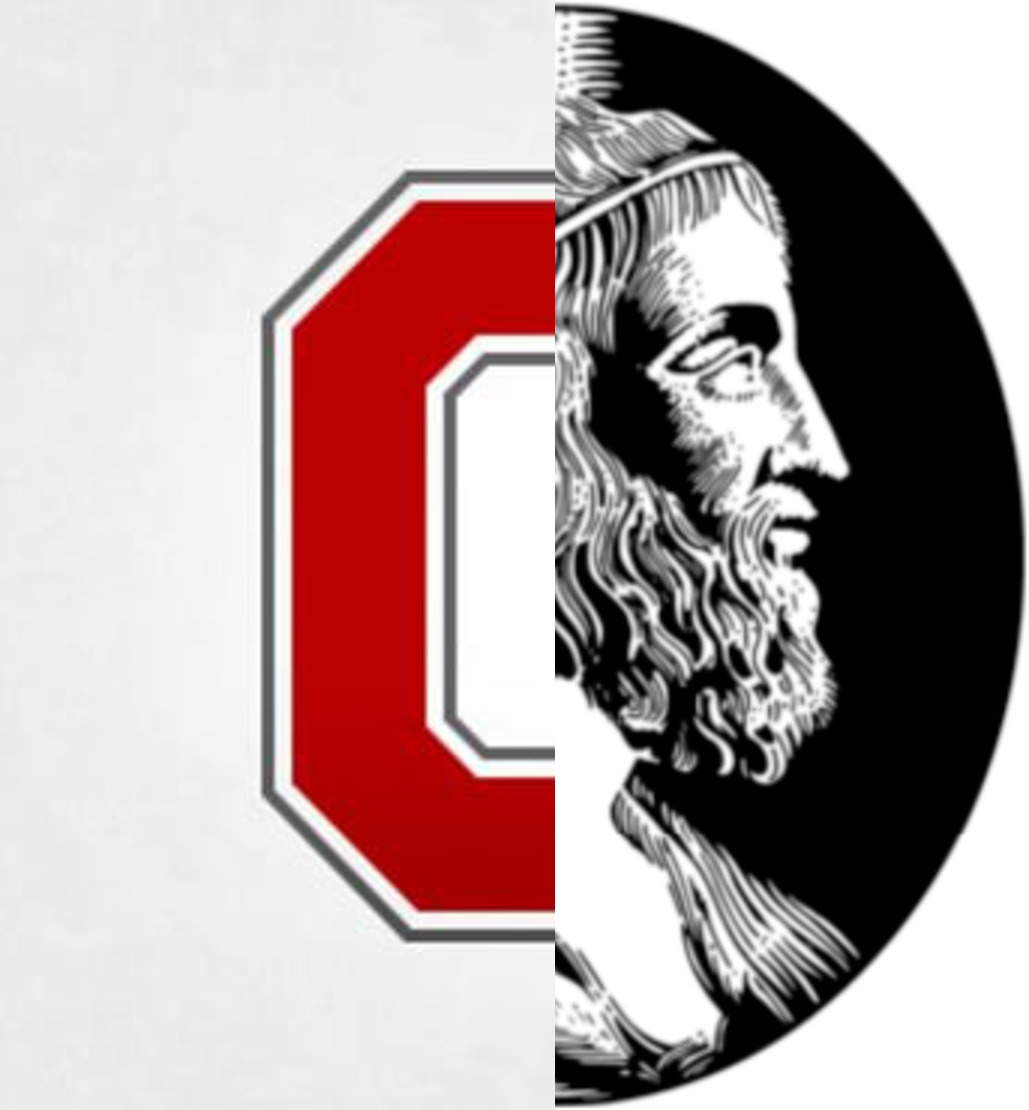}
\end{figure}
\vspace*{\fill}
\vfill

\noindent {\bf Acknowledgement:} This mini-course was supported in part by the departments of Mathematics of \emph{The Ohio State University} and of the \emph{University of S\~{a}o Paulo}, and in part by a {\sf FAPESP-OSU} 2015 Regular Research Award ({\sf FAPESP} grant: 2015/50265-6). 

\newpage

\tableofcontents

\newpage 

\listoffigures\addcontentsline{toc}{section}{List of Figures}

\newpage 

\section*{Foreword}\addcontentsline{toc}{section}{Foreword}

The present text was prepared for the mini-course \emph{MAT6702 - Topics in Lorentz Geometry}, to be taught at the University of S\~{a}o Paulo, during the week from 03/11/19 to 03/15/19. Due to time constraints, some very interesting topics (such as Lorentz boosts, the proof of the classification of matrices in $\mbox{\small ${\rm O}_1^{+\uparrow}(3,\R)$}$, and Bonnet rotations for timelike surfaces) unfortunately had to be left out, but a list of references is provided in the end. As an attempt to engage the reader actively on what is happening here, a few problems are suggested in the end of each section.

In general, the content of these notes is very introductory and meant to be a stepping stone for those interested in learning the subject without yet having advanced background, avoiding the ``heavier'' language of differentiable manifolds and assuming only some knowledge of multivariable calculus, linear algebra, and differential geometry of curves and surfaces in $\R^3$ (on the level of~\cite{dC1} or~\cite{Keti} should be enough).

For this reason, instead of focusing on the similarities between Euclidean space $\R^3$ and Lorentz-Minkowski space $\LM^3$, we will devote our little time together in trying to grasp some of the most striking differences between those ambients.

I hope you enjoy reading this, and if you learn anything new at all here, it was worth the effort.

\vskip 1.5 cm							
\begin{flushright}				
Columbus, March of 2019				
\end{flushright}					

\vskip 0.5 cm							
\begin{flushright}				
  Ivo Terek Couto
\end{flushright}					

\newpage

\section{The spaces $\R^n_\nu$}

\subsection{Basic definitions}

\begin{defn}
  Let $n > 0$ and $0 \leq \nu \leq n$ be non-negative integers. The \emph{pseudo-Euclidean space of index $\nu$}\index{Pseudo-Euclidean!space of index $\nu$} is the pair $\R^n_\nu\doteq (\R^n, \pair{\cdot,\cdot}_\nu)$, where the scalar product\linebreak[4] $\pair{\cdot,\cdot}_\nu\colon \R^n\times \R^n \to \R$ is given by \[ \pair{\vec{x},\vec{y}}_\nu \doteq x_1y_1 + \cdots + x_{n-\nu}y_{n-\nu} - x_{n-\nu+1}y_{n-\nu+1} - \cdots - x_ny_n.  \]
\end{defn}

Particular cases are the usual Euclidean space $\R^n_0 \equiv \R^n$ and the \emph{Lorentz-Minkowski space}\index{Lorentz-Minkowski space} $\LM^n \equiv \R^n_1$, whose products are then denoted simply $\pair{\cdot,\cdot}_E$ and $\pair{\cdot,\cdot}_L$, respectively.

Regarding vectors in $\R^n$ as column-vectors, one may write $\pair{\vec{x},\vec{y}}_\nu = \vec{x}^\top {\rm Id}_{n-\nu,\nu}\vec{y}$, where the \emph{identity matrix of index $\nu$}\index{Identity matrix of index $\nu$} is \[\mbox{\large ${\rm Id}_{n-\nu,\nu}$} = (\eta_{ij}^\nu)_{i,j=1}^n\doteq  \left(
    \begin{array}{c|c}
      \mbox{\large ${\rm Id}_{n-\nu}$} &  \mbox{\large $0$} \\[.7ex] \hline \\[-1.5ex]
      \mbox{\large $0$} & \mbox{\large $-{\rm Id}_\nu$}\\
    \end{array}
  \right).\]
Note that the product $\pair{\cdot,\cdot}_\nu$ is not positive-definite anymore, which is an obstacle for defining a norm $\|\cdot\|_\nu$. We will insist on trying, and setting $\|\vec{x}\|_\nu \doteq \sqrt{|\pair{\vec{x},\vec{x}}_\nu|}$ anyway. This ``fake norm'' $\|\cdot\|_\nu$ has poor properties -- we'll see a couple of them soon. Despite this perhaps-not-so-small issue, the product $\pair{\cdot,\cdot}_\nu$ has the one property that allows us to develop the theory to some extent: non-degenerability. That is to say, if $\pair{\vec{x},\vec{y}}_\nu = 0$ for every $\vec{y} \in \R^n_\nu$, we necessarily have $\vec{x}=\vec{0}$. Or in other words, the induced map $\R^n_\nu \ni \vec{x} \mapsto \pair{\vec{x},\cdot}_\nu \in (\R^n_\nu)^*$ is an isomorphism. Having lost the positivity of $\pair{\cdot,\cdot}_\nu$, it is convenient to sort vectors in $\R^n_\nu$ in three classes:

\begin{defn}[Causal character\index{Causal!character}]
  A non-zero vector $\vec{x} \in \R^n_\nu$ is called:
  \begin{itemize}
  \item \emph{spacelike} if $\pair{\vec{x},\vec{x}}_\nu > 0$.
  \item \emph{timelike} if $\pair{\vec{x},\vec{x}}_\nu < 0$.
  \item \emph{lightlike} if $\pair{\vec{x},\vec{x}}_\nu = 0$.
  \end{itemize}
The \emph{indicator}\index{Indicator!of a vector} of $\vec{x}$ is $1$, $-1$ or $0$ according to the causal type of $\vec{x}$, and it is denoted by $\epsilon_{\vec{x}}$.
\end{defn}

\begin{Ex}
  If ${\rm can} = (\vec{e}_1,\ldots,\vec{e}_n)$ is the standard basis of $\R^n_\nu$, then $\vec{e}_i$ is spacelike for $1 \leq i \leq n-\nu$ and timelike for $n-\nu+1 < i \leq n$. If $1 \leq i \leq n-\nu < j \leq n$, then $\vec{e}_i \pm \vec{e}_j$ is lightlike. In $\LM^2$ and $\LM^3$, we can actually make some sketches based in the equations $x^2-y^2 = c$ and $x^2+y^2-z^2=c$ (for positive, negative or zero $c$):

\begin{figure}[H]
  \centering
  \hspace*{1cm}
  \includegraphics[width=.25\linewidth]{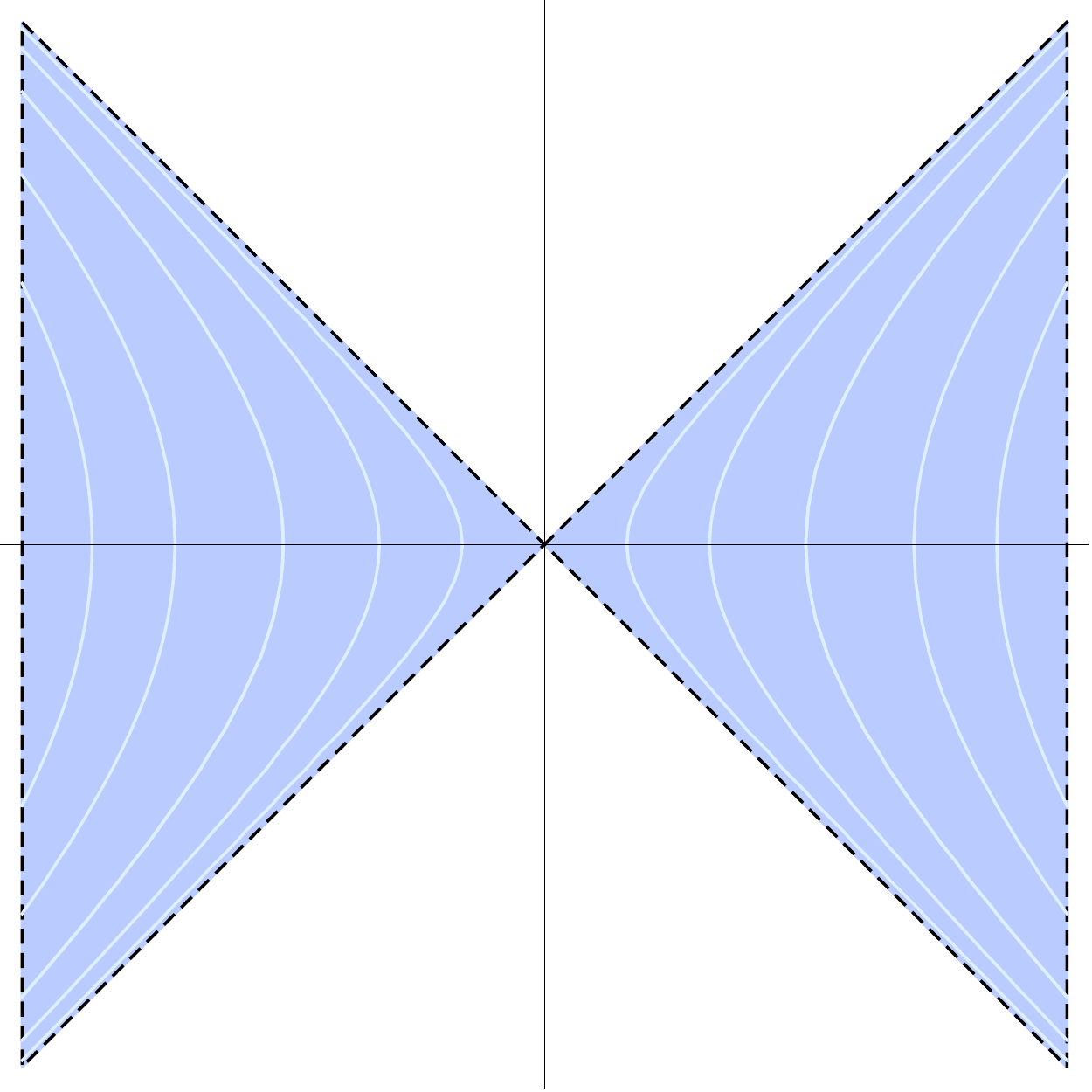}
  \includegraphics[width=.25\linewidth]{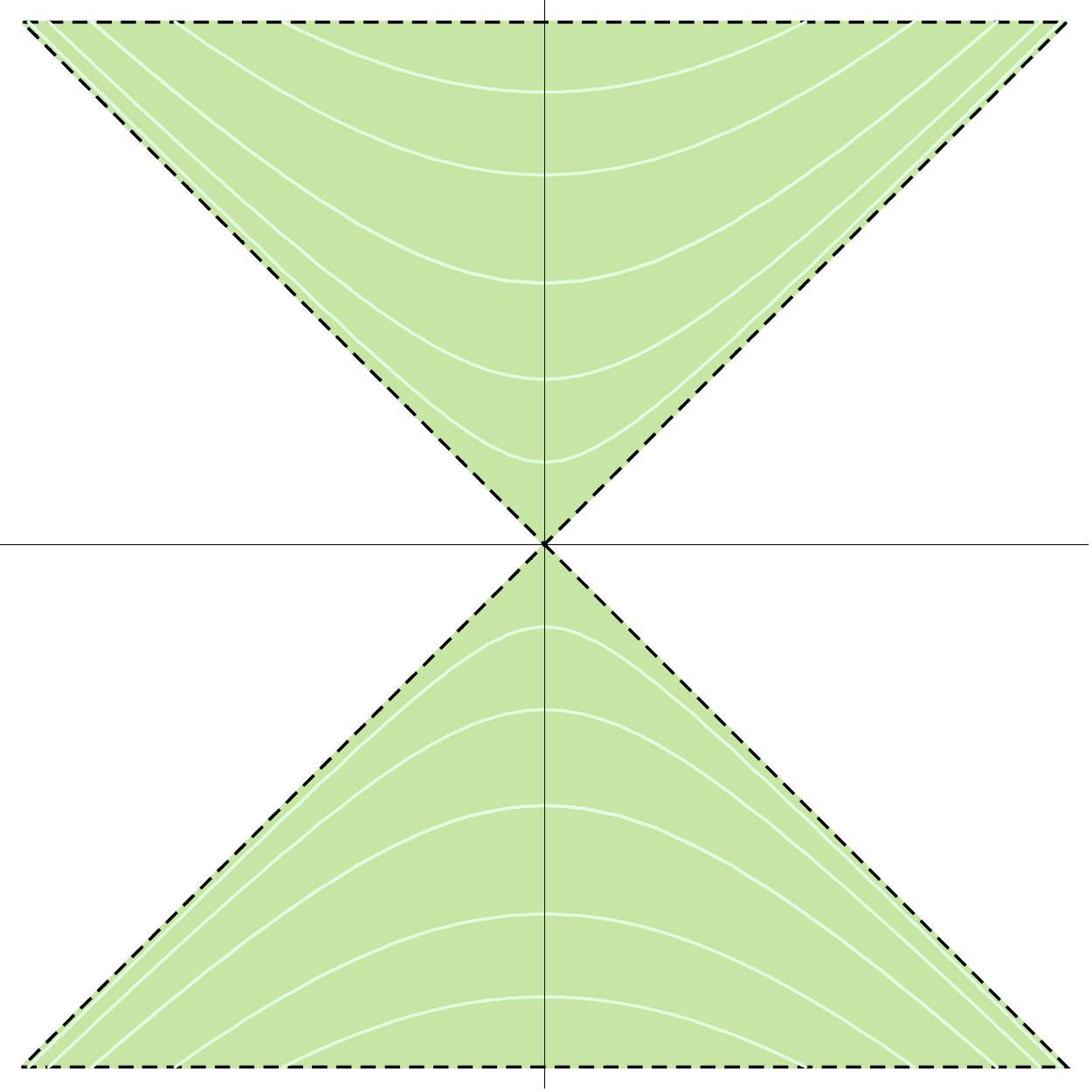}
  \includegraphics[width=.25\linewidth]{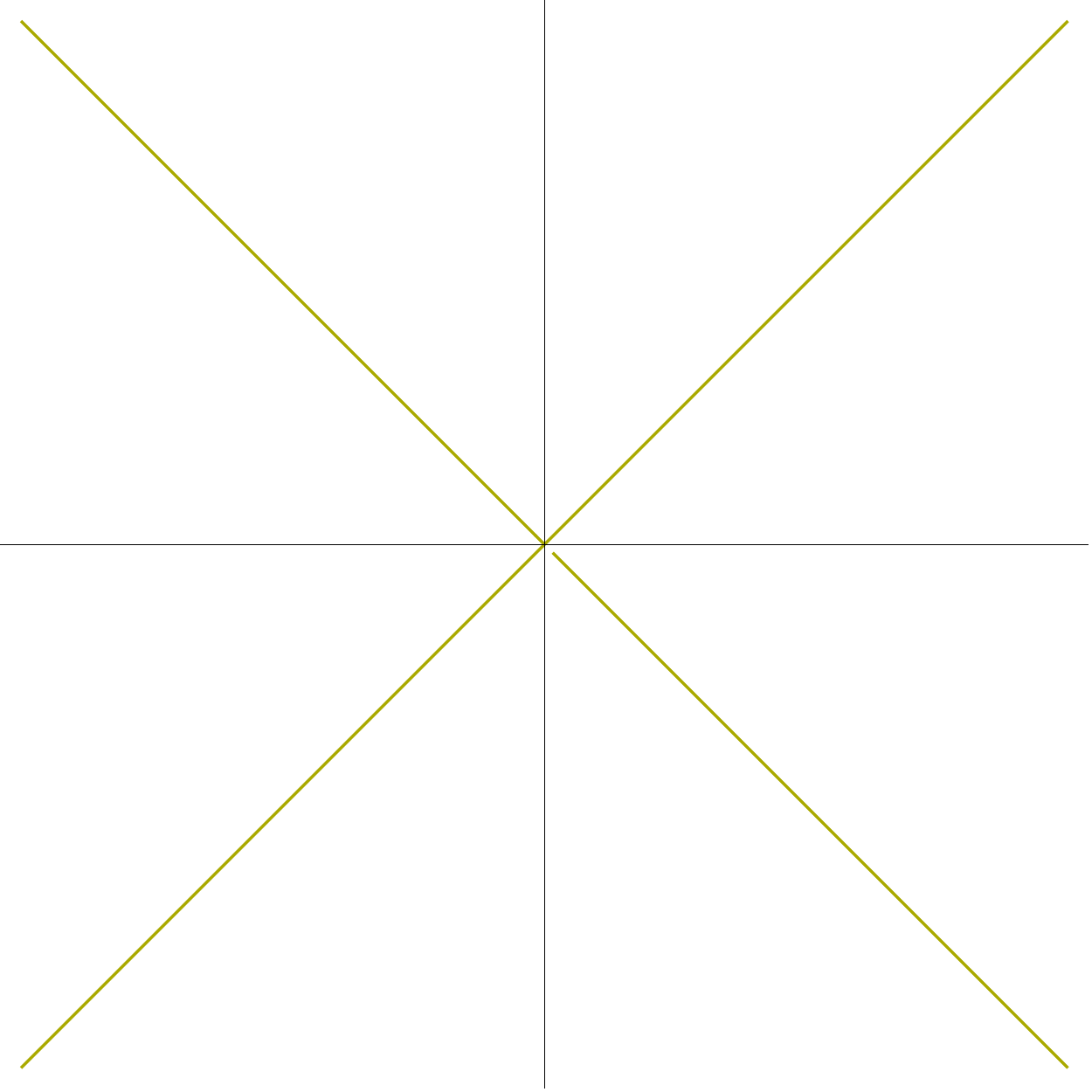}
  \caption{Causal ``regions'' in $\LM^2$.}
\end{figure}

\begin{figure}[H]
  \centering
  \begin{picture}(0,0)%
\includegraphics{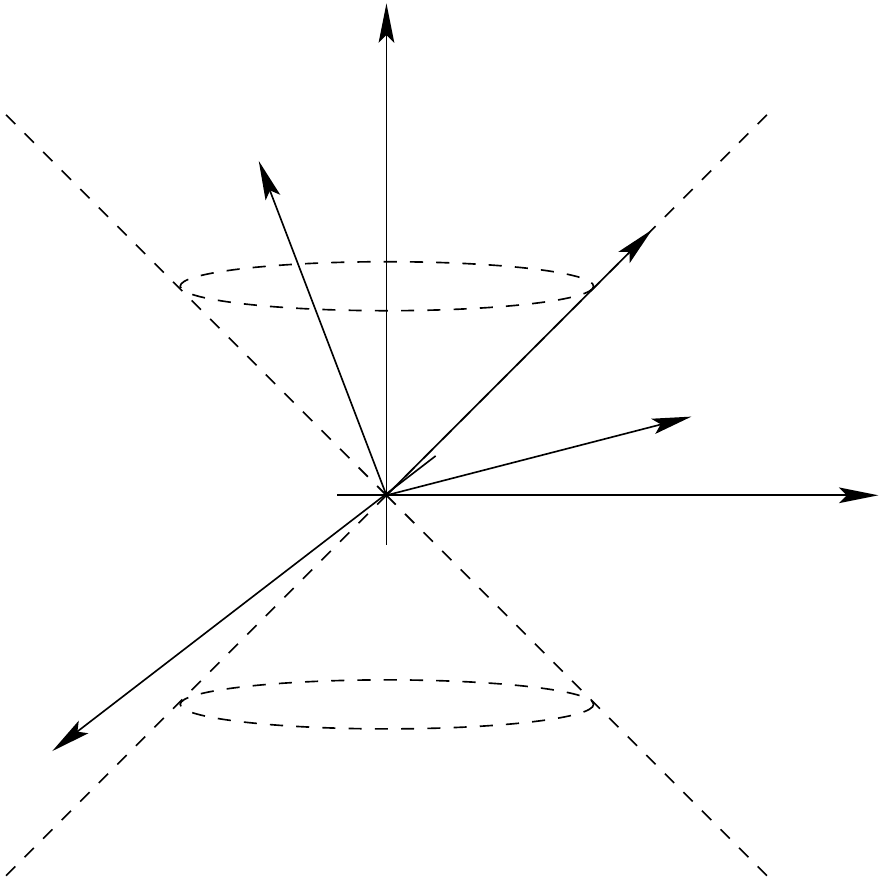}%
\end{picture}%
\setlength{\unitlength}{4144sp}%
\begingroup\makeatletter\ifx\SetFigFont\undefined%
\gdef\SetFigFont#1#2#3#4#5{%
  \reset@font\fontsize{#1}{#2pt}%
  \fontfamily{#3}\fontseries{#4}\fontshape{#5}%
  \selectfont}%
\fi\endgroup%
\begin{picture}(4029,4029)(-1766,-3853)
\put(-854,-466){\makebox(0,0)[lb]{\smash{{\SetFigFont{12}{14.4}{\rmdefault}{\mddefault}{\updefault}{\color[rgb]{0,0,0}timelike}%
}}}}
\put(1441,-1771){\makebox(0,0)[lb]{\smash{{\SetFigFont{12}{14.4}{\rmdefault}{\mddefault}{\updefault}{\color[rgb]{0,0,0}spacelike}%
}}}}
\put(1261,-1006){\makebox(0,0)[lb]{\smash{{\SetFigFont{12}{14.4}{\rmdefault}{\mddefault}{\updefault}{\color[rgb]{0,0,0}lightlike}%
}}}}
\end{picture}
  \caption{Causal types in $\LM^3$}
  \label{fig:types}
\end{figure}

\begin{figure}[H]
  \centering
  \hspace*{.5cm}
  \includegraphics[width=.250\linewidth]{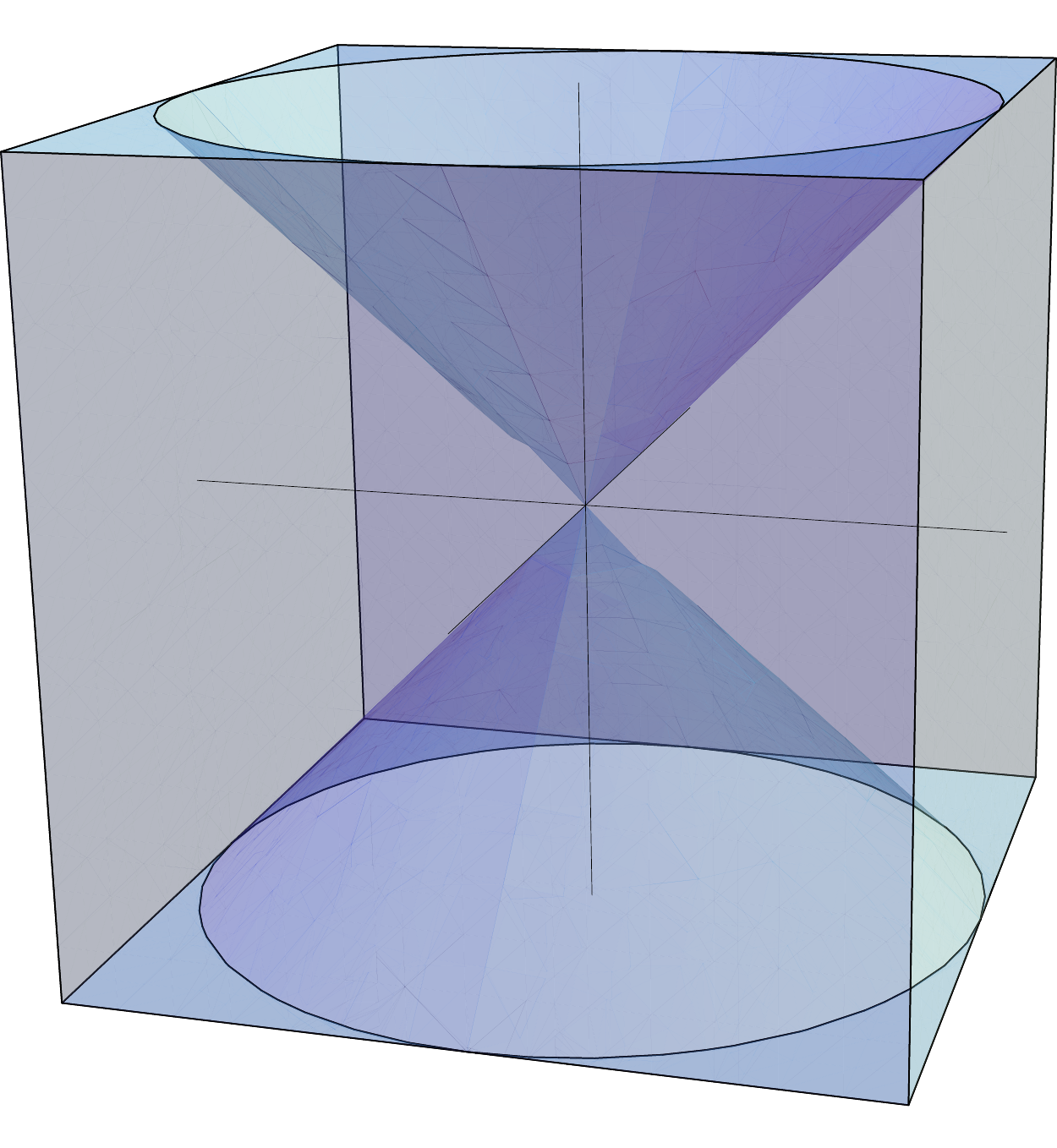}
  \includegraphics[width=.216\linewidth]{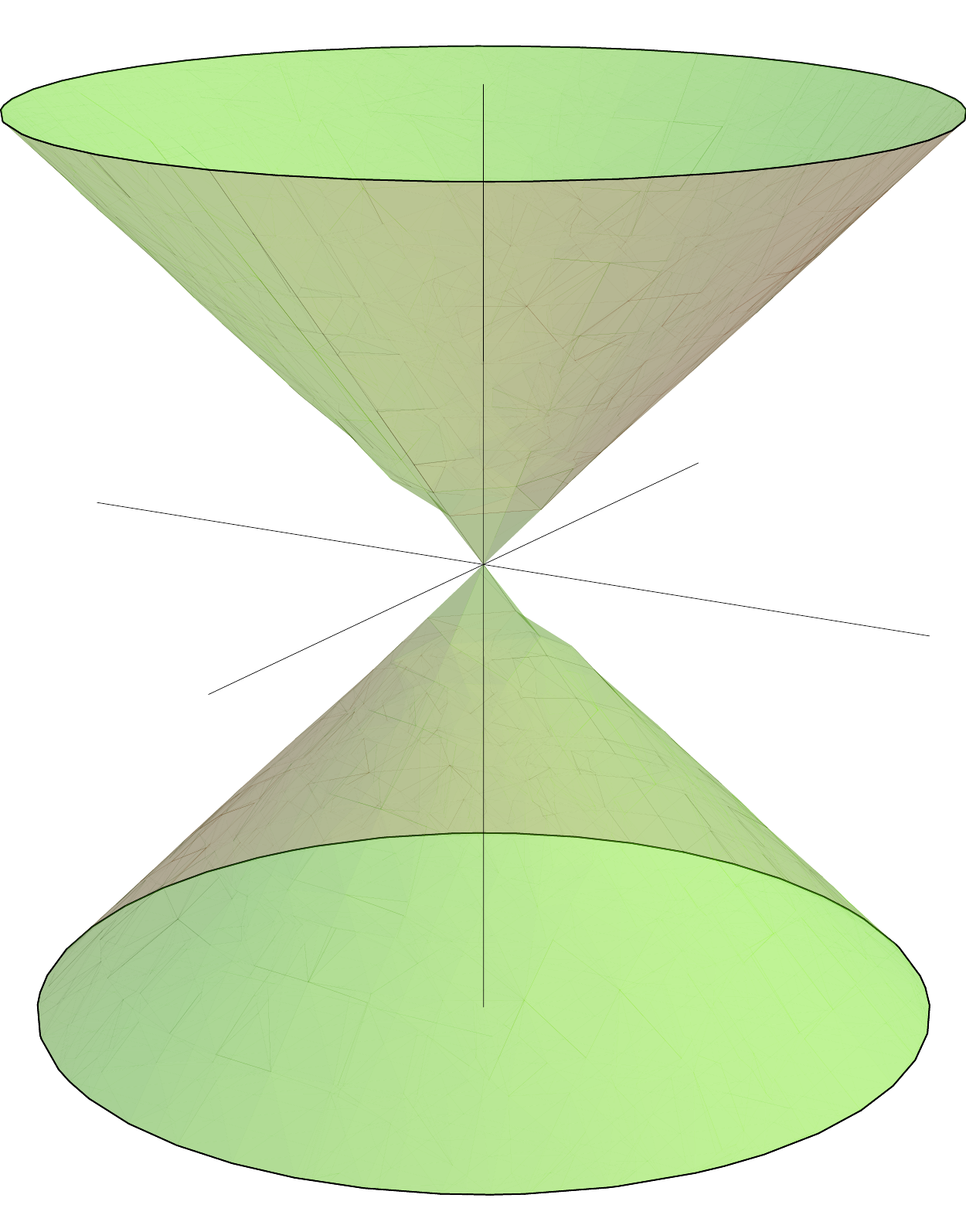}
  \includegraphics[width=.283\linewidth]{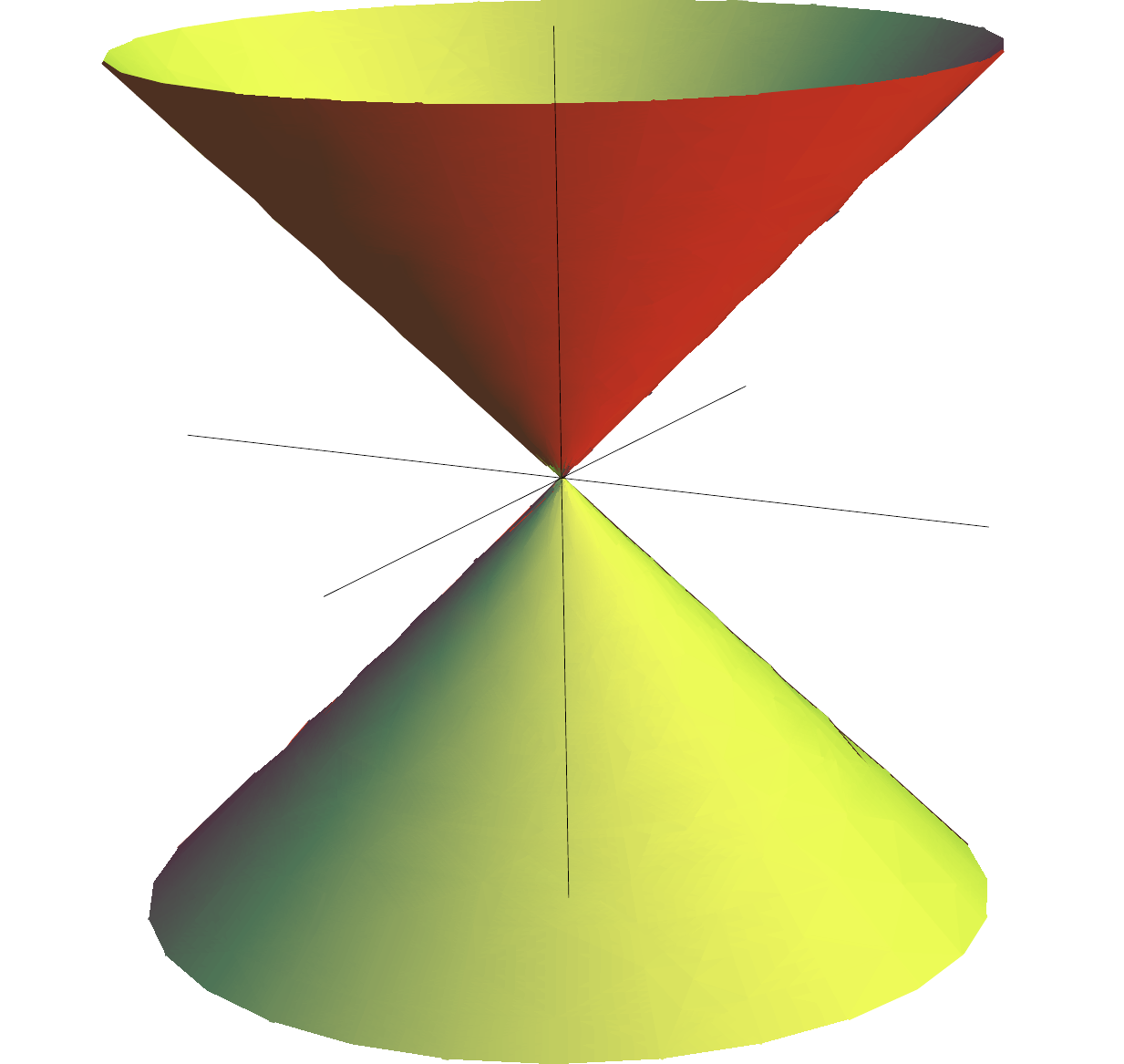}
  \caption{Causal ``regions'' in $\LM^3$.}
\end{figure}

Not surprisingly, we call the collection of all lightlike vectors in the space its \emph{light cone}.

\end{Ex}

Soon we will generalize the notion of causal character for other objects than vectors, such as subspaces, curves and surfaces. One of the fundamental concepts in geometry is the one of orthogonality. So:

\begin{defn}
  Two vectors $\vec{x},\vec{y} \in \R^n_\nu$ are \emph{$\nu$-orthogonal}\index{$\nu$-orthogonality} if $\pair{\vec{x},\vec{y}}_\nu =0$. A basis for $\R^n_\nu$ is called $\nu$-orthogonal if all its vectors are pairwise $\nu$-orthogonal, and it is said to be $\nu$-orthonormal if all its vectors have scalar square $1$ or $-1$. For $\nu = 1$, one usually uses the term ``Lorentz-orthogonal'' instead, and if there is no risk of confusion, we'll do away with the $\nu$. 
\end{defn}

\begin{defn}
  Let $S\subseteq \R^n_\nu$ be any set. Let's say that \[S^\perp = \{\vec{x} \in \R^n_\nu \mid \pair{\vec{x},\vec{y}}_\nu = 0 \mbox{ for all }\vec{y} \in S\}\] is the \emph{subspace of $\R^n_\nu$ orthogonal to $S$}\index{Orthogonal!space}.
\end{defn}

\begin{obs}
  $S^\perp$ is a vector subspace of $\R^n_\nu$ even when $S$ is not.
\end{obs}

We avoid the name ``orthogonal complement'' because when $S$ is a subspace of $\R^n_\nu$ we might not have $S\oplus S^\perp = \R^n_\nu$. For example, in $\LM^2$, the line $S$ spanned by the lightlike vector $(1,1)$ satisfies $S = S^\perp = S+S^\perp$. So a natural question should be: \emph{when} do we have $S\oplus S^\perp = \R^n_\nu$? We start with the recomforting result:

\begin{prop}
  Let $S\subseteq \R^n_\nu$ be a subspace. Then $\dim S + \dim S^\perp = n$ and $(S^\perp)^\perp = S$.
\end{prop}

\begin{dem}
  The map $\R^n_\nu \ni \vec{x} \mapsto \pair{\vec{x},\cdot}\big|_S \in S^*$ is linear, surjective (since $\pair{\cdot,\cdot}_\nu$ is non-degenerate), and its kernel is $S^\perp$. So the dimension formula follows from the rank-nullity theorem. Said formula applied twice also says that $\dim S = \dim (S^\perp)^\perp$, so $S \subseteq (S^\perp)^\perp$ implies $S = (S^\perp)^\perp$.
\end{dem}

With this, we may also conclude the:

\begin{cor}\label{cor:non-deg}
  Let $S\subseteq \R^n_\nu$ be a subspace. Then $S\oplus S^\perp=\R^n_\nu$ if and only if $S$ is non-degenerate (i.e., ${\pair{\cdot,\cdot}_\nu}\big|_S$ is non-degenerate). It also follows that $S$ is non-degenerate if and only if $S^\perp$ is also non-degenerate.
\end{cor}

\begin{dem}
  From $\dim(S+S^\perp) = \dim S + \dim S^\perp - \dim(S\cap S^\perp) = n-\dim(S\cap S^\perp)$ it follows that $S+S^\perp = \R^n_\nu$ if and only if $S \cap S^\perp = \{\vec{0}\}$, which in turn is equivalent to ${\pair{\cdot,\cdot}_\nu}\big|_S$ being non-degenerate.
\end{dem}

This means that we may define \emph{orthogonal projections}\index{Orthogonal!projection} only onto non-degenerate subspaces. Back to the previous example, we may now see what went wrong there: the line spanned by $(1,1)$ in $\LM^2$ is degenerate, since $\pair{(1,1), (\lambda,\lambda)}_L = 0$ for all $\lambda \in \R$. In $\LM^n$, we may stick to the causal type terminology previously used:

\begin{defn}
  Let $S\subseteq \LM^n$ be a non-trivial vector subspace. We say that $S$ is:
  \begin{itemize}
  \item \emph{spacelike} if $\pair{\cdot,\cdot}_L\big|_S$ is positive-definite;
  \item \emph{timelike} if $\pair{\cdot,\cdot}_L\big|_S$ is negative-definite, or indefinite and non-degenerate;
  \item \emph{lightlike} if $\pair{\cdot,\cdot}_L\big|_S$ is degenerate;
  \end{itemize}
\end{defn}

\begin{obs}
  In $\R^n_\nu$, one might also say that $S$ is timelike if $\pair{\cdot,\cdot}_\nu\big|_S$ is negative definite, but if $\nu>1$ this does not have the same physical appeal (which we'll get to in the next section) as in $\LM^n$. Moreover, one can say that $S$ is \emph{null} if $\pair{\cdot,\cdot}_\nu\big|_S = 0$, which in $\LM^n$ is the same as $S$ being lightlike and one-dimensional.
\end{obs}

Here's the relation between causal characters of subspaces and orthogonality:

\begin{teo}\label{teo:causal_perp}
  Let $S\subseteq \LM^n$ be a subspace. Then $S$ is spacelike if and only if $S^\perp$ is timelike; $S$ is lightlike if and only $S^\perp$ is also lightlike.
\end{teo}

\begin{dem}
  The second part of the result is nothing more than a restatement of Corollary \ref{cor:non-deg}, which will also be used to prove the first part. Assume that $S$ is spacelike. So $S$ is non-degenerate and we write $\LM^n = S\oplus S^\perp$. Then $S^\perp$ must necessarily contain a timelike vector because $\LM^n$ does - more precisely, take $\vec{v} \in \LM^n$ timelike and write $\vec{v} = \vec{x}+\vec{y}$ with $\vec{x}\in S$ and $\vec{y} \in S^\perp$, so that $\pair{\vec{x},\vec{x}}_L + \pair{\vec{y},\vec{y}}_L = \pair{\vec{v},\vec{v}}_L < 0$ with $\pair{\vec{x},\vec{x}}_L\geq 0$ forces $\vec{y} \in S^\perp$ to be timelike. Conversely, assume now that $S$ is timelike, and take $\vec{u} \in S$ timelike. Since $S^\perp \subseteq \vec{u}^\perp$, it suffices now to show that $\vec{u}^\perp$ is spacelike. We know again from Corollary \ref{cor:non-deg} that $\vec{u}^\perp$ is not lightlike, and if we have $\vec{v} \in \vec{u}^\perp$ timelike, the plane spanned by $\vec{u}$ and $\vec{v}$ in $\LM^n$ has dimension $2$ while being negative-definite, which is impossible.
\end{dem}

Here's another important result:

\begin{teo}
  Let $S\subseteq \R^n_\nu$ be a non-degenerate subspace. Then $S$ has an orthogonal basis.
\end{teo}

\begin{dem}
By induction. By hypothesis we may take $\vec{u} \in S$ with $\pair{\vec{u},\vec{u}}_\nu \neq 0$. Then the orthogonal complement of $\vec{u}$ in $S$ is non-degenerate and has dimension one lower. Take an orthogonal basis for this complement and add $\vec{u}$ to this list. Fill any details you may want.
\end{dem}

We can conclude this section with some results about linear independence, in general:

\begin{teo}\label{teo:light_LD}
  Let $\vec{u}_1,\ldots,\vec{u}_{\nu+1}\in \R^n_\nu$ be pairwise lightlike orthogonal vectors. Then we have that $(\vec{u}_1,\ldots,\vec{u}_{\nu+1})$ is linearly dependent.
\end{teo}

\begin{dem}
  The space $\R^n_\nu$ has a natural decomposition as $\R^n_\nu = \R^{n-\nu}\oplus \R^\nu_\nu$, so that for the standard basis ${\rm can} = (\vec{e}_1,\ldots,\vec{e}_n)$ of $\R^n_\nu$, we may decompose \[  \vec{u}_j = \vec{x}_j + \sum_{i=1}^\nu a_{ij}e_{n-\nu+i}, \qquad 1 \leq j \leq \nu+1, \]for some vectors $\vec{x}_j \in \R^{n-\nu}\oplus \{\vec{0}\}$ and real coefficients $a_{ij}$, which actually define a linear map $A\colon \R^{\nu+1}\to \R^\nu$. The condition $\pair{\vec{u}_i,\vec{u}_j}_\nu=0$  readily implies the equality $\pair{\vec{x}_i,\vec{x}_j}_\nu = \sum_{k=1}^\nu a_{ki}a_{kj}$, for all $1 \leq i,j \leq \nu+1$. For dimensional reasons, we may also choose a non-zero vector $\vec{b} = (b_i)_{i=1}^{\nu+1} \in \ker A$. Putting all of this together, we see that \[  \left\langle \sum_{i=1}^{\nu+1} b_i\vec{x}_i, \sum_{j=1}^{\nu+1} b_j\vec{x}_j\right\rangle_\nu = \sum_{i,j=1}^{\nu+1}b_ib_j \sum_{k=1}^\nu a_{ki}a_{kj} = (A\vec{b})^\top(A\vec{b}) = 0.  \]However, the combination $\sum_{i=1}^{\nu+1}b_i\vec{x}_i$ lies in the spacelike subspace $\R^{n-\nu}\oplus \{\vec{0}\}$, so the above gives $\sum_{i=1}^{\nu+1}b_j\vec{x}_j = \vec{0}$. So, $\vec{b} \in \ker A$ now gives us \[ \sum_{j=1}^{\nu+1} b_j\vec{u}_j  = \sum_{j=1}^{\nu+1}b_j\vec{x}_j + \sum_{i=1}^\nu \left(\sum_{j=1}^{\nu+1} b_ja_{ij}\right)\vec{e}_{n-\nu+1} = \vec{0}+\vec{0}=\vec{0},  \]as wanted.
\end{dem}

As a corollary, we obtain one of the most striking differences between Euclidean and Lorentzian geometry:

\begin{cor}\label{cor:lightlike_parallel}
  Two lightlike vectors in $\LM^n$ are Lorentz-orthogonal if and only if they are parallel.
\end{cor}

The previous proof might also hint that the matrix of coefficients of $\pair{\cdot,\cdot}_\nu$ with respect to a given basis (also called the \emph{Gram matrix}\index{Gram matrix} of $\pair{\cdot,\cdot}_\nu$ with respect to said basis)  will play an important role in this whole theory.

\begin{prop}
  Let $\vec{u}_1,\ldots,\vec{u}_m \in \R^n_\nu$ be vectors such that the Gram matrix $(\pair{\vec{u}_i,\vec{u}_j}_\nu)_{i,j=1}^m$ is invertible. Then $(\vec{u}_1,\ldots, \vec{u}_m)$ is linearly independent.
\end{prop}

\begin{dem}
  Write $\sum_{i=1}^m a_i\vec{u}_i =\vec{0}$ and apply $\pair{\cdot,\vec{u}_j}_\nu$ on both sides to get $\sum_{i=1}^m a_i \pair{\vec{u}_i,\vec{u}_j}_\nu = 0$. The hypothesis then implies that $a_1=\cdots=a_m=0$ as wanted.
\end{dem}

We know that for the usual Euclidean product in $\R^n$ the converse to the above result is true. It is not true, in general, in the pseudo-Euclidean spaces $\R^n_\nu$. As an extreme counter-example, take any (non-zero) lightlike vector: it is linearly independent by itself, but its $1\times 1$ Gram matrix is just $(0)$. As disappointing as this might be, this means that we'll have to add some extra conditions for this converse to hold. This leads us to what we may call the ``non-degenerability chain conditions''. Here is an example:

\begin{prop}\label{prop:converse_Gram}
  Let $(\vec{u}_1,\ldots, \vec{u}_m)$ be a $m$-uple of linearly independent vectors in $\R^n_\nu$ such that each intermediate subspace ${\rm span}(\vec{u}_1,\ldots,\vec{u}_k)$ is non-degenerate, for $1 \leq k \leq m$. Then the Gram matrix $(\pair{\vec{u}_i,\vec{u}_j}_\nu)_{i,j=1}^m$ is invertible.
\end{prop}

Another example of this non-degenerability chain condition is related to the \emph{Gram-Schmidt orthogonalization process}\index{Gram-Schmidt process}. Namely, if we start with linearly independent vectors $(\vec{u}_1,\ldots,\vec{u}_m)$ and try to produce from these vectors another set of orthogonal vectors $(\widetilde{\vec{u}_1},\ldots, \widetilde{\vec{u}_m})$ spanning the same subspace, at least in the Euclidean case we would proceed inductively, by setting \[ \widetilde{\vec{u}_{k+1}} = \vec{u}_k - \sum_{i=1}^k \frac{\pair{\vec{u}_{k+1}, \widetilde{\vec{u}_i}}}{\|\widetilde{\vec{u}_i}\|^2}\widetilde{\vec{u}_i}.  \]In the pseudo-Euclidean case, not only we need to take into account the causal character of each $\widetilde{\vec{u}_i}$, but we need to ensure that none of those vectors are lightlike. The condition that all the intermediate subspaces ${\rm span}(\vec{u}_1,\ldots,\vec{u}_k)$ are non-degenerate, for $1 \leq k \leq m$, is again precisely what we need to safely do \[ \widetilde{\vec{u}_{k+1}} = \vec{u}_k - \sum_{i=1}^k \epsilon_{\widetilde{\vec{u}_i}}\frac{\pair{\vec{u}_{k+1}, \widetilde{\vec{u}_i}}_\nu}{\|\widetilde{\vec{u}_i}\|_\nu^2}\widetilde{\vec{u}_i}  \]in $\R^n_\nu$. Usually it is a bad idea to insist on using the ``fake norm'' $\|\cdot\|_\nu$: we'll try to avoid the absolute values the most we can. So we may alternatively write \[ \widetilde{\vec{u}_{k+1}} = \vec{u}_k - \sum_{i=1}^k \frac{\pair{\vec{u}_{k+1}, \widetilde{\vec{u}_i}}_\nu}{\pair{\widetilde{\vec{u}_i},\widetilde{\vec{u}_i}}_\nu}\widetilde{\vec{u}_i}  \]instead, which automatically takes into account the indicators of the $\widetilde{\vec{u}}_i$. For the proof of Proposition \ref{prop:converse_Gram} above and more details about the adapted Gram-Schmidt process, see~\cite{TL}. When we start discussing curve theory, we will see that we'll have three classes of curves: the admissible curves, the lightlike curves and the semi-lightlike curves. The latter two require some special treatment precisely because they fail to respect a certain non-degenerability chain condition (but you might have guessed this by now). We move on.

\subsection{Pseudo-orthogonal transformations}

When studying the geometry of any scalar product, it is essential to understand the transformations of the ambient space which preserve said product:

\begin{defn}
  A linear transformation $\Lambda\colon \R^n_\nu \to \R^n_\nu$ such that $\pair{\Lambda\vec{x},\Lambda\vec{y}}_\nu = \pair{\vec{x},\vec{y}}_\nu$ for all $\vec{x},\vec{y} \in \R^n_\nu$ is called a \emph{pseudo-orthogonal transformation}\index{Pseudo-orthogonal transformation}. We denote the collection of these transformations, maybe not surprisingly, by ${\rm O}_\nu(n,\R)$. When $\nu = 1$, $\Lambda$ is called a \emph{Lorentz transformation} and ${\rm O}_1(n,\R)$ is called the \emph{Lorentz group}.
\end{defn}

Let's get the following simple characterization out of the way:

\begin{prop}
  Let $\Lambda\colon \R^n_\nu \to \R^n_\nu$ be a linear transformation. Then $\Lambda \in {\rm O}_\nu(n,\R)$ if and only if $\Lambda^\top {\rm Id}_{n-\nu,\nu}\Lambda = {\rm Id}_{n-\nu,\nu}$. It follows that $\det \Lambda = \pm 1$, and so $\Lambda$ is an isomorphism.
\end{prop}

\begin{obs}
  Another way to state the above is saying that the rows and columns of $\Lambda$ form orthonormal bases of $\R^n_\nu$. This proposition also implies that ${\rm O}_\nu(n,\R)$ is a group closed under matrix transposition (proof?).
\end{obs}

\begin{Ex}
  Given $\varphi > 0$, the \emph{hyperbolic rotation}\index{Hyperbolic!rotation} $R_\varphi^h\colon \LM^2 \to \LM^2$ given by \[  R_\varphi^h(x,y) = (x\cosh\varphi + y\sinh\varphi, x\sinh \varphi + y\cosh \varphi)  \]is a Lorentz transformation, whose inverse is naturally $(R^h_\varphi)^{-1} = R^h_{-\varphi}$ (you should check this if you don't immediately believe it, it is instructive). Up to a couple of signs, this is actually the only Lorentz transformation in dimension $2$. We'll come back to that in Theorem \ref{teo:o12}.
\end{Ex}

In general, in the same way that a rigid motion of $\R^n$ is always the composition of an orthogonal map and a translation, the corresponding notion of ``rigid motion'' in $\R^n_\nu$ also has this property. Rewriting the definition of a rigid motion in $\R^n$ without employing $\|\cdot\|$ leads to the:

\begin{defn}
  A \emph{pseudo-Euclidean isometry}\index{Pseudo-Euclidean!isometry} in $\R^n_\nu$ is a map $F\colon \R^n_\nu\to \R^n_\nu$ such that \[  \pair{F(\vec{x})-F(\vec{y}), F(\vec{x})-F(\vec{y})}_\nu = \pair{\vec{x}-\vec{y},\vec{x}-\vec{y}}_\nu, \]for all $\vec{x},\vec{y} \in \R^n_\nu$. The collection of such maps is denoted by ${\rm E}_\nu(n,\R)$. When $\nu = 1$, $F$ is called a \emph{Poincar\'{e} transformation} and ${\rm P}(n,\R) = {\rm E}_1(n,\R)$ is called the \emph{Poincar\'{e} group}.
\end{defn}

To justify the name ``Poincar\'{e} group'', one has to check that pseudo-Euclidean isometries are indeed invertible, and that its inverse is also a pseudo-Euclidean isometry. One possible way to do this is actually going over and beyond, and classifying these maps. We can even say that the above definition was written precisely so that the same strategy used in proving that every rigid motion in $\R^n$ is the composition of a translation and an orthogonal map works. As such, we won't provide a full proof, but the main steps:

\begin{enumerate}[(i)]
\item show that if $F \in {\rm E}_\nu(n,\R)$ is such that $F(\vec{0}) = \vec{0}$, then $F \in {\rm O}_\nu(n,\R)$ (using a polarization formula for $\pair{\cdot,\cdot}_\nu$ and the result of Problem \ref{prob:auto_lin} ahead);
\item apply (i) for $\Lambda = F-F(\vec{0})$, where $F \in {\rm E}_\nu(n,\R)$ is now any pseudo-Euclidean isometry, and conclude that $F = T_{F(\vec{0})}\circ \Lambda$, where $T_{F(\vec{0})}$ denotes translation by $F(\vec{0})$;
\item check that $T_{\vec{a}_1}\circ \Lambda_1= T_{\vec{a}_2}\circ \Lambda_2$ implies $\vec{a}_1 = \vec{a}_2$ and $\Lambda_1 = \Lambda_2$, for all $\vec{a}_1,\vec{a}_2 \in \R^n_\nu$ and $\Lambda_1,\Lambda_2 \in {\rm O}_\nu(n,\R)$, by simply evaluating both sides of the assumed equality at $\vec{0}$.
\end{enumerate}

See Problem \ref{prob:semi_direct} in the end of the chapter for another point of view about this.

The pseudo-Euclidean space has a natural decomposition as $\R^n_\nu = \R^{n-\nu}\oplus \R^\nu_\nu$, as we have explored before in the proof of Theorem \ref{teo:light_LD}. This allows us to understand the structure of ${\rm O}_\nu(n,\R)$, by writing any $\Lambda$ in block-form as

\[\mbox{\Large
    \renewcommand{\arraystretch}{1.2}
    $\Lambda =
    \left(
      \begin{array}{c|c}
        \Lambda_S & B \\
        \hline
        C & \Lambda_T
      \end{array}
    \right)$},\] where $\Lambda_S \in {\rm Mat}(n-\nu, \R)$ e $\Lambda_T \in
{\rm Mat}(\nu,
\R)$ are to be called the \emph{spatial and temporal parts of $\Lambda$}\index{Spatial and temporal parts of a pseudo-orthogonal map}. Since $\Lambda$ is an isomorphism and preserves causal types, we have that $\Lambda_S$ e
$\Lambda_T$ are also non-singular. The blocks $\Lambda_S$ and $\Lambda_T$ are intimately related:

\begin{teo}\label{teo:det_lorentz}
 $\det \Lambda_S =  \det \Lambda_T \det \Lambda$.
\end{teo}

\begin{dem}
  Let ${\rm can} = (\vec{e}_i)_{i=1}^n$ be the usual basis for $\R^n_\nu$,
  and also consider the orthonormal basis of $\R^n_\nu$ formed by the columns of $\Lambda$, namely, $\mathcal{B} = \big( \Lambda\vec{e}_1,\ldots,
  \Lambda\vec{e}_n\big)$. Write $\Lambda$ explicitly as $\Lambda = (\lambda_{ij})_{1 \leq i,j \leq n}$. Let's ``delete'' the block $B$, defining a linear map
  $T\colon \R^n_\nu \to \R^n_\nu$ by \[ T(\Lambda\vec{e}_j) =
\begin{cases}
  \Lambda \vec{e}_j,&\mbox{if } 1 \leq j \leq n-\nu \quad\mbox{e} \\ \sum_{i=n-\nu+1}^n \lambda_{ij}\vec{e}_i,&\mbox{if } n-\nu < j \leq n.
\end{cases}
\]
We immediately have $[\Lambda]_{{\rm can},{\mathcal B}} = {\rm Id}_n$
and
\[\mbox{\Large
    \renewcommand{\arraystretch}{1.2}
    $[T]_{\mathcal{B},{\rm can}} =
    \left(
      \begin{array}{c|c}
        \Lambda_S & 0 \\
        \hline
        C & \Lambda_T
      \end{array}
    \right)$}.\]Compute now the matrix $[T]_{\mathcal
  B}$. The expression $T(\Lambda\vec{e}_j)= \Lambda\vec{e}_j$, which holds for the indices $1 \leq
j \leq
n-\nu$, tells us that the upper left and lower left blocks of $[T]_{\mathcal B}$ are, respectively, ${\rm Id}_{n-\nu}$ and
$0$. To compute the determinant of $[T]_{\mathcal
  B}$ by blocks, we need the last $\nu$ components of $T(\Lambda
\vec{e}_j)$ in the base ${\mathcal B}$, for $n-\nu < j \leq
n$.  Using the shorthand $\epsilon_k \doteq
\epsilon_{\vec{e}_k}$, we have:
\begin{align*}
  T(\Lambda \vec{e}_j) &= \sum_{i=n-\nu+1}^n \lambda_{ij}\vec{e}_i =\sum_{i=n-\nu+1}^n \lambda_{ij}\sum_{k=1}^n \epsilon_k \pair{\vec{e}_i, \Lambda\vec{e}_k}_\nu\Lambda\vec{e}_k \\ &= \sum_{i=n-\nu+1}^n \sum_{k=1}^n \sum_{\ell=1}^n \epsilon_k \lambda_{ij}\lambda_{\ell k}\pair{\vec{e}_i,\vec{e}_\ell}_\nu \Lambda\vec{e}_k \\ &= \sum_{k=1}^n \left(\sum_{i=n-\nu+1}^n \sum_{\ell=1}^n \epsilon_k \lambda_{ij}\lambda_{\ell k}\eta_{i\ell}^{\nu}  \right)\Lambda\vec{e}_k.
\end{align*}The desired last $\nu$ components correspond to $n-\nu < k \leq n$, and in these conditions, we have that the entries of the lower right block of $[T]_{\mathcal B}$ are given by \[
  \sum_{i=n-\nu+1}^n \sum_{\ell=1}^n - \lambda_{ij}\lambda_{\ell
    k}(-\delta_{i\ell}) =
  \sum_{i=n-\nu+1}^n\lambda_{ij}\lambda_{ik},  \]which we may recognize as the definition of the matrix product between $\Lambda_T^\top$ and
$\Lambda_T$. We obtain:
\[
\mbox{\Large $[T]_{\mathcal B}$} = 
\left(
    \begin{array}{c|c}
      \mbox{\Large ${\rm Id}_{n-\nu}$} &  \mbox{\Large $\ast$} \\[.7ex] \hline \\[-1.5ex]
      \mbox{\Large $0$} & \mbox{\Large $\Lambda_T^\top\Lambda_T$}\\
    \end{array}
  \right).\]In particular, it follows that $\det T=
(\det \Lambda_T)^2$. Moreover: \[ [T\Lambda]_{\mathcal{B}} =
  [T]_{\mathcal{B},{\rm can}}[\Lambda]_{{\rm can},\mathcal{B}}= 
\left(
    \begin{array}{c|c}
      \mbox{\Large $\Lambda_S$} &  \mbox{\Large $0$} \\[.7ex] \hline \\[-1.5ex]
      \mbox{\Large $C$} & \mbox{\Large $\Lambda_T$}\\
    \end{array}
  \right).\]Thus \[ (\det \Lambda_T)^2 \det \Lambda = \det T \det \Lambda
  = \det(T\Lambda) = \det \Lambda_T\det \Lambda_S,  \]and finally $\det \Lambda_S = \det \Lambda_T\det \Lambda$, as wanted.
\end{dem}

With this result in our hands, we may label the elements in ${\rm O}_\nu(n,\R)$ by the signs of the determinants of its spatial and temporal parts. This gives us a partition of ${\rm O}_\nu(n,\R)$:
\begin{align*}
  {\rm O}_\nu^{+\uparrow}(n,\R) &\doteq \{ \Lambda \in {\rm O}_\nu(n,\R) \mid \det \Lambda_S > 0 \,\,\,\mbox{e  }\det \Lambda_T > 0  \} \\  {\rm O}_\nu^{+\downarrow}(n,\R) &\doteq \{ \Lambda \in {\rm O}_\nu(n,\R) \mid \det \Lambda_S > 0 \mbox{  e  }\det \Lambda_T < 0  \}  \\  {\rm O}_\nu^{-\uparrow}(n,\R) &\doteq \{ \Lambda \in {\rm O}_\nu(n,\R) \mid \det \Lambda_S < 0 \,\,\,\mbox{e  }\det \Lambda_T > 0  \}  \\  {\rm O}_\nu^{-\downarrow}(n,\R) &\doteq \{ \Lambda \in {\rm O}_\nu(n,\R) \mid \det \Lambda_S < 0 \,\,\,\mbox{e  }\det \Lambda_T < 0  \} 
\end{align*}

Then we may say that the elements of ${\rm O}_\nu^{+ \bullet}(n,\R)$
\emph{preserve the orientation of space}, while the elements of ${\rm O}_\nu^{\bullet \uparrow}(n,\R)$ \emph{preserve the orientation of time} (i.e., they are \emph{or\-tho\-chro\-nous}\index{Orthochronous map}). We know that $\det\Lambda>0$ means that $\Lambda$ preserves the algebraic orientation of the vector space $\R^n_\nu$, but on the other hand, if $\Lambda\in{\rm O}_\nu(n,\R)$ and $\det\Lambda_S>0$ then $\Lambda$ preserves the spatial orientation of the spacelike subspaces\footnote{Now read this sentence again. Slowly.} of $\R^n_\nu$. Using convenient diagonal matrices with only $1$'s and $-1$'s, we conclude the:

\begin{cor}
  ${\rm O}_\nu^{+\downarrow}(n,\R)$, ${\rm O}_\nu^{-\downarrow}(n,\R)$ and ${\rm O}_\nu^{-\uparrow}(n,\R)$ are cosets of ${\rm O}_\nu^{+\uparrow}(n,\R)$.
\end{cor}

This means that we may focus our attention to the identity component ${\rm O}_\nu^{+\uparrow}(n,\R)$. In low dimensions, we have the following classifications:

\begin{teo}\label{teo:o12}
  \[ {\rm O}_1^{+\uparrow}(2,\R) = \left\{
    \begin{pmatrix}
      \cosh\varphi & \sinh\varphi \\ \sinh\varphi & \cosh\varphi
    \end{pmatrix}\in {\rm Mat}(2,\R) \mid \varphi \in \R
\right\}.  \]
\end{teo}

\begin{dem}
  Any $\Lambda = (\lambda_{ij})_{i,j=1}^2 \in {\rm O}_1^{+\uparrow}(2,\R)$ satisfies \[
  \begin{cases}
    \lambda_{11}^2-\lambda_{21}^2 &= 1 \\ \lambda_{12}^2 - \lambda_{22}^2 &= -1, \quad\mbox{and} \\ \lambda_{11}\lambda_{12} - \lambda_{21}\lambda_{22} &= 0
  \end{cases}
\]with $\lambda_{11},\lambda_{22} \geq 1$. So we get unique $t,s \in \R_{\geq 0}$ with $\lambda_{11} = \cosh t$ and $\lambda_{22} = \cosh s$. The above equations imply that $|\lambda_{21}| = \sinh t$ and $|\lambda_{12}| = \sinh s$. The additional condition $\det \Lambda =1$ gives $\lambda_{12}\lambda_{21} = \cosh t\cosh s-1 \geq 0$, so $\lambda_{12}$ and $\lambda_{21}$ have the same sign. No matter which sign, the the third equation above now says that \[ 0 = \cosh t\sinh s - \sinh t \cosh s = \sinh(s-t) \implies s=t. \]Then $\Lambda$ is one of the following matrices, for $t>0$: \[
\begin{pmatrix}
  \cosh t & \sinh t \\ \sinh t & \cosh t
\end{pmatrix}\qquad\mbox{or}\qquad
\begin{pmatrix}
  \cosh t & -\sinh t \\ -\sinh t & \cosh t
\end{pmatrix}.\]
\end{dem}

A somewhat similar strategy also gives us the classification in dimension $3$:

\begin{teo}
Any $\Lambda \in {\rm O}_1^{+\uparrow}(3,\R)$ is conjugate to one of the following matrices: \[
\begin{pmatrix}
  1 & 0 & 0 \\ 0 & \cosh \varphi & \sinh \varphi \\ 0 & \sinh \varphi & \cosh \varphi
\end{pmatrix}, \quad
\begin{pmatrix}
  \cos \theta & -\sin \theta & 0 \\ \sin \theta & \cos \theta & 0 \\ 0 & 0 & 1
\end{pmatrix}, \quad\mbox{or}\quad \begin{pmatrix} 
                     1 &-\theta &\theta \\[.5ex]
                     \theta & 1 - \nicefrac{\theta^2}{2} &
                     \nicefrac{\theta^2}{2} \\[.5ex]
                     \theta & \nicefrac{-\theta^2}{2} & 1+
                     \nicefrac{\theta^2}{2} \end{pmatrix},
\] where $\varphi,\theta \in \R$. The transformation $\Lambda$ is called hyperbolic, elliptic or parabolic, depending on its conjugacy class. 
\end{teo}

\begin{obs}\mbox{}
  \begin{itemize}
  \item  One can prove that any $\Lambda \in {\rm O}_1^{+\uparrow}(3,\R)$ has at least one unit eigenvector, say $\vec{v}$. The causal character of $\vec{v}$ decides what is the class of $\Lambda$. Namely, $\Lambda$ is hyperbolic if $\vec{v}$ is spacelike (so $\Lambda$ acts as an hyperbolic rotation in the timelike plane $\vec{v}^\perp$), elliptic if $\vec{v}$ is timelike (so $\Lambda$ acts as a Euclidean rotation in the spacelike plane $\vec{v}^\perp$), and parabolic if $\vec{v}$ is lightlike (so $\Lambda$ has that shear-like action in the null line defined by $\vec{v}$).

  \item  This terminology is also useful in establishing the classification of helices in $\LM^3$ (\emph{Lancret's theorem}\index{Lancret's theorem}), according to the causal type of the helix's axis. 
  \end{itemize}
\end{obs}

\subsection{Relation with Special Relativity}

Here we will motivate the names ``spacelike'', ``timelike'' and ``lightlike'', and try to give some relation between what we have done so far and the mathematics used in Special Relativity. We focus on Lorentz-Minkowski space $\LM^4$, whose points are, in this setting, called \emph{events}. Fixing the inertial frame given by the standard basis of $\LM^4$, we write the coordinates in $\LM^4$ as $(x,y,z,t)$. Assume that a particle with positive mass moves in spacetime from event $\vec{p}$ to event $\vec{q}$, through some time interval $\Delta t \neq 0$, and let \[\vec{v} = \vec{q} - \vec{p} = (\Delta x, \Delta y, \Delta z, \Delta t)\] be the spacetime displacement vector. The fact that the particle may not move at a speed greater than the speed of light $c$ may be written as \[ \left(\frac{\Delta x}{\Delta t}\right)^2 + \left(\frac{\Delta y}{\Delta t}\right)^2 +\left(\frac{\Delta z}{\Delta t}\right)^2  < c^2.    \]So, if we let $\widetilde{\vec{v}} = (\Delta x/\Delta t, \Delta y/\Delta t, \Delta z/\Delta t)$ be the velocity vector of the worldline of the particle, in $\R^3 \cong \R^3 \oplus \{0\} \subseteq \LM^4$, the above means that $\|\widetilde{\vec{v}}\|_E<c$. We henceforth set the so called \emph{geometric units}\index{Geometric units}, where $c=1$. With this in mind, computing
\begin{align*}
  \pair{\vec{v},\vec{v}}_L &= (\Delta x)^2+(\Delta y)^2 + (\Delta z)^2 - (\Delta t)^2 \\ &= (\Delta t)^2 \left( \left(\frac{\Delta x}{\Delta t}\right)^2 + \left(\frac{\Delta y}{\Delta t}\right)^2 +\left(\frac{\Delta z}{\Delta t}\right)^2-1\right) \\ &= (\Delta t)^2(\|\widetilde{\vec{v}}\|_E^2-1) \\ &= (\Delta t)^2(\|\widetilde{\vec{v}}\|_E+1)(\|\widetilde{\vec{v}}\|_E-1)
\end{align*}
we see that:
\begin{enumerate}[(1)]
\item if $\vec{v}$ is timelike, then $\|\widetilde{\vec{v}}\|_E<1$, and so the event $\vec{p}$ may influence event $\vec{q}$ if $\Delta t > 0$, and the other way around if $\Delta t < 0$, e.g., via the propagation of a material wave.
\item if $\vec{v}$ is lightlike and $\Delta t \neq 0$, then $\|\widetilde{\vec{v}}\|_E=1$ and so the influence between the events can only be given via the propagation of some eletromagnectic wave, or by the emission of some light signal sent by one of the events and reaching the other.
\item if $\vec{v}$ is spacelike with $\Delta t \neq 0$, there is no influence relation between the events, since $\|\widetilde{\vec{v}}\|_E>1$ means that the speed necessary for a particle starting at one event to reach the spatial location of the other must be greater than the speed of light, which is impossible: not even a photon or neutrino is fast enough to experience both events. Both of them are not inside, or even in the boundary, of the other's lightcone.
\begin{figure}[H]
  \centering
  \includegraphics[height=5cm]{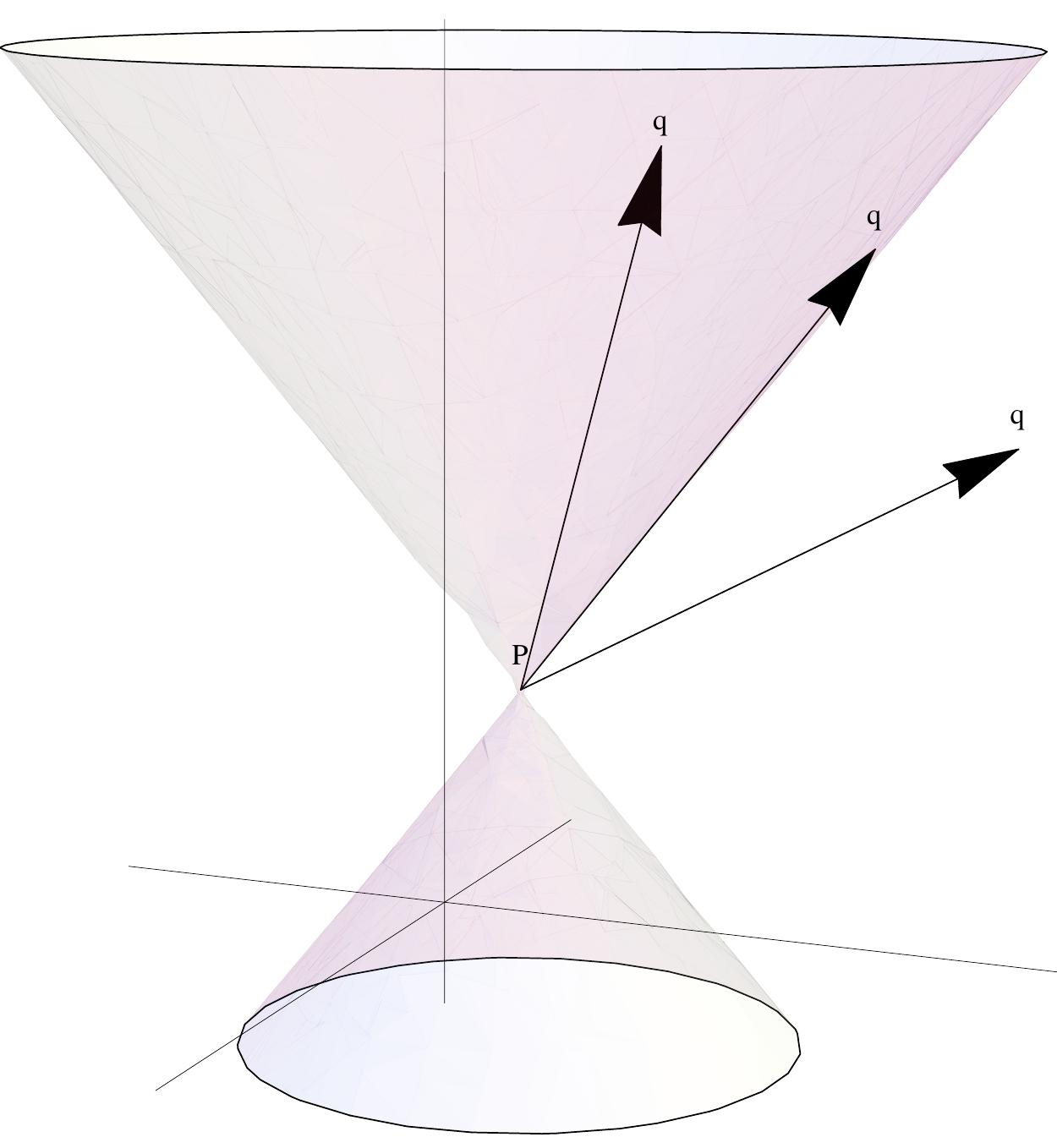}
  \caption[Physical interpretation for causal characters.]{Physical interpretation for causal characters.}
\end{figure}
\end{enumerate}

Let's try and make more precise this notion of causal influence. For this, we need to know what does it mean for a vector to point to the future (or past):

\begin{defn}
Let $\vec{e}_n = (0,\ldots,0,1) \in \LM^n$. A timelike or lightlike vector $\vec{v} \in \LM^n$ is \emph{future-directed} (resp. \emph{past-directed}) if $\pair{\vec{v},\vec{e}_n}_L< 0$ (resp. $\pair{\vec{v},\vec{e}_n}_L> 0$).  
\end{defn}

\begin{defn}[$\ll$ and $\preccurlyeq$] Given $\vec{p} \in \LM^n$, we define the \emph{timecone and lightcone centered at $\vec{p}$} by
  \[  C_T(\vec{p}) = \{ \vec{q} \in \LM^n \mid \vec{q}-\vec{p}\mbox{ is timelike}\}\quad\mbox{and}\quad C_L(\vec{p}) = \{ \vec{q} \in \LM^n \mid \vec{q}-\vec{p}\mbox{ is lightike}\}.\]Naturally, using the previous definition we may divide those in \emph{future cones} $C_T^+(\vec{p})$ and $C_L^+(\vec{p})$, and \emph{past cones} $C_T^-(\vec{p})$ and $C_L^-(\vec{p})$. We'll say that $\vec{p}$ \emph{chronologically preceds}\index{Chronological and causal precedences} $\vec{q}$ (resp. \emph{causally preceds} $\vec{q}$) if $\vec{q} \in C_T^+(\vec{p})$ (resp. $\vec{q} \in C_T^+(\vec{p})\cup C_L^+(\vec{p})$). These relations will be denoted by $\vec{p} \ll \vec{q}$ and $\vec{p} \preccurlyeq \vec{q}$.
\end{defn}

Let's list some properties of these relations:

\begin{prop}
  Given $\vec{p}, \vec{u},\vec{v} \in \LM^n$, we have that:
  \begin{enumerate}[(i)]
  \item if $\vec{u},\vec{v} \in C_T^+(\vec{p})$, then $\pair{\vec{u}-\vec{p},\vec{v}-\vec{p}}_L<0$;
  \item if $\vec{u},\vec{v} \in C_T(\vec{p})$ and $\pair{\vec{u}-\vec{p},\vec{v}-\vec{p}}_L < 0$, then both $\vec{u}$ and $\vec{v}$ are in $C_T^+(\vec{p})$ or $C_T^-(\vec{p})$;
  \item $\ll$ and $\preccurlyeq$ are transitive.
  \end{enumerate}
\end{prop}

Geometrically, they're easy to understand, but their proofs rely on technicalities with hyperbolic trigonometric functions. You are welcome to try and prove them, but you can check the proofs on~\cite{Na} or~\cite{TL}, and more general results in the contexts of spacetimes in General Relativity may be found on~\cite{BEE},~\cite{Haw} and~\cite{Pen}.

Now, we have previously mentioned that the ``fake norm'' $\|\cdot\|_L$ has poor properties, which is mainly due to the fact that it is not induced by a positive-definite inner product. In this context, here's probably the best we can get:

\begin{prop}[Backwards Cauchy-Schwarz]\label{prop:backwards_cs}
  Let $\vec{u},\vec{v} \in \LM^n$ be timelike vectors. Then $|\pair{\vec{u},\vec{v}}_L| \geq \|\vec{u}\|_L\|\vec{v}\|_L$. Furthermore, equality holds if and only if $\vec{u}$ and $\vec{v}$ are proportional.
\end{prop}

\begin{dem}
  Write $\LM^n = \R\vec{u}\oplus \vec{u}^\perp$ and write $\vec{v} = \lambda\vec{u}+\vec{u}_0$, with $\lambda \in \R$ and $\vec{u}_0$ spacelike and Lorentz-orthogonal to $\vec{u}$. On one hand, we have $\pair{\vec{v},\vec{v}}_L = \lambda^2\pair{\vec{u},\vec{u}}_L+\pair{\vec{u}_0,\vec{u}_0}_L$. On the other, we compute:
  \begin{align*}
    \pair{\vec{u},\vec{v}}_L^2 &= \pair{\vec{u}, \lambda\vec{u}+\vec{u}_0}_L^2 \\ &= \lambda^2 \pair{\vec{u},\vec{u}}_L^2 \\ &= \big(\pair{\vec{v},\vec{v}}_L - \pair{\vec{u}_0,\vec{u}_0}_L\big)\pair{\vec{u},\vec{u}}_L \\ &\geq \pair{\vec{v},\vec{v}}_L\pair{\vec{u},\vec{u}}_L > 0,
  \end{align*}using that $\vec{u}_0$ is spacelike and $\vec{u}$ is timelike. The result follow by taking roots. Note that the equality holds if and only if $\vec{u}_0 = \vec{0}$, which is equivalent to $\vec{u}$ and $\vec{v}$ being proportional.
\end{dem}

With this we may define the hyperbolic angle between timelike vectors, both future-directed or past-directed, in the same fashion one defines the angle between vectors in a vector space with a positive-definite inner product. Since the image of $\cosh$ is the interval $[1,+\infty[$, there is a unique $\varphi \geq 0$ such that $\pair{\vec{u},\vec{v}}_L = -\|\vec{u}\|_L\|\vec{v}\|_L\cosh\varphi$. Another consequence is the:

\begin{cor}[Backwards triangle inequality]
  Let $\vec{u},\vec{v} \in \LM^n$ timelike vectors, both future-directed or past-directed. Then $\|\vec{u}+\vec{v}\|_L \geq \|\vec{u}\|_L+\|\vec{v}\|_L$.
\end{cor}

As a general strategy in Mathematics, once we have defined something (here, $\ll$ and $\preccurlyeq$), it is natural to turn our attention to the mappings related to what we have defined. So we write the:

\begin{defn}
  A map $F\colon \LM^n \to \LM^n$ is called a \emph{causal automorphism}\index{Causal!automorphism} if it is bijective and both $F$ and $F^{-1}$ preserve $\preccurlyeq$, that is: \[   \vec{x}\preccurlyeq \vec{y} \iff F(\vec{x}) \preccurlyeq F(\vec{y})\qquad\mbox{and}\qquad  \vec{x}\preccurlyeq \vec{y} \iff F^{-1}(\vec{x}) \preccurlyeq F^{-1}(\vec{y}). \]
\end{defn}

\begin{obs}
  It can be shown that preserving $\preccurlyeq$ is the same as preserving $\ll$.
\end{obs}

Obvious examples of causal automorphisms are positive homotheties, translations, and orthochronous Lorentz transformations. Amazingly, that's all of them:

\begin{teo}[Alexandrov-Zeeman]
  Let $n \geq 3$ and $F\colon \LM^n \to \LM^n$ be a causal automorphism. Then there is a positive constant $c>0$, an orthochronous Lorentz transformation $\Lambda$, and a vector $\vec{a} \in \LM^n$ such that \[ F(\vec{x}) = c\Lambda(\vec{x})+\vec{a},\qquad \mbox{ for all }\vec{x} \in \LM^n.  \]Moreover, this decomposition is unique.
\end{teo}

The proof of this theorem is actually difficult (except maybe for the uniqueness part\footnote{{\bf Proof:} assume $c_1\Lambda_1(\vec{x})+\vec{a}_1 = c_2\Lambda_2(\vec{x})+\vec{a}_2$ for all $\vec{x} \in \LM^n$, according to the statement of the theorem. Evaluate at $\vec{0}$ to get $\vec{a}_1=\vec{a}_2$. Cancel the translation to get $c_1\Lambda_1(\vec{x}) = c_2\Lambda_2(\vec{x})$ for all $\vec{x} \in \LM^n$. Take the scalar square of both sides to get $c_1^2=c_2^2$. Since $c_1,c_2>0$, it follows that $c_1=c_2$. We conclude that $\Lambda_1=\Lambda_2$.$\hfill \qed$}), employing a mix of results from the linear algebra we have seen so far, Darboux's fundamental theorem of geometry (regarding certain doubly-ruled surfaces), and lifting properties of some maps. You can see more details in~\cite{Na}, for example. The result is false for $n=2$, in view of some deeper results about the conformal structure of $\LM^2$ -- a counter-example is discussed in~\cite{TL}. Furthermore, this theorem has also a topological flavor: the usual topology in $\LM^n$ does not properly capture the causal features of this spacetime, in contrast with the so called \emph{Zeeman topology}\index{Zeeman topology}, whose homeomorphisms are precisely the causal automorphisms here discussed. For more about these topologies, you may consult~\cite{Na2}.

\subsection{Cross product}

With a new scalar product $\pair{\cdot,\cdot}_\nu$, comes together a new notion of \emph{cross product}:

\begin{defn}
  The \emph{index $\nu$ cross product}\index{Cross product} of $\vec{v}_1,\ldots, \vec{v}_{n-1} \in \R^n_\nu$ is the unique vector $\vec{v} \in \R^n_\nu$ such that $\pair{\vec{v},\vec{x}}_\nu = \det(\vec{x},\vec{v}_1,\ldots,\vec{v}_{n-1})$, for all $\vec{x} \in \R^n_\nu$. The existence and uniqueness of such $\vec{v}$ is ensured by the non-degenerability of $\pair{\cdot,\cdot}_\nu$. We then denote $\vec{v}$ by $\vec{v}_1 \times \cdots \times \vec{v}_{n-1}$, the index $\nu$ being understood.
\end{defn}

\begin{obs}
  Just like we denote the scalar products of $\R^3$ and $\LM^3$ by $\pair{\cdot,\cdot}_E$ and $\pair{\cdot,\cdot}_L$, we'll follows this convention for cross products, using $\times_E$ and $\times_L$, respectively.
\end{obs}

Just from the definition, we the cross product inherits some immediate properties from $\det$, registered in the:

\begin{prop}
  The index $\nu$ cross product in $\R^n_\nu$ is $(n-1)$-multilinear, totally skew-symmetric, and orthogonal to each of its arguments. If $n=3$, it additionaly satisfies the identity $\pair{\vec{v}_1\times \vec{v}_2,\vec{v}_3}_\nu = \pair{\vec{v}_1,\vec{v}_2\times \vec{v}_3}_\nu$, for all $\vec{v}_1,\vec{v}_2,\vec{v}_3 \in \R^3_\nu$ (comma commutes with $\times$).
\end{prop}

As important as these properties are, they still do not tell us how to explicitly compute cross products. Just like when you first learned about cross products in $\R^3$, we'll keep using formal determinants with a convenient Laplace expansion along the first row:

\begin{prop}
  Let $\mathcal{B} = (\vec{u}_i)_{i=1}^n$ be a positive orthonormal basis for $\R^n_\nu$ and let be given vectors $\vec{v}_j = \sum_{i=1}^n v_{ij}\vec{u}_i \in \R^n_\nu$, for $1 \leq j \leq n-1$. Using the shorthand $\epsilon_i \doteq \epsilon_{\vec{u}_i}$ for the indicators of the elements in $\mathcal{B}$, we have:
\[ \vec{v}_1 \times \cdots \times \vec{v}_{n-1} =
\begin{vmatrix}
  \epsilon_1 \vec{u}_1 & \cdots & \epsilon_{n}\vec{u}_{n} \\ v_{11} & \cdots & v_{n1} \\ \vdots & \ddots & \vdots \\ v_{1,n-1} & \cdots & v_{n,n-1} 
\end{vmatrix}.
\]
\end{prop}

\begin{prop}\label{prop:prods_det}
  Let $\vec{u}_1,\ldots,\vec{u}_{n-1},\vec{v}_1,\ldots,\vec{v}_{n-1} \in \R^n_\nu$. Then we have \[ \pair{\vec{u}_1\times\cdots\times\vec{u}_{n-1},\vec{v}_1\times\cdots\times\vec{v}_{n-1}}_\nu = (-1)^\nu \det\big( (\pair{\vec{u}_i,\vec{v}_j}_\nu)_{1 \leq i,j \leq n-1} \big).  \]
\end{prop}

\begin{dem}
  If $(\vec{u}_i)_{i=1}^{n-1}$ or $(\vec{v}_j)_{j=1}^{n-1}$ is linearly dependent, there's nothing to do. Assume then that both are linearly independent. Since both sides of the proposed equality are linear in each of the $2n-2$ variables, and both the cross product and the determinant are totally skew-symmetric, we may assume without loss of generality that $\vec{u}_k = \vec{e}_{i_k}$ and $\vec{v}_\ell = \vec{e}_{j_\ell}$, where $(\vec{e}_i)_{i=1}^n$ is the standard basis for $\R^n_\nu$ and \[ 1 \leq i_1 < \cdots < i_{n-1}\leq n \quad\mbox{e}\quad 1 \leq j_1 < \cdots < j_{n-1}\leq n.  \]
We will proceed with the analysis in cases, in terms of the indices $i^*$ and $j^*$ being omitted in each of the two $(n-1)$-uples of indices considered.
\begin{itemize}
  \item If $i^* \neq j^*$, both sides vanish. To wit, the left hand side equals $\pair{\vec{e}_{i^*},\vec{e}_{j^*}}_\nu = 0$, and the determinant on the right hand side has the $i^*$-th row and the $j^*$-th column only with zeros.\\[-.8cm]
\item If $1 \leq  i^* = j^* \leq n-\nu$, the left hand side equals $\pair{\vec{e}_{i^*},\vec{e}_{i^*}}_\nu = 1$, and the right hand side equals $(-1)^\nu \det {\rm Id}_{n-1,\nu} = (-1)^\nu(-1)^\nu = 1$. \\[-.8cm]
\item If $n-\nu < i^* = j^* \leq n$, the left hand side equals $\pair{\vec{e}_{i^*},\vec{e}_{i^*}}_\nu = -1$, and the right hand side equals $(-1)^\nu \det {\rm Id}_{n-1,\nu-1} = (-1)^\nu (-1)^{\nu-1} = -1$.
\end{itemize}
\end{dem}
\vspace{-.8cm}
\begin{cor}[Lagrange's Identities]\label{cor:lagrange}\index{Lagrange's Identities}
  Let $\vec{u},\vec{v} \in \R^3_\nu$. Then:
  \begin{align*}
    \|\vec{u} \times_E \vec{v}\|^2_E &=\|\vec{u}\|^2_E\|\vec{v}\|^2_E -
                                       \pair{\vec{u},\vec{v}}^2_E,\\
    \pair{\vec{u}\times_L \vec{v}, \vec{u} \times_L \vec{v}}_L &=
                                                                 -\pair{\vec{u},\vec{u}}_L
                                                                 \pair{\vec{v},\vec{v}}_L
                                                                 +
                                                                 \pair{\vec{u},\vec{v}}^2_L.
  \end{align*}
\end{cor}

The orientation of the bases chosen for $\R^n_\nu$ will be very important for defining convenient frames along lightlike and semi-lightlike curves in the next chapter. So we might as well discuss this now in a bit greater generality. We follow the convention that the standard basis for $\R^n_\nu$ is, of course, positive.

If $\vec{v}_1,\ldots,\vec{v}_{n-1}\in\R^n_\nu$ are linearly independent, do not span a lightlike hyperplane, and we denote $\vec{v}=\vec{v}_1\times\cdots\times\vec{v}_{n-1}$,
then ${\mathcal B}=\big(\vec{v}_1,\ldots,\vec{v}_{n-1},\vec{v}\big)$ is a basis for $\R^n_\nu$, and it would natural to ask ourselves when such basis is positive or negative. The answer is in the determinant of the matrix having these vectors in rows or columns. We have\vspace{-.1cm}
\[\det(\vec{v}_1,\ldots,\vec{v}_{n-1},\vec{v})=(-1)^{n-1}\det(\vec{v},\vec{v}_1,\ldots,\vec{v}_{n-1})=(-1)^{n-1}\pair{\vec{v},\vec{v}}_\nu\]\vspace{-.1cm}
and, hence, positiveness of the basis $\mathcal{B}$ depends not only on the parity of $n$, but also on the causal character of $\vec{v}$. Explicitly: if $\vec{v}$ is spacelike, $\mathcal{B}$
is positive if $n$ is odd, and negative if $n$ is even; if $\vec{v}$ is timelike, $\mathcal{B}$ is positive if $n$ is even, and negative if $n$ is odd. 

In particular, for $n=3$ we may represent all the possible cross products between the elements in the standard basis of $\R^3_\nu$ by the following diagrams:

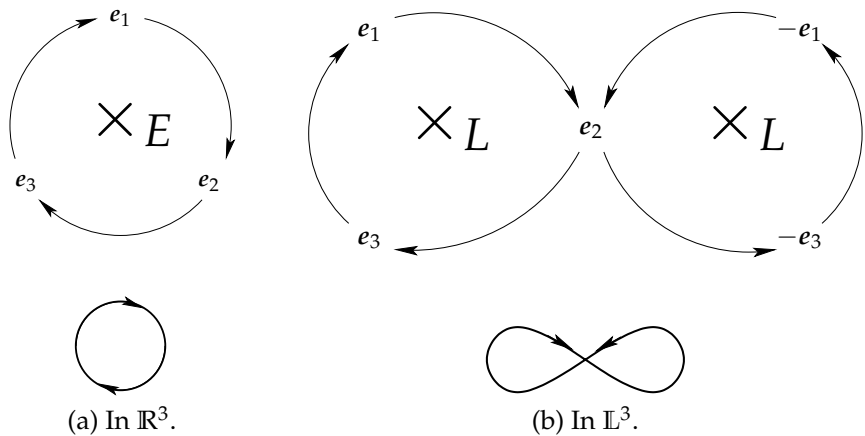
\begin{figure}[H]
  \centering \subfloat[In
  $\R^3$.]{\begin{picture}(0,0)%
\includegraphics{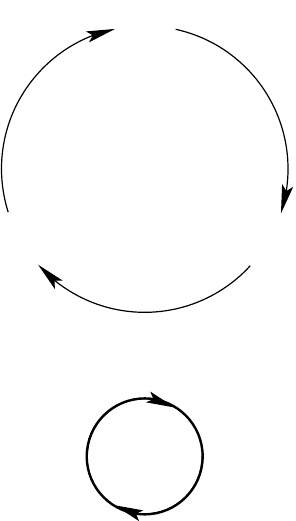}%
\end{picture}%
\setlength{\unitlength}{3025sp}%
\begingroup\makeatletter\ifx\SetFigFont\undefined%
\gdef\SetFigFont#1#2#3#4#5{%
  \reset@font\fontsize{#1}{#2pt}%
  \fontfamily{#3}\fontseries{#4}\fontshape{#5}%
  \selectfont}%
\fi\endgroup%
\begin{picture}(1837,3242)(444,-2715)
\put(1081,-601){\makebox(0,0)[lb]{\smash{{\SetFigFont{9}{10.8}{\rmdefault}{\mddefault}{\updefault}{\color[rgb]{0,0,0}{\Huge $\times_{\mbox{{\Large $E$}}}$}}%
}}}}
\put(496,-1006){\makebox(0,0)[lb]{\smash{{\SetFigFont{9}{10.8}{\rmdefault}{\mddefault}{\updefault}{\color[rgb]{0,0,0}$\vec{e}_3$}%
}}}}
\put(1261,344){\makebox(0,0)[lb]{\smash{{\SetFigFont{9}{10.8}{\rmdefault}{\mddefault}{\updefault}{\color[rgb]{0,0,0}$\vec{e}_1$}%
}}}}
\put(1981,-1006){\makebox(0,0)[lb]{\smash{{\SetFigFont{9}{10.8}{\rmdefault}{\mddefault}{\updefault}{\color[rgb]{0,0,0}$\vec{e}_2$}%
}}}}
\end{picture}}\label{fig:cross_r3} \qquad
  \subfloat[In
  $\LM^3$.]{\begin{picture}(0,0)%
\includegraphics{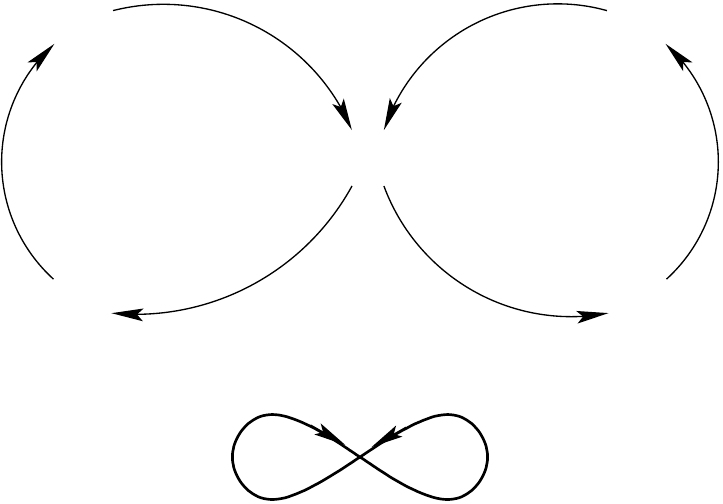}%
\end{picture}%
\setlength{\unitlength}{3356sp}%
\begingroup\makeatletter\ifx\SetFigFont\undefined%
\gdef\SetFigFont#1#2#3#4#5{%
  \reset@font\fontsize{#1}{#2pt}%
  \fontfamily{#3}\fontseries{#4}\fontshape{#5}%
  \selectfont}%
\fi\endgroup%
\begin{picture}(4063,2830)(444,-2438)
\put(811,209){\makebox(0,0)[lb]{\smash{{\SetFigFont{10}{12.0}{\rmdefault}{\mddefault}{\updefault}{\color[rgb]{0,0,0}$\vec{e}_1$}%
}}}}
\put(811,-1321){\makebox(0,0)[lb]{\smash{{\SetFigFont{10}{12.0}{\rmdefault}{\mddefault}{\updefault}{\color[rgb]{0,0,0}$\vec{e}_3$}%
}}}}
\put(1171,-556){\makebox(0,0)[lb]{\smash{{\SetFigFont{10}{12.0}{\rmdefault}{\mddefault}{\updefault}{\color[rgb]{0,0,0}{\Huge $\times_{\mbox{{\Large $L$}}}$}}%
}}}}
\put(2431,-511){\makebox(0,0)[lb]{\smash{{\SetFigFont{10}{12.0}{\rmdefault}{\mddefault}{\updefault}{\color[rgb]{0,0,0}$\vec{e}_2$}%
}}}}
\put(3871,209){\makebox(0,0)[lb]{\smash{{\SetFigFont{10}{12.0}{\rmdefault}{\mddefault}{\updefault}{\color[rgb]{0,0,0}$-\vec{e}_1$}%
}}}}
\put(3871,-1321){\makebox(0,0)[lb]{\smash{{\SetFigFont{10}{12.0}{\rmdefault}{\mddefault}{\updefault}{\color[rgb]{0,0,0}$-\vec{e}_3$}%
}}}}
\put(3331,-556){\makebox(0,0)[lb]{\smash{{\SetFigFont{10}{12.0}{\rmdefault}{\mddefault}{\updefault}{\color[rgb]{0,0,0}{\Huge $\times_{\mbox{{\Large $L$}}}$}}%
}}}}
\end{picture}}\label{fig:cross_l3}
  \caption{Understanding the cross products in $\R^3_\nu$.}
  \label{fig:cross}
\end{figure}

The cross products are obtained by following the arrows. For example, we have $\vec{e}_2\times_E\vec{e}_3=\vec{e}_1$ and
$\vec{e}_1\times_L\vec{e}_2=-\vec{e}_3$. In $\R^3$, following the arrows in the opposite direction, we obtain the results with the swapped sign (since $\times_E$ is skew), e.g.,
$\vec{e}_1\times_E\vec{e}_3=-\vec{e}_2$. In $\LM^3$ this does not work anymore due to the presence of causal characters: note that \linebreak[4]$\vec{e}_3\times_L\vec{e}_2=-\vec{e}_1\neq -(-\vec{e}_1) =
\vec{e}_1$. The cross products in $\LM^3$ which cannot be obtained directly from the above diagram may be obtained by using that $\times_L$ is skew, after finding the up-to-sign correct product on the diagram.

\begin{obs}
  These diagrams remain valid using any positive orthonormal basis of the space, provided that in $\LM^3$ the timelike vector (corresponding to $\vec{e}_3$) is the last one.
\end{obs}

We'll conclude the chapter stating two general facts from linear algebra, which will be necessary for giving adequate definitions for the Gaussian and mean curvatures of a surface later:

\begin{lem}
Let $B\colon \R^n_\nu\times \R^n_\nu \to Z$ be a bilinear map, where $Z$ is any vector space. If $(\vec{v}_i)_{i=1}^n$ and $(\vec{w}_i)_{i=1}^n$ are orthonormal bases for $\R^n_\nu$, then we have:
\begin{enumerate}[(i)]
\item $\sum_{i=1}^n \epsilon_{\vec{v}_i}B(\vec{v}_i,\vec{v}_i) = \sum_{i=1}^n \epsilon_{\vec{w}_i}B(\vec{w}_i,\vec{w}_i)$;
\item $\det\big((B(\vec{v}_i,\vec{v}_j))_{i,j=1}^n\big)=\det\big((B(\vec{w}_i,\vec{w}_j))_{i,j=1}^n\big)$, provided $Z=\R$.
\end{enumerate}
These quantities (which are then invariant under change of basis) are denoted by ${\rm tr}_{\pair{\cdot,\cdot}_\nu}B$ and $\det_{\pair{\cdot,\cdot}_\nu}B$.
\end{lem}

\newpage
\section*{Problems}\addcontentsline{toc}{subsection}{Problems for the first section}

\begin{problem}[Triangles of light]
  Show that in $\LM^n$ we cannot have lightlike vectors $\vec{u}_1,\vec{u}_2$ and $\vec{u}_3$ with $\vec{u}_1+\vec{u}_2+\vec{u}_3=\vec{0}$ and $\{\vec{u}_1,\vec{u}_2,\vec{u}_3\}$ linearly independent. Try to generalize.
\end{problem}

\begin{problem}
  Show that if $S\subseteq \LM^n$ is lightlike then $\dim(S\cap S^\perp)=1$, and conclude that if $S$ is a lightlike hyperplane, then $S^\perp\subseteq S$. Give an example of a subspace $S\subseteq \R^n_\nu$ with $\nu > 1$ and $\dim(S \cap S^\perp) \geq 2$.
\end{problem}

\begin{problem}[Sylvester's Law of Inertia]
  Show that every orthonormal basis for $\R^n_\nu$ must necessarily have $n-\nu$ spacelike vectors, $\nu$ timelike vectors, and no lightlike vectors.
  \begin{dica}
    There's a proof in~\cite{ON2}, which you should try to at least understand if you cannot come up with a solution on your own.
  \end{dica}
\end{problem}

\begin{problem}\label{prob:auto_lin}
  Show that if a map $\Lambda\colon \R^n_\nu \to \R^n_\nu$ preserves $\pair{\cdot,\cdot}_\nu$, then it is automatically linear (and hence in ${\rm O}_\nu(n,\R)$).
\end{problem}

\begin{problem}\label{prob:semi_direct}
  Consider the semi-direct product ${\rm O}_\nu(n,\R) \ltimes \R^n_\nu$ with operation $\ast$ given by \[ (A,\vec{v})\ast (B,\vec{w}) = (AB,A\vec{w}+\vec{v}).  \]Prove that this operation is indeed associative with identity element $({\rm Id}_n, \vec{0})$, compute $(A,\vec{v})^{-1}$ for any $(A,\vec{v}) \in {\rm O}_\nu(n,\R)\ltimes \R^n_\nu$, and show that $\Phi\colon {\rm O}_\nu(n,\R)\ltimes \R^n_\nu\to  {\rm E}_\nu(n,\R)$ given by $\Phi(A,\vec{v}) = T_{\vec{v}}\circ A$ is a group isomorphism.
\end{problem}

\begin{problem}
  Let $\Lambda \in {\rm O}_1(n,\R)$ be a Lorentz transformation.
  \begin{enumerate}[(a)]
  \item Show that a non-lightlike eigenvector must have $1$ or $-1$ as associated eigenvalue.
  \item Show that the product of the eigenvalues associated with two linearly independent lightlike vectors is $1$.
  \item If $W\subseteq \LM^n$ is an eigenspace of $\Lambda$ containing a non-lightlike vector, show that every other eigenspace of $\Lambda$ is orthogonal to $W$.
  \item If $W\subseteq \LM^n$ is a subspace, show that $W$ is $\Lambda$-stable (i.e., $\Lambda[W]\subseteq W$) if and only if $W^\perp$ is $\Lambda$-stable.
  \end{enumerate}
\end{problem}

\begin{problem}[Margulis Invariant]
  Let $F \in {\rm P}(3,\R)$ be a hyperbolic Poincar\'{e} transformation, given by $F(\vec{x}) = \Lambda\vec{x}+\vec{w}$, with $\Lambda \in {\rm O}_1^{+\uparrow}(3,\R)$ and $\vec{w} \in \LM^3$.
  \begin{enumerate}[(a)]
  \item Show that $\Lambda$ has three positive eigenvalues $1/\lambda < 1<\lambda$. The eigenspaces associated to $\lambda$ and $1/\lambda$ are automatically null lines.
  \item Let $\vec{v}_\lambda$ and $\vec{v}_{1/\lambda}$ be future-directed eigenvectors associated to $\lambda$ and $1/\lambda$, and $\vec{v}_1$ be a unit eigenvector associated to $1$ such that the base $\mathcal{B} = (\vec{v}_\lambda,\vec{v}_1,\vec{v}_{1/\lambda})$ is positive. Show that $F$ leaves invariant a unique (affine) line parallel to $\vec{v}_1$, and acts on such line by translation. That is, show that there are $\vec{p} \in \LM^3$ and $\alpha_F \in \R$ such that \[  F(\vec{p}+t\vec{v}_1) = \vec{p}+t\vec{v}_1 + \alpha_F\vec{v}_1, \]for all $t \in \R$. We say that $\alpha_F$ is the \emph{Margulis invariant}\index{Margulis invariant} of $F$.
  \item Show that $\alpha_F = \pair{\vec{w},\vec{v}_1}_L$ and use this to show that if $F_1,F_2 \in {\rm P}(3,\R)$ are hyperbolic and conjugate by an element of ${\rm O}_1(3,\R)$, then $\alpha_{F_1} = \alpha_{F_2}$ (thus justifying the name ``invariant'').
  \item Show that for every non-zero integer $n$, $\alpha_{F^n} = n\alpha_F$. Be careful with the case $n<0$, and pay close attention to the orientation of the eigenbasis associated to $\Lambda^{-1}$.
  \end{enumerate}
\end{problem}

\begin{problem}
  Let $\vec{x}\in \LM^3$ be a spacelike vector. Show that $T \doteq \vec{x}\times \_\colon \LM^3 \to \LM^3$ is diagonalizable, and the directions of the null lines given by the intersection of the timelike plane $\vec{x}^\perp$ with the lightcone of $\LM^3$ are eigenvectors of $T$.
\end{problem}

\begin{problem}[Lorentz $\gamma$ factor]
  Let $\vec{v} = (\Delta x_1,\ldots, \Delta x_{n-1},\Delta t) \in \LM^n$ be the displacement vector of a particle, moving between two events in spacetime. Show that the hyperbolic angle $\varphi$ between $\vec{v}$ and $\vec{e}_n$ is characterized by \[ \gamma \doteq \cosh\varphi = \frac{1}{\sqrt{1-\|\widetilde{\vec{v}}\|_E^2}},  \]where $\widetilde{\vec{v}} = (\Delta x_1/\Delta t,\ldots, \Delta x_{n-1}/\Delta t) \in \R^{n-1}$ is the velocity vector of the particle's trajectory in $\R^{n-1}$. Show also that $\|\widetilde{\vec{v}}\|_E = \tanh\varphi$.
\end{problem}

\begin{problem}[Coordinate-free index raising]
  Let $B\colon \R^n_\nu\times \R^n_\nu\to \R$ be a bilinear map. There is a unique linear operator $T\colon \R^n_\nu \to \R^n_\nu$ such that $B(\vec{x},\vec{y}) = \pair{T\vec{x},\vec{y}}_\nu$ for all $\vec{x},\vec{y} \in \R^n_\nu$. Show that ${\rm tr}_{\pair{\cdot,\cdot}_\nu}B = {\rm tr}\,T$ and $\det_{\pair{\cdot,\cdot}_\nu}B = (-1)^\nu \det T$.
\end{problem}

\hrulefill

\newpage

\section{Curve theory in $\LM^3$}

\begin{obs}
  All curves and functions will be assumed of class $\mathcal{C}^\infty$ (even though most of the time $\mathcal{C}^3$ or $\mathcal{C}^4$ is enough), and $I$ will always denote an open interval in $\R$.
\end{obs}

\subsection{Admissible curves and the Frenet Trihedron}

We know from classical differential geometry in Euclidean space $\R^3$ that: 
\begin{itemize}
\item any regular curve $\vec{\alpha}\colon I\to \R^3$ admits a reparametrization with unit speed, so we may assume without loss of generality that $\|\vec{\alpha}'(s)\|_E=1$;
\item we may define, for each $s \in I$, a positive orthonormal frame $(\subb{T}{\alpha}(s), \subb{N}{\alpha}(s), \subb{B}{\alpha}(s))$ for $\R^3$, pictured as attached to the point $\vec{\alpha}(s)$ --  these vectors are called the \emph{tangent}, \emph{normal} and \emph{binormal} vectors to $\vec{\alpha}$ at $s$, and they form the so-called \emph{Frenet Trihedron}\index{Frenet Trihedron} of $\vec{\alpha}$ at $s$;
\item there are functions $\kappa_{\vec{\alpha}}\colon I \to \R_{\geq 0}$ and $\tau_{\vec{\alpha}}\colon I \to \R$, called the \emph{curvature} and \emph{torsion} of $\vec{\alpha}$, such that \[
  \begin{pmatrix}
    \subb{T}{\alpha}'(s) \\ \subb{N}{\alpha}'(s) \\ \subb{B}{\alpha}'(s)
  \end{pmatrix} =
  \begin{pmatrix}
    0 & \kappa_{\vec{\alpha}}(s) & 0 \\ -\kappa_{\vec{\alpha}}(s) & 0 & \tau_{\vec{\alpha}}(s) \\ 0 & -\tau_{\vec{\alpha}}(s) & 0
  \end{pmatrix}  \begin{pmatrix}
    \subb{T}{\alpha}(s) \\ \subb{N}{\alpha}(s) \\ \subb{B}{\alpha}(s)
  \end{pmatrix},
\]for all $s \in I$.
\end{itemize}

With this data, one states and proves the Fundamental Theorem of Curves in $\R^3$, which basically says that up to rigid motions of $\R^3$, $\vec{\alpha}$ itself is determined by the functions $\kappa_{\vec{\alpha}}$ and $\tau_{\vec{\alpha}}$. More precisely:

\begin{teo}\label{teo:fundR3}
  Let $\kappa,\tau\colon I \to \R$ be given functions with $\kappa > 0$, $\vec{p}_0 \in \R^3$, $s_0 \in I$ and $(\vec{T}_0,\vec{N}_0,\vec{B}_0)$ a positive orthonormal basis for $\R^3$. Then there exists a unique unit speed regular curve $\vec{\alpha}\colon I\to \R^3$ such that:
  \begin{itemize}
  \item $\vec{\alpha}(s_0) = \vec{p}_0$;
  \item $(\subb{T}{\alpha}(s_0),\subb{N}{\alpha}(s_0), \subb{B}{\alpha}(s_0)) = (\vec{T}_0,\vec{N}_0,\vec{B}_0)$;
  \item $\kappa_{\vec{\alpha}}(s) = \kappa(s)$ and $\tau_{\vec{\alpha}}(s) = \tau(s)$ for all $s \in I$.
  \end{itemize}
\end{teo}

The proof consists, briefly speaking, in solving the Frenet system for $\vec{\alpha}$. From this point onwards, we focus our attention on three-dimensional Lorentz-Minkowski space $\LM^3$. Recall that the definition of a (parametrized) regular curve does not really depend on the scalar product we have equipped the ambient space with. And in the same way that a parametrized curve $\vec{\alpha}\colon I \to \LM^3$ is regular if $\vec{\alpha}'(t) \neq \vec{0}$ for all $t$ in $I$ (which is the same as saying that $\{\vec{\alpha}'(t)\}$ is linearly independent for all $t \in I$), we may take one step further and say that $\vec{\alpha}$ is \emph{biregular}\index{Biregular curve} if $\{\vec{\alpha}'(t),\vec{\alpha}''(t)\}$ is linearly independent for all $t \in I$. You might be (correctly) guessing what a $k$-regular curve in $\R^n_\nu$ is, by now.

Silly as this may seem, this together with a non-degenerability chain condition (such as the ones used to relate linear independence of a set of vectors with invertibility of the associated Gram matrix, or the one which allows us to perform the Gram-Schmidt orthogonalization process) is precisely what we need to adapt the classical curve theory developed in $\R^3$ for $\LM^3$.

\begin{defn}
  A curve $\vec{\alpha}\colon I \to \LM^3$ is called \emph{admissible}\index{Admissible curve} if it is biregular, and for each $t\in I$ both the \emph{tangent line} spanned by $\vec{\alpha}'(t)$ and the \emph{osculating plane}\index{Osculating plane} ${\rm span}(\vec{\alpha}'(t),\vec{\alpha}''(t))$ are non-degenerate.
\end{defn}

We might as well define the notion of causal character for curves now:

\begin{defn}
  Let $\vec{\alpha}\colon I \to \R^n_\nu$ be a regular curve and $t_0 \in I$. We say that $\vec{\alpha}$ is:
  \begin{enumerate}[(i)]
  \item \emph{spacelike at $t_0$} if $\vec{\alpha}'(t_0)$ is a spacelike vector;
  \item \emph{timelike at $t_0$} if $\vec{\alpha}'(t_0)$ is a timelike vector;
  \item \emph{lightlike at $t_0$} if $\vec{\alpha}'(t_0)$ is a lightlike vector.
  \end{enumerate}
If the causal type of $\vec{\alpha}'(t)$ is the same for all $t \in I$ according to the above, we attribute said causal type to $\vec{\alpha}$ itself. If this is the case for curves in $\LM^3$, we also define:
\begin{enumerate}[(i)]\setcounter{enumi}{3}
\item the \emph{indicator}\index{Indicator!of a curve} $\epsilon_{\vec{\alpha}}$ of $\vec{\alpha}$ to be $1$, $-1$ or $0$ if $\vec{\alpha}$ is spacelike, timelike or lightlike, respectively.
\item the \emph{coindicator}\index{Coindicator of a curve} $\eta_{\vec{\alpha}}$ of $\vec{\alpha}$ to be $1$, $-1$ or $0$ if the osculating planes are spacelike, timelike or lightlike, respectively.
\end{enumerate}
\end{defn}

With this out of the way, let's analyze the recipe described for curves in $\R^3$. First, we need a good parametrization for the curve. It turns out that for this first step, regularity is almost enough. Here's a general statement:

\begin{prop}
  Let $\vec{\alpha}\colon I \to \R^n_\nu$ be a regular curve, which is not lightlike (at any point). Then $\vec{\alpha}$ admits a reparametrization with unit speed.
\end{prop}

\begin{dem}
  Fix $t_0 \in I$ and define $s\colon I \to \R$ by \[ s(t) \doteq \int_{t_0}^t \|\vec{\alpha}'(u)\|_\nu\,\d{u}.  \]By the Fundamental Theorem of Calculus and the given hypotheses, we have that $s'(t)  = \|\vec{\alpha}'(t)\|_\nu > 0$. So $s$ is an increasing diffeomorphism from $I$ into $J \doteq s[I]$, with inverse $h\colon J \to I$. Then $\widetilde{\vec{\alpha}} \doteq \vec{\alpha}\circ h$ has unit speed.
\end{dem}

\begin{obs}
  For timelike curves in $\LM^n$, we call such parameter the \emph{proper time} of $\vec{\alpha}$ and denote it by $\mathfrak{t}$. Physically, the condition $\|\vec{\alpha}'(\mathfrak{t})\|_L=1$ says that if $\vec{\alpha}$ represents the trajectory of an observer carrying a clock, then $\mathfrak{t} - \mathfrak{t}_0$ is the time lapse measured by such observer between the events $\vec{\alpha}(\mathfrak{t}_0)$ and $\vec{\alpha}(\mathfrak{t})$.
\end{obs}

Now, the admissibility condition allows us to apply Corollary \ref{cor:non-deg} (p. \pageref{cor:non-deg}) for the osculating planes to the curve and write the:

\begin{defn}
Let $\vec{\alpha}\colon I \to \LM^3$ be a unit speed admissible curve. The \emph{tangent vector} to $\vec{\alpha}$ at $s$ is $\subb{T}{\alpha}(s)\doteq \vec{\alpha}'(s)$. Since the conditions on $\vec{\alpha}$ ensure that the \emph{curvature} of $\vec{\alpha}$ at s, $\kappa_{\vec{\alpha}}(s) \doteq \|\vec{\alpha}''(s)\|_L$, never vanishes, we may define the \emph{normal vector} to $\vec{\alpha}$ at $s$ to be the unique unit vector $\subb{N}{\alpha}(s)$ such that $\subb{T}{\alpha}'(s) = \kappa_{\vec{\alpha}}(s)\subb{T}{\alpha}(s)$. Then let the \emph{binormal vector} to $\vec{\alpha}$ at $s$, $\subb{B}{\alpha}(s)$, be the unique unit vector such that $(\subb{T}{\alpha}(s), \subb{N}{\alpha}(s), \subb{B}{\alpha}(s))$ is a positive orthonormal basis for $\LM^3$.
\end{defn}

\begin{obs}\mbox{}
  \begin{itemize}
  \item Note that $\pair{\vec{\alpha}'(s),\vec{\alpha}'(s)}_L = \epsilon_{\vec{\alpha}}$ implies $2\pair{\vec{\alpha}''(s),\vec{\alpha}'(s)}_L=0$, so indeed the vectors $\subb{T}{\alpha}(s)$ and $\subb{N}{\alpha}(s)$ are orthogonal.
  \item It follows from our previous discussion regarding orientability of bases in $\R^n_\nu$ that $\subb{B}{\alpha}(s) = (-1)^\nu\epsilon_{\vec{\alpha}}\eta_{\vec{\alpha}}\subb{T}{\alpha}(s)\times \subb{N}{\alpha}(s)$ (of course, we're interested in what will happen for $\nu=1$ here) -- this can be also checked by applying Lagrange's identity (Corollary \ref{cor:lagrange}, p. \pageref{cor:lagrange}) together with the definition of the index $\nu$ cross product as the vector representing the linear functional induced by $\det$ and its arguments.
  \end{itemize}
\end{obs}

In the same setting as the above definition, the torsion $\tau_{\vec{\alpha}}$ of $\vec{\alpha}$ will be the unique function such that \[
  \begin{pmatrix}
    \subb{T}{\alpha}'(s) \\ \subb{N}{\alpha}'(s) \\ \subb{B}{\alpha}'(s)
  \end{pmatrix} =
  \begin{pmatrix}
    0 & \kappa_{\vec{\alpha}}(s) & 0 \\ -\epsilon_{\vec{\alpha}}\eta_{\vec{\alpha}}\kappa_{\vec{\alpha}}(s) & 0 & \tau_{\vec{\alpha}}(s) \\ 0 & (-1)^{\nu+1}\epsilon_{\vec{\alpha}}\tau_{\vec{\alpha}}(s) & 0
  \end{pmatrix}  \begin{pmatrix}
    \subb{T}{\alpha}(s) \\ \subb{N}{\alpha}(s) \\ \subb{B}{\alpha}(s)
  \end{pmatrix},
\]for all $s \in I$. Setting $\nu=0$ and $\epsilon_{\vec{\alpha}}=\eta_{\vec{\alpha}}=1$, we recover the usual Frenet equations in $\R^3$. This means that the theory for admissible curves can be developed simultaneously in both ambients $\R^3$ and $\LM^3$. Here's a more powerful version of Theorem \ref{teo:fundR3}:

\begin{teo}\label{teo:fundL3}
  Let $\kappa,\tau\colon I \to \R$ be given functions with $\kappa > 0$, $\vec{p}_0 \in \R^3_\nu$, $s_0 \in I$ and $(\vec{T}_0,\vec{N}_0,\vec{B}_0)$ a positive orthonormal basis for $\R^3_\nu$. Then there exists a unique unit speed admissible curve $\vec{\alpha}\colon I\to \R^3_\nu$ such that:
  \begin{itemize}
  \item $\vec{\alpha}(s_0) = \vec{p}_0$;
  \item $(\subb{T}{\alpha}(s_0),\subb{N}{\alpha}(s_0), \subb{B}{\alpha}(s_0)) = (\vec{T}_0,\vec{N}_0,\vec{B}_0)$;
  \item $\kappa_{\vec{\alpha}}(s) = \kappa(s)$ and $\tau_{\vec{\alpha}}(s) = \tau(s)$ for all $s \in I$.
  \end{itemize}
\end{teo}

A detailed proof of this version of the Fundamental Theorem of Curves, and also how to adapt what was summarized here for admissible curves not necessarily having unit speed, see~\cite{TL}.

\subsection{Curves with lightlike osculating plane}

We will continue to work with biregular curves (without further comments). In particular, we are excluding null lines. Let's say that a unit speed non-lightlike and non-admissible curve is \emph{semi-lightlike} (observe that such curves are automatically spacelike). That is to say, a non-admissible curve is either lightlike or semi-lightlike, according to whether the tangent line or the osculating plane is degenerate. Or equivalently, a lightlike curve has $(\epsilon_{\vec{\alpha}}, \eta_{\vec{\alpha}}) = (0,1)$ while a semi-lightlike curve has $(\epsilon_{\vec{\alpha}}, \eta_{\vec{\alpha}}) = (1,0)$.

It is possible to treat both lightlike and semi-lightlike curves simultaneously. However, there is an issue we must solve first: lightlike curves obviously do not admit reparametrizations with unit speed. One can also understand the reason for this bearing in mind that lightlike curves may be seen as worldlines of photons or neutrinos -- the proper time measured by it is zero, and so it cannot be used as the curve parameter. If we cannot have $\|\vec{\alpha}'(t)\|_L=1$, we'll move on to the next best thing: $\|\vec{\alpha}''(t)\|_L=1$. More precisely:

\begin{lem}
  Let $\vec{\alpha}\colon I \to \LM^n$ be a lightlike curve with $\|\vec{\alpha}''(t)\|_L \neq 0$ for all $t \in I$. Then $\vec{\alpha}$ admits an \emph{arc-photon reparametrization}\index{Arc-photon reparametrization}. Namely, there is an open interval $J\subseteq \R$ and a diffeomorphism $h\colon J \to I$ such that $\widetilde{\vec{\alpha}} = \vec{\alpha} \circ h$ satisfies $\|\widetilde{\vec{\alpha}}''(\phi)\|_L=1$ for all $\phi \in J$.
\end{lem}

\begin{dem}
  Let's check what such $h$ must satisfy, and see if said conditions are actually enough to define it. We should have $\widetilde{\vec{\alpha}}(\phi) = \vec{\alpha}(h(\phi))$ for all $\phi \in J$, and differentiating everything twice we get \[ \widetilde{\vec{\alpha}}''(\phi) = \vec{\alpha}''(h(\phi))h'(\phi)^2 + \vec{\alpha}'(h(\phi))h''(\phi). \]Since $\vec{\alpha}$ is lightlike, $\vec{\alpha}''(h(\phi))$ is orthogonal to $\vec{\alpha}'(h(\phi))$, and the given condition \linebreak[4]$\|\vec{\alpha}''(h(\phi))\|_L \neq 0$ says that $\vec{\alpha}''(h(\phi))$ is spacelike. So, taking scalar squares on both sides yields \[  1 = \pair{\vec{\alpha}''(h(\phi)),\vec{\alpha}''(h(\phi))}_L h'(\phi)^4, \]which readily implies that $h'(\phi) = \|\vec{\alpha}''(h(\phi))\|_L^{-1/2}$. This is a first order differential equation which depends continuously on $h$, and given $\phi_0 \in J$ and $t_0 \in I$, there is a unique solution $h$ with $h(\phi_0) = t_0$. For this $h$, define $\widetilde{\vec{\alpha}} = \vec{\alpha}\circ h$. This is the desired reparametrization.
\end{dem}

\begin{Ex}
Consider the helix $\vec{\alpha}\colon \R \to \LM^3$ given by $\vec{\alpha}(t) = (r\cos t, r \sin t, rt)$, where $r>0$ is fixed. Since $\vec{\alpha}'(t) = (-r\sin t, r \cos t, r)$ is a lightlike vector for all $t \in \R$, $\vec{\alpha}$ itself is lightlike. Moreover, $\vec{\alpha}''(t) = (-r\cos t, -r\sin t,0)$ satisfies $\|\vec{\alpha}''(t)\|_L = r \neq 0$ for all $t \in \R$. So, there is an arc-photon reparametrization. The differential equation to solve becomes just $h'(\phi) = 1/\sqrt{r}$. It follows that \[  \widetilde{\vec{\alpha}}(\phi) = \left( r\cos\left(\frac{\phi}{\sqrt{r}}\right),r\sin\left(\frac{\phi}{\sqrt{r}}\right) ,\sqrt{r}\phi\right), \qquad \phi \in \R,  \]is an arc-photon reparametrization of $\vec{\alpha}$.
\end{Ex}

When treating both types of curves at the same time, we will omit the distinguished parameter $s$ or $\phi$, to avoid notation clutter. The next step is, like before, to define an adapted frame for each point in the curve. But in this case, an orthonormal frame does not carry geometric information about the curve's acceleration vector. If we cannot normalize the acceleration vector... we just don't do it. We start with the:

\begin{defn}
  Let $\vec{\alpha}:I\to\LM^3$ be a lightlike or semi-lightlike curve. We define the \emph{tangent} and \emph{normal} vectors to the curve by
  \[\subb{T}{\alpha}\doteq\vec{\alpha}'\quad\mbox{and}\quad\subb{N}{\alpha}\doteq
    \vec{\alpha}'',\] respectively.
\end{defn}

We have given up on orthonormality, but not on positiveness. To complete the frame, we need to find a third vector $\subb{B}{\alpha}$, again to be called the \emph{binormal vector}, such that the basis $(\subb{T}{\alpha},\subb{N}{\alpha},\subb{B}{\alpha})$ is positive at each point of the curve.

In general, we may define the orientation of a basis $(\vec{v},\vec{w})$ for a lightlike plane in terms of a choice of a euclidean-normal vector $\vec{n}$ to the plane. More precisely, we'll say that $(\vec{v},\vec{w})$ is \emph{positive} if $(\vec{v},\vec{w},\vec{n})$ is a positive basis for $\LM^3$, with $\vec{n}$ future-directed. If $\vec{v}$ is lightlike and $\vec{w}$ is unit (and spacelike), we have that the cross product $\vec{v}
\times_L \vec{w}$ is also lightlike, and hence proportional to $\vec{v}$. Writing $\vec{v}\times_L \vec{w} = \lambda \vec{v}$ for some $\lambda \in \R$,
we may geometrically analyze the sign of $\lambda$ as follows:

\begin{figure}[H]
  \centering
  \subfloat[$(\vec{v},\vec{w})$ positive ($\lambda<0$)]
  {\begin{picture}(0,0)%
\includegraphics{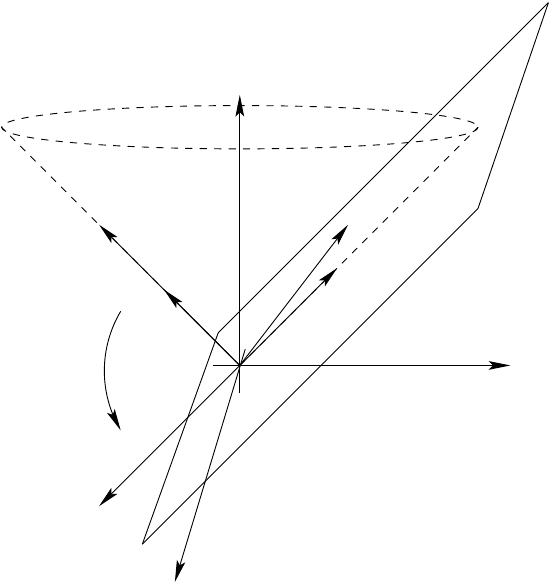}%
\end{picture}%
\setlength{\unitlength}{2279sp}%
\begingroup\makeatletter\ifx\SetFigFont\undefined%
\gdef\SetFigFont#1#2#3#4#5{%
  \reset@font\fontsize{#1}{#2pt}%
  \fontfamily{#3}\fontseries{#4}\fontshape{#5}%
  \selectfont}%
\fi\endgroup%
\begin{picture}(4569,4839)(259,-4123)
\put(3061,-1726){\makebox(0,0)[lb]{\smash{{\SetFigFont{7}{8.4}{\rmdefault}{\mddefault}{\updefault}{\color[rgb]{0,0,0}$\vec{v}$}%
}}}}
\put(3151,-1141){\makebox(0,0)[lb]{\smash{{\SetFigFont{7}{8.4}{\rmdefault}{\mddefault}{\updefault}{\color[rgb]{0,0,0}$\vec{w}$}%
}}}}
\put(1756,-1726){\makebox(0,0)[lb]{\smash{{\SetFigFont{7}{8.4}{\rmdefault}{\mddefault}{\updefault}{\color[rgb]{0,0,0}$\vec{n}$}%
}}}}
\put(811,-1186){\rotatebox{315.0}{\makebox(0,0)[lb]{\smash{{\SetFigFont{7}{8.4}{\rmdefault}{\mddefault}{\updefault}{\color[rgb]{0,0,0}$\vec{v}\times_E\vec{w}$}%
}}}}}
\put(901,-3481){\rotatebox{45.0}{\makebox(0,0)[lb]{\smash{{\SetFigFont{7}{8.4}{\rmdefault}{\mddefault}{\updefault}{\color[rgb]{0,0,0}$\vec{v}\times_L\vec{w}$}%
}}}}}
\end{picture}}\qquad
  \subfloat[$(\vec{v},\vec{w})$ negative ($\lambda>0$)]
  {\begin{picture}(0,0)%
\includegraphics{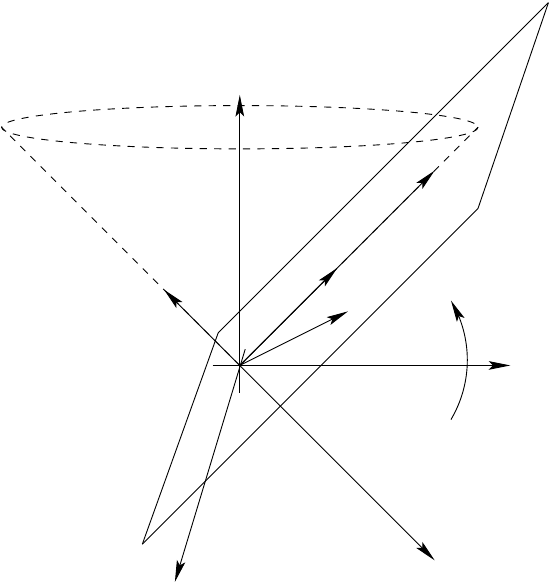}%
\end{picture}%
\setlength{\unitlength}{2279sp}%
\begingroup\makeatletter\ifx\SetFigFont\undefined%
\gdef\SetFigFont#1#2#3#4#5{%
  \reset@font\fontsize{#1}{#2pt}%
  \fontfamily{#3}\fontseries{#4}\fontshape{#5}%
  \selectfont}%
\fi\endgroup%
\begin{picture}(4569,4839)(259,-4123)
\put(1756,-1726){\makebox(0,0)[lb]{\smash{{\SetFigFont{7}{8.4}{\rmdefault}{\mddefault}{\updefault}{\color[rgb]{0,0,0}$\vec{n}$}%
}}}}
\put(2881,-1501){\makebox(0,0)[lb]{\smash{{\SetFigFont{7}{8.4}{\rmdefault}{\mddefault}{\updefault}{\color[rgb]{0,0,0}$\vec{v}$}%
}}}}
\put(2926,-1861){\makebox(0,0)[lb]{\smash{{\SetFigFont{7}{8.4}{\rmdefault}{\mddefault}{\updefault}{\color[rgb]{0,0,0}$\vec{w}$}%
}}}}
\put(2971,-3391){\rotatebox{315.0}{\makebox(0,0)[lb]{\smash{{\SetFigFont{7}{8.4}{\rmdefault}{\mddefault}{\updefault}{\color[rgb]{0,0,0}$\vec{v}\times_E\vec{w}$}%
}}}}}
\put(3601,-1276){\rotatebox{45.0}{\makebox(0,0)[lb]{\smash{{\SetFigFont{7}{8.4}{\rmdefault}{\mddefault}{\updefault}{\color[rgb]{0,0,0}$\vec{v}\times_L\vec{w}$}%
}}}}}
\end{picture}}
  \caption{Orientations for a lightlike plane.}
  \label{fig:base_luz}
\end{figure}
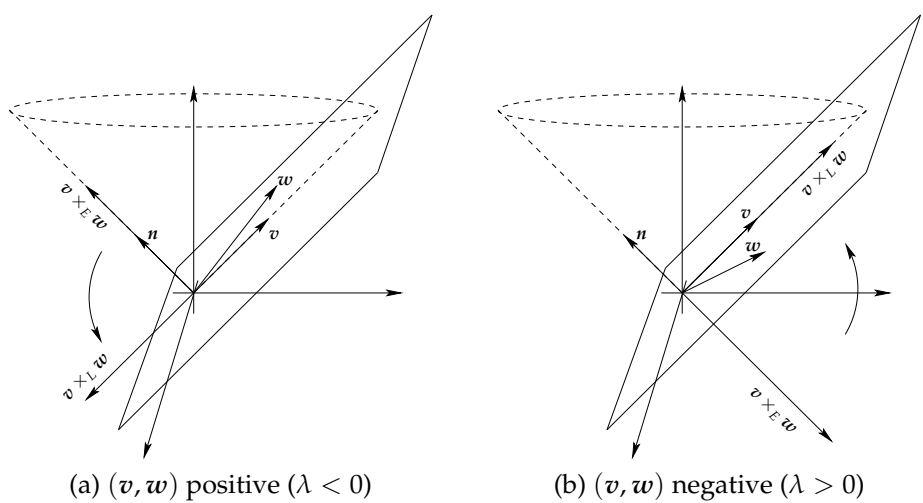

This way, if $(\vec{v},\vec{w})$ is positive then $\lambda<0$ and, similarly, if $(\vec{v},\vec{w})$ is negative we have $\lambda>0$.

Back to defining $(\subb{T}{\alpha},\subb{N}{\alpha},\subb{B}{\alpha})$: we may assume (by reparametrizing $\vec{\alpha}$ if necessary) that the bases $(\subb{T}{\alpha},\subb{N}{\alpha})$ of the osculating planes are positive. In this case, to determine the vector $\subb{B}{\alpha}$, to be lightlike, we need to also prescribe the values of $\pair{\subb{T}{\alpha},\subb{B}{\alpha}}_L$ and
$\pair{\subb{N}{\alpha},\subb{B}{\alpha}}_L$. In view of the above, one of these values should be $0$ (so that $\subb{B}{\alpha}$ is Lorentz-orthogonal to the spacelike vector) and the other $-1$ (so that $\subb{B}{\alpha}$ is not proportional to the other lightlike vector, preserving linear independence). Which of these products will be $0$ and which will be $-1$ should naturally depend on the causal type of $\vec{\alpha}$ itself. Choosing lightlike $\subb{B}{\alpha}$ such that $\pair{\subb{T}{\alpha},\subb{B}{\alpha}}_L=-\subg{\eta}{\alpha}$ and
$\pair{\subb{N}{\alpha},\subb{B}{\alpha}}_L=-\subg{\epsilon}{\alpha}$, we can treat all the cases simultaneously. So:

\begin{prop}\label{prop:tripla_pos_luz}
  Let $\vec{\alpha}:I\to\LM^3$ be a lightlike or semi-lightlike curve. The triple $(\subb{T}{\alpha},\subb{N}{\alpha},\subb{B}{\alpha})$ is a positive basis for $\LM^3$, for each point in $\vec{\alpha}$.
\end{prop}

\begin{dem}
  Our goal is to show that $\det(\subb{T}{\alpha},\subb{N}{\alpha},\subb{B}{\alpha})>0$. Let's do the case $\subg{\epsilon}{\alpha}=0$ and $\subg{\eta}{\alpha}=1$. Writing $\subb{B}{\alpha}(\phi)$ as in terms of the basis $\big(\subb{T}{\alpha}(\phi),
  \subb{N}{\alpha}(\phi),\subb{T}{\alpha}(\phi)\times_E
  \subb{N}{\alpha}(\phi)\big)$, we see that the only relevant component of
  $\subb{B}{\alpha}(\phi)$ for the determinant we're going to compute is the one in the direction of $\subb{T}{\alpha}(\phi)\times_E \subb{N}{\alpha}(\phi)$ -- call it  $\mu(\phi) \subb{T}{\alpha}(\phi)\times_E\subb{N}{\alpha}(\phi)$. Then
  \[\det(\subb{T}{\alpha}(\phi),\subb{N}{\alpha}(\phi),\subb{B}{\alpha}(\phi))=\mu(\phi)\underbrace{\det(\subb{T}{\alpha}(\phi),\subb{N}{\alpha}(\phi),\subb{T}{\alpha}(\phi)\times_E\subb{N}{\alpha}(\phi))}_{>0},\]so that we only have to verify that $\mu(\phi)>0$. From Figure \ref{fig:base_luz}, we may write that $\subb{T}{\alpha}(\phi)\times_L \subb{N}{\alpha}(\phi) = \lambda(\phi)
  \subb{T}{\alpha}(\phi)$ for a certain coefficient $\lambda(\phi) <0$ (since $(\subb{T}{\alpha}(\phi),\subb{N}{\alpha}(\phi))$ is positive). Applying ${\rm Id}_{2,1}$ on this equality, it follows that  \begin{align*}
    \subb{T}{\alpha}(\phi)\times_E \subb{N}{\alpha}(\phi)={\rm Id}_{2,1}\big(\subb{T}{\alpha}(\phi)&\times_L \subb{N}{\alpha}(\phi)\big) =
    \lambda(\phi) {\rm Id}_{2,1}\subb{T}{\alpha}(\phi)\implies \\ &\implies
    \pair{\subb{T}{\alpha}(\phi)\times_E\subb{N}{\alpha}(\phi),{\rm Id}_{2,1}\subb{T}{\alpha}(\phi)}_E<0.
  \end{align*}
  Finallly, since $\subb{T}{\alpha}(\phi)$ and $\subb{N}{\alpha}(\phi)$ are Lorentz-orthogonal to $\subb{T}{\alpha}(\phi)$, we have that
  \[-1=\pair{\subb{B}{\alpha}(\phi),\subb{T}{\alpha}(\phi)}_L
    =\mu(\phi)\pair{\subb{T}{\alpha}(\phi)\times_E\subb{N}{\alpha}(\phi),{\rm Id}_{2,1}\subb{T}{\alpha}(\phi)}_E,\]
  and so we conclude that
  $\mu(\phi)>0$.
\end{dem}

The triple $\big(\subb{T}{\alpha},\subb{N}{\alpha},\subb{B}{\alpha}\big)$ is then called the \emph{Cartan Trihedron}\index{Cartan!Trihedron} of $\vec{\alpha}$.

Geometrically, when the curve is lightlike, the situation is as follows:
the vector $\subb{N}{\alpha}(\phi)$ is spacelike, and so its orthogonal complement is a timelike plane which intersects the lightcone of $\LM^3$ in two null lines, with exactly one of them in the direction of $\subb{T}{\alpha}(\phi)$. The binormal vector is then in the direction of the \emph{other} null line in $\subb{N}{\alpha}(\phi)^\perp$, being determined by the equation $\pair{\subb{B}{\alpha}(\phi),
  \subb{T}{\alpha}(\phi)}_L=-1$. A similar interpretation can be made for semi-lightlike curves.

Now, recall that the Frenet equations arise when we write the derivatives of the vectors in the frame as a combination of the frame elements themselves. The equations were then a quick consequence of the general formula for the orthonormal expansion of a given vector -- formula that we no longer have in this setting. Here's what we have instead:

\begin{lem}\label{lem:expansao_luz}
  Let $\vec{\alpha}\colon I \to \LM^3$ and $\vec{v} \in \LM^3$. So:
  \begin{enumerate}[(i)]
  \item if $\vec{\alpha}$ is lightlike, we have \[\vec{v} = -\pair{\vec{v},\subb{B}{\alpha}(\phi)}_L\subb{T}{\alpha}(\phi)+\pair{\vec{v},\subb{N}{\alpha}(\phi)}_L\subb{N}{\alpha}(\phi)-\pair{\vec{v},\subb{T}{\alpha}(\phi)}_L\subb{B}{\alpha}(\phi),\]for all $\phi \in I$;
  \item if $\vec{\alpha}$ is semi-lightlike, we have \[\vec{v} = \pair{\vec{v},\subb{T}{\alpha}(s)}_L\subb{T}{\alpha}(s)-\pair{\vec{v},\subb{B}{\alpha}(s)}_L\subb{N}{\alpha}(s)-\pair{\vec{v},\subb{N}{\alpha}(s)}_L\subb{B}{\alpha}(s),\]for all $s \in I$.
  \end{enumerate}
\end{lem}

\begin{obs}
  One possible mnemonic is: switch the position and sign only of the coefficients corresponding to lightlike directions.
\end{obs}

\begin{dem}
  We will treat both cases at once, noting the relations $\subg{\epsilon}{\alpha}^n = \subg{\epsilon}{\alpha}$,
  $\subg{\eta}{\alpha}^n = \subg{\eta}{\alpha}$ for all $n \geq 1$,
  $\subg{\epsilon}{\alpha}\subg{\eta}{\alpha} = 0$ and
  $\subg{\epsilon}{\alpha}+\subg{\eta}{\alpha}=1$, which follow from the fact that the only possibilities of pairs are $(\subg{\epsilon}{\alpha},\subg{\eta}{\alpha}) = (1,0)$ and $(\subg{\epsilon}{\alpha},\subg{\eta}{\alpha}) = (0,1)$. Moreover, recall that we are still assuming that $(\subb{T}{\alpha},\subb{N}{\alpha})$ is positive. That being said, write $\vec{v} = a\subb{T}{\alpha}+b\subb{N}{\alpha} +
  c\subb{B}{\alpha}$. Taking all possible products with the elements of the Cartan Trihedron and organizing the results in a matrix, we get \[
\begin{pmatrix}
  \pair{\vec{v},\subb{T}{\alpha}}_L \\ \pair{\vec{v},\subb{N}{\alpha}}_L \\ \pair{\vec{v},\subb{B}{\alpha}}_L
\end{pmatrix} =
\begin{pmatrix}
  \subg{\epsilon}{\alpha} & 0 & -\subg{\eta}{\alpha} \\ 0 & \subg{\eta}{\alpha} & -\subg{\epsilon}{\alpha} \\ -\subg{\eta}{\alpha} & -\subg{\epsilon}{\alpha} & 0
\end{pmatrix}
\begin{pmatrix}
  a \\ b \\ c
\end{pmatrix}.
\] From the relations mentioned previously, the inverse of this coefficient matrix exists, and it is just the original matrix, so that: \[
\begin{pmatrix}
  a \\ b \\ c
\end{pmatrix}=
\begin{pmatrix}
  \subg{\epsilon}{\alpha} & 0 & -\subg{\eta}{\alpha} \\
  0 & \subg{\eta}{\alpha} & -\subg{\epsilon}{\alpha} \\
  -\subg{\eta}{\alpha} & -\subg{\epsilon}{\alpha} & 0
\end{pmatrix}
\begin{pmatrix}
  \pair{\vec{v},\subb{T}{\alpha}}_L \\
  \pair{\vec{v},\subb{N}{\alpha}}_L\\
  \pair{\vec{v},\subb{B}{\alpha}}_L
\end{pmatrix}.
\]We are done.
\end{dem}

Before this lemma comes into play, we have the:
\begin{defn}
  Let $\vec{\alpha}\colon I \to \LM^3$ be a lightlike or semi-lightlike curve. The \emph{pseudo-torsion}\index{Pseudo-torsion} of $\vec{\alpha}$ is the function $\subg{\ptau}{\alpha}:I\to \R$ given by $\subg{\ptau}{\alpha} \doteq -\pair{\subb{N}{\alpha}', \subb{B}{\alpha}
  }_L$.
\end{defn}  

\begin{obs}
  The function $\sub{\ptau}{\alpha}$ is also called the \emph{Cartan curvature}\index{Cartan!curvature|seealso{pseudo-torsion}} of $\vec{\alpha}$.
\end{obs}

\begin{teo}\label{teo:ref_cartan_eqs}
  Let $\vec{\alpha}\colon I \to \LM^3$ be a lightlike or semi-lightlike curve. Then we have that \[
    \begin{pmatrix}
      \subb{T}{\alpha}' \\ \subb{N}{\alpha}' \\ \subb{B}{\alpha}'
    \end{pmatrix} =
    \begin{pmatrix}
      0 & 1 & 0 \\
      \subg{\eta}{\alpha}\subg{\ptau}{\alpha} &
      \subg{\epsilon}{\alpha}\subg{\ptau}{\alpha} & \subg{\eta}{\alpha} \\
      \subg{\epsilon}{\alpha} & \subg{\eta}{\alpha}\subg{\ptau}{\alpha} &
      -\subg{\epsilon}{\alpha}\subg{\ptau}{\alpha}
    \end{pmatrix}
    \begin{pmatrix}
      \subb{T}{\alpha} \\ \subb{N}{\alpha} \\ \subb{B}{\alpha}
    \end{pmatrix}.
  \]
\end{teo}

\begin{obs}
  Explicitly, the coefficient matrices when $\vec{\alpha}$ is lightlike or semi-lightlike are, respectively,
  \[ \begin{pmatrix}
      0 & 1 & 0 \\
      \subg{\ptau}{\alpha}(\phi) &0 & 1 \\
      0 & \subg{\ptau}{\alpha}(\phi) & 0
    \end{pmatrix}\quad\mbox{e}\quad
    \begin{pmatrix}
      0 & 1 & 0 \\
      0 & \subg{\ptau}{\alpha}(s) & 0 \\
      1 & 0 & -\subg{\ptau}{\alpha}(s)
    \end{pmatrix}. \]
\end{obs}

\begin{dem}
  The first equation is the very definition of the normal vector. For the second one, we apply Lemma \ref{lem:expansao_luz} regarding $\subb{N}{\alpha}'$ as a column vector to get

  \begin{align*}
    \subb{N}{\alpha}' &=
                       \begin{pmatrix}
                         \subg{\epsilon}{\alpha} & 0 &-\subg{\eta}{\alpha}
                         \\
                         0 & \subg{\eta}{\alpha}
                         &-\subg{\epsilon}{\alpha}\\
                         -\subg{\eta}{\alpha} &
                         -\subg{\epsilon}{\alpha} & 0
                       \end{pmatrix}
     \begin{pmatrix}
       \pair{\subb{N}{\alpha}',\subb{T}{\alpha}}_L \\
       \pair{\subb{N}{\alpha}',\subb{N}{\alpha}}_L \\
       \pair{\subb{N}{\alpha}',\subb{B}{\alpha}}_L
     \end{pmatrix}\\ 
                     &=\begin{pmatrix}
                       \subg{\epsilon}{\alpha} & 0 & -\subg{\eta}{\alpha}
                       \\
                       0 & \subg{\eta}{\alpha} &
                       -\subg{\epsilon}{\alpha}\\
                       -\subg{\eta}{\alpha} & -\subg{\epsilon}{\alpha} & 0
                     \end{pmatrix}
                     \begin{pmatrix}
                       -\subg{\eta}{\alpha} \\ 0 \\ -\subg{\ptau}{\alpha}
                     \end{pmatrix}=
    \begin{pmatrix}
      \subg{\eta}{\alpha}\subg{\ptau}{\alpha} \\
      \subg{\epsilon}{\alpha}\subg{\ptau}{\alpha} \\
      \subg{\eta}{\alpha}
    \end{pmatrix},
  \end{align*}and so we have the second row of the sought coefficient matrix. Similarly for $\subb{B}{\alpha}'$, we have
  \begin{align*}
    \subb{B}{\alpha}' &=
                       \begin{pmatrix}
                         \subg{\epsilon}{\alpha} & 0 &
                         -\subg{\eta}{\alpha} \\
                         0 & \subg{\eta}{\alpha} &
                         -\subg{\epsilon}{\alpha} \\
                         -\subg{\eta}{\alpha} & -\subg{\epsilon}{\alpha} & 0
                       \end{pmatrix}
    \begin{pmatrix}
      \pair{\subb{B}{\alpha}',\subb{T}{\alpha}}_L \\
      \pair{\subb{B}{\alpha}',\subb{N}{\alpha}}_L \\
      \pair{\subb{B}{\alpha}',\subb{B}{\alpha}}_L
    \end{pmatrix} \\ 
                     &=\begin{pmatrix}
                       \subg{\epsilon}{\alpha} & 0 & -\subg{\eta}{\alpha}
                       \\
                       0 & \subg{\eta}{\alpha} & -\subg{\epsilon}{\alpha}
                       \\
                       -\subg{\eta}{\alpha} & -\subg{\epsilon}{\alpha} & 0
                     \end{pmatrix}
    \begin{pmatrix}
      \subg{\epsilon}{\alpha} \\ \subg{\ptau}{\alpha} \\ 0
    \end{pmatrix}  =
    \begin{pmatrix}
      \subg{\epsilon}{\alpha}  \\ \subg{\eta}{\alpha}\subg{\ptau}{\alpha} \\ - \subg{\epsilon}{\alpha}\subg{\ptau}{\alpha} 
    \end{pmatrix},
  \end{align*} and we obtain the last row.
\end{dem}

\begin{Ex}
  Let $r>0$ and consider again the curve $\vec{\alpha}\colon \R \to \LM^3$ given by  \[\vec{\alpha}(\phi) = \left(r \cos
      \left(\frac{\phi}{\sqrt{r}}\right), r \sin
      \left(\frac{\phi}{\sqrt{r}}\right), \sqrt{r}\phi \right),  \]which is lightlike with arc-photon parameter. We readily have
   \begin{align*}
     \subb{T}{\alpha}(\phi) &=\vec{\alpha}'(\phi) = \left(-\sqrt{r}
                             \sin\left(\frac{\phi}{\sqrt{r}}\right),
                             \sqrt{r}\cos\left(\frac{\phi}{\sqrt{r}}\right),
                             \sqrt{r} \right) \mbox{    and}\\
     \subb{N}{\alpha}(\phi) &= \vec{\alpha}''(\phi) = \left( -\cos
                             \left(\frac{\phi}{\sqrt{r}}\right), -
                             \sin\left(\frac{\phi}{\sqrt{r}}\right),0\right).
   \end{align*}
   To compute $\subb{B}{\alpha}(\phi)$, note that the cross product
   \[ \subb{T}{\alpha}(\phi) \times_E \subb{N}{\alpha}(\phi) = \left(
       \sqrt{r} \sin\left(\frac{\phi}{\sqrt{r}}\right),-\sqrt{r}
       \cos\left(\frac{\phi}{\sqrt{r}}\right) ,\sqrt{r}\right), \]seen in $\LM^3$, is lightlike and future-directed, so that the basis $(\subb{T}{\alpha}(\phi), \subb{N}{\alpha}(\phi))$ of the osculating plane is always positive (so there is no need to further reparametrize $\vec{\alpha}$). Furthermore, in this case, we have one particularity:
   $\subb{T}{\alpha}(\phi) \times_E \subb{N}{\alpha}(\phi)$ is also Lorentz-orthogonal to $\subb{N}{\alpha}(\phi)$. This implies that
   $\subb{B}{\alpha}(\phi)$ must be a positive multiple of the
   $\subb{T}{\alpha}(\phi) \times_E \subb{N}{\alpha}(\phi)$. To obtain
   $\pair{\subb{B}{\alpha}(\phi), \subb{T}{\alpha}(\phi)}_L = -1$, if suffices to take
   \[ \subb{B}{\alpha}(\phi) = \left( \frac{1}{2\sqrt{r}}
       \sin\left(\frac{\phi}{\sqrt{r}}\right),- \frac{1}{2\sqrt{r}}
       \cos\left(\frac{\phi}{\sqrt{r}}\right) ,
       \frac{1}{2\sqrt{r}}\right). \]Finally, we have:
   \[ \subg{\ptau}{\alpha}(\phi) = -\pair{\subb{N}{\alpha}'(\phi),
       \subb{B}{\alpha}(\phi)}_L = -\frac{1}{2r}
     \sin^2\left(\frac{\phi}{\sqrt{r}}\right)-\frac{1}{2r}
     \cos^2\left(\frac{\phi}{\sqrt{r}}\right) + 0 = -\frac{1}{2r}.\]
   \begin{figure}[H]
     \centering
     \includegraphics[width=.20\textwidth]{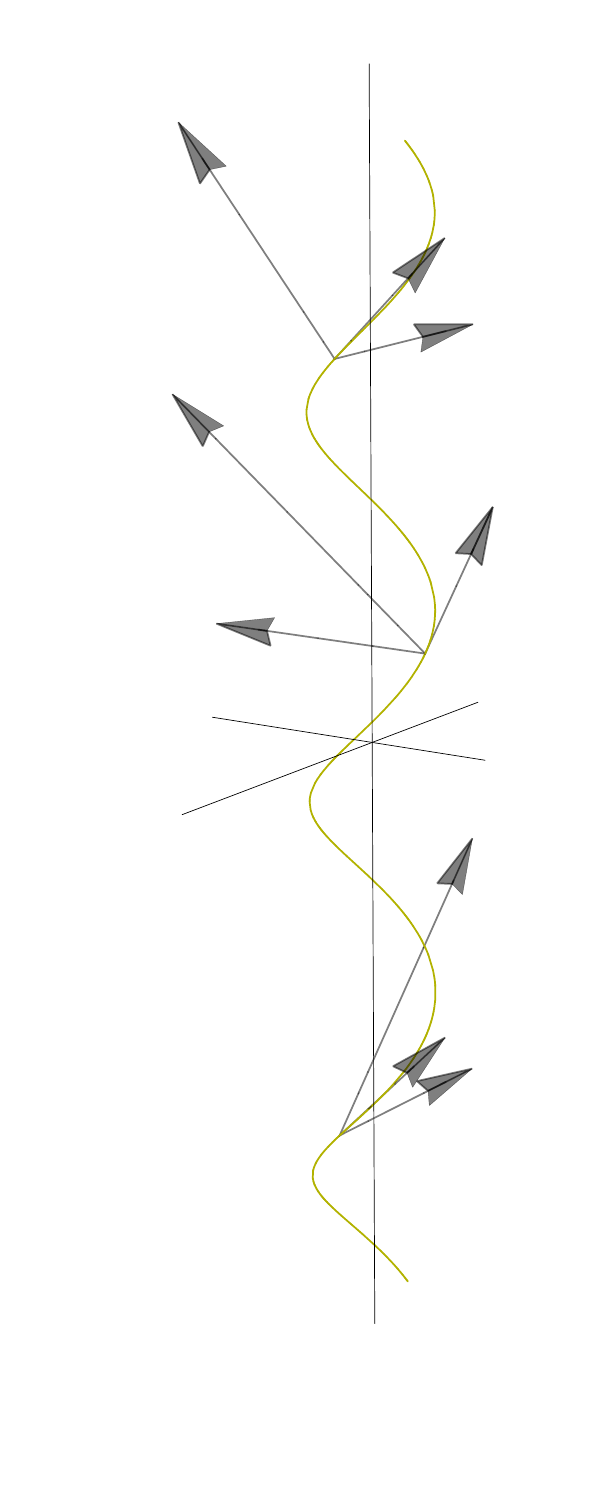}
     \caption[Cartan Trihedron for a lightlike helix.]{Cartan Trihedron for $\vec{\alpha}$ with $r=1/4$.}
     \label{fig:helice_cartan}
   \end{figure}
 \end{Ex}

From the names ``pseudo-torsion'' and ``Cartan curvature'', you might be guessing that $\sub{\ptau}{\alpha}$  should be halfway between $\subg{\kappa}{\alpha}$ and $\subg{\tau}{\alpha}$. The next two results actually show how far $\sub{\ptau}{\alpha}$ actually is from $\subg{\tau}{\alpha}$:

\begin{teo}
  The only plane lightlike curves in $\LM^3$ are null lines.
\end{teo}

\begin{dem}
  Clearly null lines are plane curves, and if $\vec{\alpha}$ is not a null line, then it has an arc-photon reparametrization. It then suffices to check that if
  $\vec{\alpha}\colon I \to \LM^3$ is a lightlike curve with arc-photon parameter and $\pair{\vec{\alpha}(\phi)- \vec{p},\vec{v}}_L=0$ for all $\phi \in I$, and certain $\vec{p},\vec{v} \in \LM^3$, then $\vec{v} = \vec{0}$. To wit, differentiating the given expression thrice we obtain:
 \[ \pair{\subb{T}{\alpha}(\phi),\vec{v}}_L =
   \pair{\subb{N}{\alpha}(\phi),\vec{v}}_L = \subg{\ptau}{\alpha}(\phi)
   \pair{\subb{T}{\alpha}(\phi),\vec{v}}_L +
   \pair{\subb{B}{\alpha}(\phi),\vec{v}}_L=0. \]If follows from Lemma
 \ref{lem:expansao_luz} (p. \pageref{lem:expansao_luz}) that $\vec{v} = \vec{0}$ as wanted.
\end{dem}

\begin{Ex}
 Let $f\colon I \to \R$ be a smooth function with positive second derivative, and consider $\vec{\alpha}\colon I \to \LM^3$ given by $\vec{\alpha}(s) = \left(s, f(s), f(s)
    \right)$. We have that $\vec{\alpha}$ is semi-lightlike with
  $\subb{T}{\alpha}(s) = \vec{\alpha}'(s) =( 1, f'(s), f'(s))$ and $\subb{N}{\alpha}(s) = \vec{\alpha}''(s) = \left(0,f''(s),f''(s) \right)$.
  Also, $\subb{T}{\alpha}(s) \times_E \subb{N}{\alpha}(s) = (0,-f''(s),f''(s))$ is a future-directed lightlike vector, so that
  $(\subb{T}{\alpha}(s),\subb{N}{\alpha}(s))$ is positive. We look for a lightlike vector $\subb{B}{\alpha}(s)=(a(s),b(s),c(s))$, Lorentz-orthogonal to $\subb{T}{\alpha}(s)$ and such that
  $\pair{\subb{B}{\alpha}(s),\subb{N}{\alpha}(s)}_L=-1$. Explicitly,
  we have the system:\[
    \begin{cases}
      a(s)^2+b(s)^2-c(s)^2 &= 0 \\ a(s)+f'(s)(b(s)-c(s)) &= 0 \\ f''(s)(b(s)-c(s)) &= -1
    \end{cases}
  \]By substituting the third equation in the second one we obtain $a(s) = f'(s)/f''(s)$. With this, the first equation becomes
  \[(b(s)-c(s))(b(s)+c(s)) = b(s)^2-c(s)^2 = -\frac{f'(s)^2}{f''(s)^2}
    \implies b(s)+c(s) = \frac{f'(s)^2}{f''(s)}, \]after using the third equation again. We then obtain
  \[ \subb{B}{\alpha}(s) = \frac{1}{2f''(s)}\left( 2 f'(s),
      f'(s)^2-1,f'(s)^2+1\right).  \]Finally, we compute
\[  \subg{\ptau}{\alpha}(s) = -\pair{\subb{N}{\alpha}'(s),
    \subb{B}{\alpha}(s)}_L = \frac{f'''(s)}{f''(s)}. \]
  In particular, note that $\vec{\alpha}$ is contained in the (lightlike) plane
   $\Pi \colon y-z = 0$, but we may choose functions $f$ for which the pseudo-torsion does not vanish.
\end{Ex}

The above example shows that, in general, the pseudo-torsion of a semi-lightlike curve is not a measure of how much the curve deviates from being a plane curve. One might wonder next whether the sign of $\subg{\ptau}{\alpha}$ says something about how the curve crosses its own osculating planes (just like $\subg{\tau}{\alpha}$ does in $\R^3$). Again, the answer is a resounding no. Let $\vec{\alpha}\colon I \to \LM^3$ be lightlike and assume that $0 \in I$ and $\vec{\alpha}(0) = \vec{0}$. Taylor expansion gives
\[ \vec{\alpha}(\phi) = \phi \vec{\alpha}'(0) + \frac{\phi^2}{2}
  \vec{\alpha}''(0) + \frac{\phi^3}{6}\vec{\alpha}'''(0) +
  \vec{R}(\phi), \]where $\vec{R}(\phi)/\phi^3 \to \vec{0}$ as $\phi \to 0$. Organizing this in terms of the Cartan Trihedron $\mathcal{F} = \big(\subb{T}{\alpha}(0), \subb{N}{\alpha}(0),
\subb{B}{\alpha}(0) \big)$, we see that the components of
$\vec{\alpha}(\phi) - \vec{R}(\phi)$ are
\[ \vec{\alpha}(\phi)-\vec{R}(\phi) = \left(\phi +
    \subg{\ptau}{\alpha}(0) \frac{\phi^3}{6}, \frac{\phi^2}{2} ,
    \frac{\phi^3}{6} \right)_{\mathcal{F}}.  \]
\begin{figure}[H]
  \centering
  \includegraphics[width=.5\textwidth]{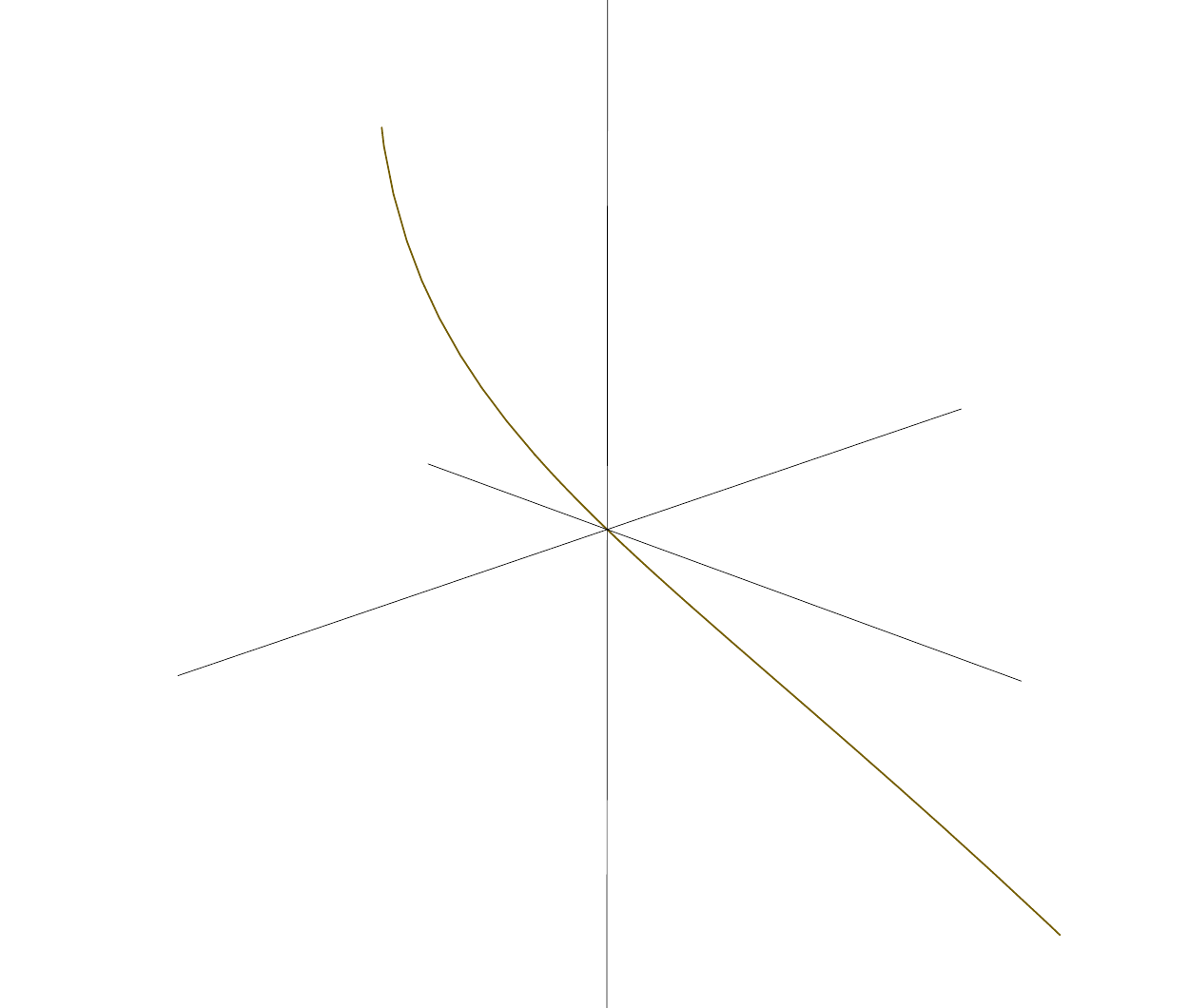}
  \caption{A ``test'' lightlike curve $\vec{\alpha}$.}
\end{figure}

Projecting, independent of the sign of $\subg{\ptau}{\alpha}(0)$, we get:
\begin{figure}[H]
  \centering
  \subfloat[Projection in the normal plane]{\includegraphics[width=.4\textwidth]{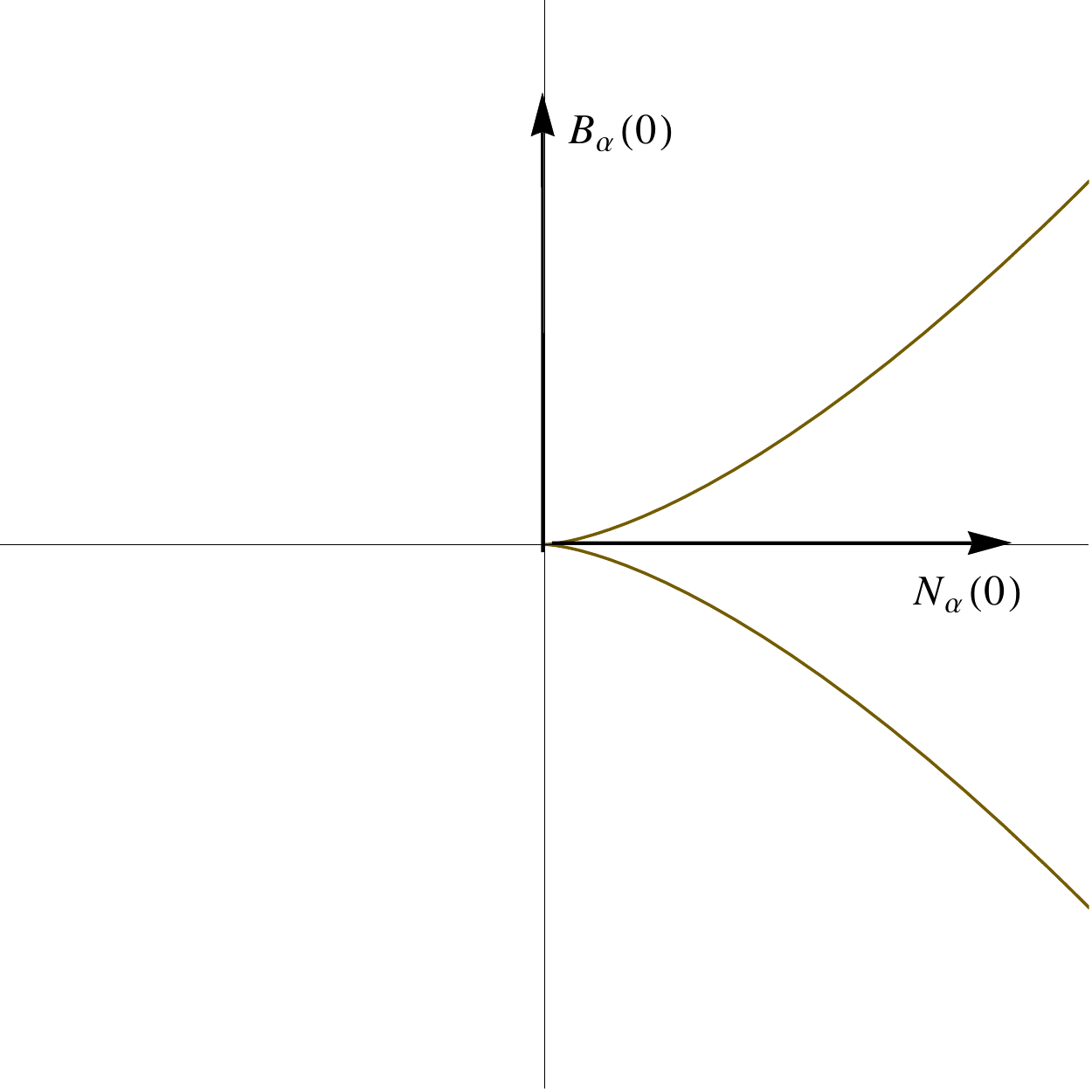}}\qquad
  \subfloat[Projection in the rectifying plane]{\includegraphics[width=.4\textwidth]{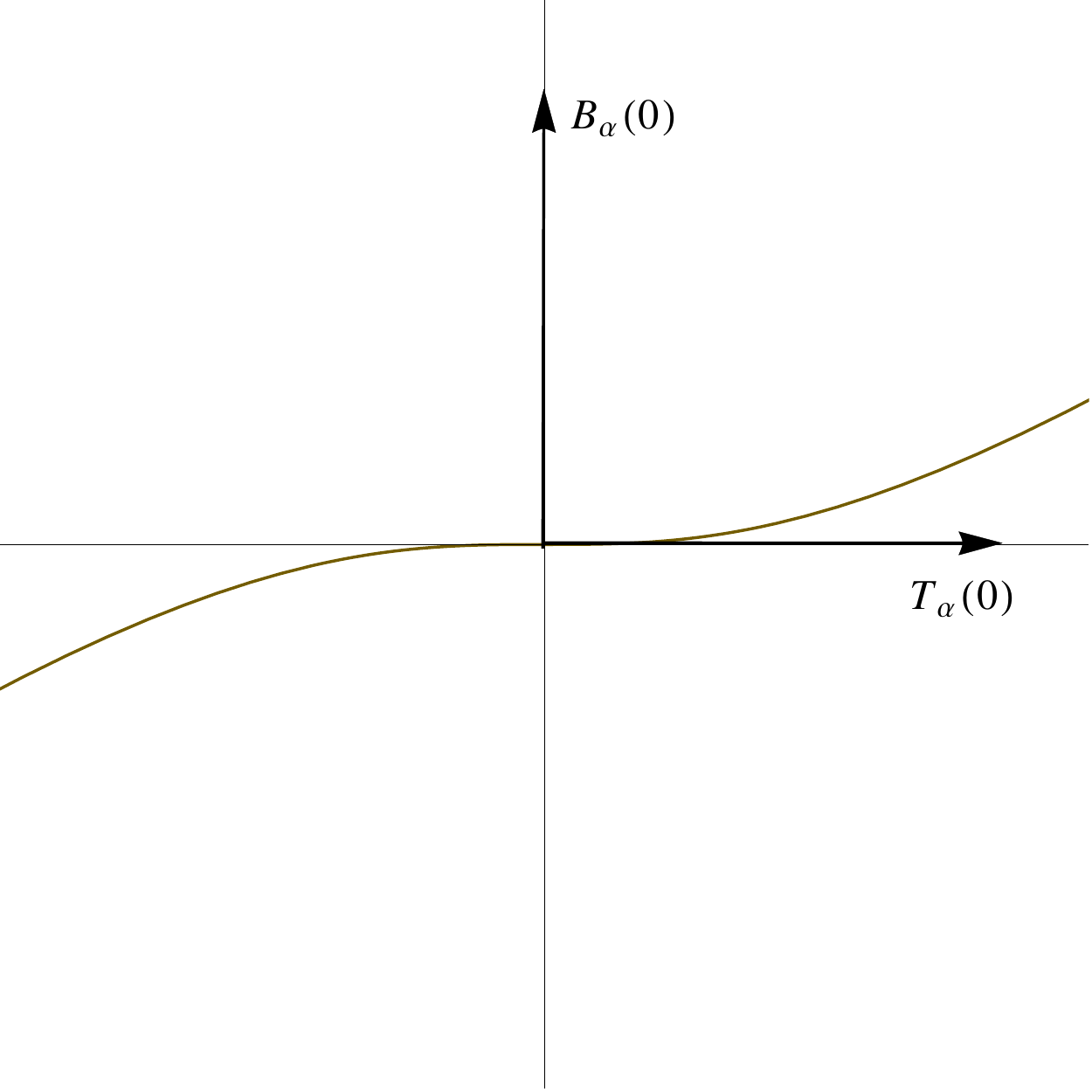}}
  \caption{Projections in the coordinate planes of the Cartan Trihedron.}
  \label{fig:interpret_ptorcao}
\end{figure}
It might be worth noting here that even though the vectors of the Cartan Trihedron are not mutually orthogonal, we may still picture them as in the above figures, bearing in mind that only their linear independence and the assumed positive orientation are relevant to concluding information about how $\vec{\alpha}$ crosses the osculating plane. We conclude that no matter the sign of the pseudo-torsion, any lightlike curve crosses its osculating planes in the direction of the binormal vector.

If $\vec{\alpha}$ is semi-lightlike instead, a similar calculation gives
\[ \vec{\alpha}(s) - \vec{R}(s) = \left( s, \frac{s^2}{2}+
    \subg{\ptau}{\alpha}(0) \frac{s^3}{6},0 \right)_{\mathcal{F}}, \]which hints at a much more extreme situation:

\begin{teo}\label{teo:semi_luz_plana}
  Every semi-lightlike curve is plane and contained in a lightlike plane.
\end{teo}

\begin{dem}
If $\vec{\alpha}\colon I \to \LM^3$ is semi-lightlike, we seek $\vec{p},\vec{v} \in \LM^3$, with lightlike $\vec{v}$, such that $\pair{\vec{\alpha}(s)-\vec{p}, \vec{v}}_L = 0$ for all $s \in I$. If this condition is satisfied, differentiating twice gives $\pair{\subb{N}{\alpha}(s),\vec{v}}_L = 0$, and we conclude that $\vec{v}$ should be proportional to $\subb{N}{\alpha}(s)$ (two Lorentz-orthogonal lightlike vectors are parallel by Corollary \ref{cor:lightlike_parallel}, p. \pageref{cor:lightlike_parallel}). Motivated by this, we seek a smooth function $\lambda\colon I \to \R$ such that $\vec{v} = \lambda(s)\subb{N}{\alpha}(s)$ is constant. This would lead us to \[ \vec{0} = (\lambda'(s) + \subg{\ptau}{\alpha}(s)\lambda(s))\subb{N}{\alpha}(s) , \]for all $s \in I$. Define $\vec{v}$ in such a way, by taking \[\lambda(s) = \exp\left(-\int_{s_0}^s \subg{\ptau}{\alpha}(\xi)\,{\dd}\xi\right),\]where $s_0 \in I$ is fixed. By construction, $\vec{v}$ is constant and then we just take $\vec{p} = \vec{\alpha}(s_0)$. This being understood, the justificative that such $\vec{p}$ e $\vec{v}$ satisfy everything we need is the usual: consider $f\colon I \to \R$ given by $f(s) = \pair{\vec{\alpha}(s)-\vec{\alpha}(s_0), \vec{v}}_L$. Clearly $f(s_0) = 0$ and $f'(s) = \pair{\subb{T}{\alpha}(s),\vec{v}}_L = 0$ for all $s \in I$.
\end{dem}

Back to the given Taylor expansion, we see that its only relevant projection is $\vec{\gamma}\colon I \to \R^2$ given by
\[ \vec{\gamma}(s) = \left( s, \frac{s^2}{2}+ \subg{\ptau}{\alpha}(0)
    \frac{s^3}{6}\right), \]and it would be natural to seek a relation between the curvature of $\vec{\gamma}$ at $0$ (as a plane curve) and the pseudo-torsion
$\subg{\ptau}{\alpha}(0)$. There is a crucial detail here, however, which will stop us from pursuing this question further: since the osculating plane is degenerate, the ``metric'' to be used in this $\R^2$ is not $\pair{\cdot,\cdot}_E$ nor $\pair{\cdot,\cdot}_L$, but the ill-behaved product
$\langle\!\langle(x_1,x_2),(y_1,y_2)\rangle\!\rangle \doteq x_1y_1$. In view of this, the expression
\[ \frac{\det(\vec{\gamma}'(s),
    \vec{\gamma}''(s))}{\|\vec{\gamma}'(s)\|^3} = 1 +
  \subg{\ptau}{\alpha}(0)s \]may no longer be seen as the curvature\footnote{Recall here that if $\vec{\gamma}\colon I \to \R^2$ is a regular plane curve in the Euclidean plane, not necessarily with unit speed, then its curvature is given by $\subg{\kappa}{\gamma}(t) = \det(\vec{\gamma}'(t),\vec{\gamma}''(t))/\|\vec{\gamma}'(t)\|^3$.}of $\vec{\gamma}$, since $\langle\!\langle\cdot,\cdot\rangle\!\rangle$ is degenerate. Even worse, there is no reasonable notion of curvature here, since every curve of the form $(s,f(s))$, where $f$ is a smooth function, can be mapped into the $x$ axis via $F\colon \R^2 \to \R^2$ given by $F(x,y) = (x,y-f(x))$. The derivative $DF(x,y)$ is a linear map which preserves
$\langle\!\langle \cdot,\cdot \rangle\!\rangle$, and so $F$ is a ``rigid motion'' of the degenerate plane. That is to say, all the graphs of smooth functions are then congruent. Now, since every spacelike curve may be parametrized as a graph over the $x$ axis and the lightlike curves are vertical lines, we conclude that it is not possible to assign a geometric invariant which distinguishes those curves.

Despite all these technical issues, the pseudo-torsion is powerful enough by itself to classify all lightlike and semi-lightlike curves in $\LM^3$ up to Poincaré transformations.

\begin{teo}\label{teo:fund_luz}\index{Fundamental Theorem!for lightlike and semi-lightlike curves}
  Let $\ptau\colon I \to \R$ be a continuous function, $\vec{p}_0 \in
  \LM^3$, $s_0,\phi_0 \in I$ and $(\vec{T}_0,\vec{N}_0,\vec{B}_0)$ a positive basis for $\LM^3$ such that $\vec{B}_0$ is a lightlike vector and $(\vec{T}_0,\vec{N}_0)$ is a positive basis for a lightlike plane. Then:
  \begin{enumerate}[(i)]
  \item if $\vec{T}_0$ is lightlike, $\vec{N}_0$ is unit spacelike and $\pair{\vec{T}_0,\vec{B}_0}_L = -1$, there is a unique lightlike curve $\vec{\alpha}\colon I \to \LM^3$ with arc-photon parameter such that
    \begin{itemize}
    \item $\vec{\alpha}(\phi_0) = \vec{p}_0$;
    \item $(\subb{T}{\alpha}(\phi_0), \subb{N}{\alpha}(\phi_0), \subb{B}{\alpha}(\phi_0)) = (\vec{T}_0,\vec{N}_0,\vec{B}_0)$;
    \item $\subg{\ptau}{\alpha}(\phi) = \ptau(\phi)$ for all $\phi \in I$.
    \end{itemize}

  \item if $\vec{T}_0$ is unit spacelike, $\vec{N}_0$ is lightlike and $\pair{\vec{N}_0,\vec{B}_0}_L = -1$, there is a unique unit speed semi-lightlike curve $\vec{\alpha}\colon I \to \LM^3$ such that
    \begin{itemize}
    \item $\vec{\alpha}(s_0) = \vec{p}_0$;
    \item $(\subb{T}{\alpha}(s_0), \subb{N}{\alpha}(s_0), \subb{B}{\alpha}(s_0)) = (\vec{T}_0,\vec{N}_0,\vec{B}_0)$;
    \item $\subg{\ptau}{\alpha}(s) = \ptau(s)$ for all $s \in I$.
    \end{itemize}
  \end{enumerate}
\end{teo}

\begin{dem}
  We will treat case (i). In a similar way done in the proof of the classical version of this result in $\R^3$, consider the following initial-value-problem in $\R^9$:
  \[
    \begin{cases}
      \begin{pmatrix}
        \vec{T}'(\phi)\\\vec{N}'(\phi)\\\vec{B}'(\phi)
      \end{pmatrix}=
      \begin{pmatrix}
        0& 1 &0\\
        \ptau(\phi) &0& 1\\
        0&\ptau(\phi)&0
      \end{pmatrix}
      \begin{pmatrix}
        \vec{T}(\phi)\\ \vec{N}(\phi)\\\vec{B}(\phi)
      \end{pmatrix}\\[.5ex]
      \mbox{e}\quad
      \big(\vec{T}(\phi_0),\vec{N}(\phi_0),\vec{B}(\phi_0)\big)=\big(\vec{T}_0,\vec{N}_0,\vec{B}_0\big).
    \end{cases}
  \]
  Such a system of linear ordinary differential equations has a unique globally defined solution
  $\big(\vec{T}(\phi),\vec{N}(\phi),\vec{B}(\phi)\big)$. We claim that this solution still satisfies, for all $\phi \in I$, the same conditions as in $\phi_0$. Namely, we will have that $\vec{T}(\phi)$ and
  $\vec{B}(\phi)$ are lightlike, $\vec{N}(\phi)$ is unit spacelike and Lorentz-orthogonal to $\vec{B}(\phi)$, and  $\pair{\vec{T}(\phi),\vec{B}(\phi)}_L=-1$. To wit, we now consider the following initial-value-problem for $\vec{a}:I\to\R^6$:
\[    \begin{cases}
      \vec{a}'(\phi)&=A(\phi)\vec{a}(\phi),\\
      \vec{a}(\phi_0)&=\big(0,1,0,0,-1,0\big),
    \end{cases}\]
 where\[ A(\phi)=\begin{pmatrix}
    0 & 0 & 0 & 2  & 0 & 0 \\ 0 & 0 & 0 & 2\ptau(\phi)
    & 0 & 2 \\
    0 & 0 & 0 & 0 & 0 & 2\ptau(\phi) \\
    \ptau(\phi) & 1 & 0 &
    0 & 1 & 0 \\ 0 & 0 & 0 &
    \ptau(\phi) & 0 & 1 \\ 0 &
    \ptau(\phi) & 1 & 0 &
    \ptau(\phi) & 0
  \end{pmatrix}.\]
 If the components of $\vec{a}(\phi)$ are all the possible products between the frame vectors\footnote{In order, $\vec{a} =
  \big(\pair{\vec{T},\vec{T}}_L,
  \pair{\vec{N},\vec{N}}_L ,
  \pair{\vec{B},\vec{B}}_L,\pair{\vec{T},\vec{N}}_L,\pair{\vec{T},\vec{B}}_L,\pair{\vec{N},\vec{B}}_L\big)$.}
 $\vec{T}(\phi)$, $\vec{N}(\phi)$ and $\vec{B}(\phi)$, we conclude that the unique solution with the given initial values is the constant vector
$\vec{a}_0=\big(0,1,0,0,-1,0\big)$, from where the claim follows. We may then define\[ \vec{\alpha}(\phi) \doteq \vec{p}_0 + \int_{\phi_0}^\phi
  \vec{T}(\xi)\,{\dd}\xi.  \]To finish the proof, we must verify that this $\vec{\alpha}$ is lightlike, has arc-photon parameter, and $\subg{\ptau}{\alpha} = \ptau$. Clearly we have
$\vec{\alpha}(\phi_0) = \vec{p}_0$ and
$\vec{\alpha}'(\phi) = \vec{T}(\phi)$, whence $\vec{\alpha}$ is lightlike. Differentiating again, we obtain $\vec{\alpha}''(\phi) = \vec{N}(\phi)$,
so that $\vec{\alpha}$ has an arc-photon parameter. This way, $\subb{T}{\alpha}(\phi) = \vec{T}(\phi)$ and
$\subb{N}{\alpha}(\phi)=\vec{N}(\phi)$, and the positivity of these bases ensure that $\subb{B}{\alpha}(\phi) = \vec{B}(\phi)$
too. Now, differentiating $\subb{N}{\alpha}(\phi) = \vec{N}(\phi)$ yields
\[ \subg{\ptau}{\alpha}(\phi)\subb{T}{\alpha}(\phi) +
  \subb{B}{\alpha}(\phi) = \ptau(\phi)\vec{T}(\phi) + \vec{B}(\phi), \]and from all the equalities seen so far it follows that
$\subg{\ptau}{\alpha}(\phi) = \ptau(\phi)$ for all $\phi \in I$. The uniqueness of such $\vec{\alpha}$ is verified in the same way as in the proof of the classical theorem: the Cartan Trihedron for another curve $\vec{\beta}$ will satisfy the same initial-value-problem, implying that $\subb{T}{\alpha} = \subb{T}{\beta}$, and so $\vec{\alpha}(\phi_0) = \vec{\beta}(\phi_0)$ gives $\vec{\alpha}=\vec{\beta}$.
\end{dem}

\begin{cor}\label{cor:uniqueness}
  Two curves, both lightlike or semi-lightlike and with the same pseudo-torsion, whose osculating planes are positively oriented, are congruent by a positive Poincaré transformation of $\LM^3$.
\end{cor}

\subsection{Lancret's theorem and classification of helices}

Here, as an application of the fundamental theorems seen so far, we can classify helices in $\R^3_\nu$. In $\R^3$, you should remember that a helix is a curve admiting a direction which makes a constant angle with all the curve's tangent lines. In $\LM^3$, a priori we can only speak of the hyperbolic angle between two timelike vectors pointing both to the future or to the past, defined just after Proposition \ref{prop:backwards_cs} (p. \pageref{prop:backwards_cs}). We would like to work with a definition of helix that works on both ambients simultaneously. Here's one:

\begin{defn}
  Let $\vec{\alpha}\colon I \to \R^3_\nu$ be a regular curve. We will say that $\vec{\alpha}$ is a \emph{helix}\index{Helix} if there is a non-zero vector $\vec{v} \in \R^3_\nu$ such that $\pair{\subb{T}{\alpha}(t), \vec{v}}$ is constant. Furthermore, in $\LM^3$, we will say that the helix is
  \begin{enumerate}[(i)]
  \item \emph{hyperbolic} if $\vec{v}$ is spacelike;
  \item \emph{elliptic} if $\vec{v}$ is timelike;
  \item \emph{parabolic} if $\vec{v}$ is lightlike.
  \end{enumerate}
The direction defined by $\vec{v}$ is called the \emph{helical axis} of $\vec{\alpha}$.
\end{defn}

For admissible curves, we have the:

\begin{teo}[Lancret]\label{teo:lancret_adm}
  Let $\vec{\alpha}:I\to\R^3_\nu$ be a unit speed admissible curve. Then $\vec{\alpha}$ is a helix if and only if the ratio
  $\subg{\tau}{\alpha}(s)/\subg{\kappa}{\alpha}(s)$ is constant.
\end{teo}

\begin{dem}
  Assume that $\vec{\alpha}$ is a helix whose helical axis is given by a vector
  $\vec{v}$. If we define $c\doteq\pair{\sub{T}{\alpha}(s), \vec{v}}$, then $\pair{\subg{\kappa}{\alpha}(s)\sub{N}{\alpha}(s),
      \vec{v}}=0$ readily implies $\pair{\sub{N}{\alpha}(s), \vec{v}}=0$, since we have $\subg{\kappa}{\alpha}(s)\neq 0$. Differentiating again, we get \[-\subg{\epsilon}{\alpha}\subg{\eta}{\alpha}\subg{\kappa}{\alpha}(s)c+\subg{\tau}{\alpha}(s)\pair{\sub{B}{\alpha}(s),
      \vec{v}}=0.\] To see that the ratio $\subg{\tau}{\alpha}(s)/\subg{\kappa}{\alpha}$ is constant, it suffices to verify that $\pair{\sub{B}{\alpha}(s),\vec{v}}$ is a non-zero constant. To wit:
  \[\frac{{\dd}}{{\dd}s}
    \pair{\sub{B}{\alpha}(s),\vec{v}}=
    (-1)^{\nu+1}\subg{\epsilon}{\alpha}\subg{\tau}{\alpha}(s)
    \pair{\sub{N}{\alpha}(s),\vec{v}}=0.\]Now, if $\pair{\sub{B}{\alpha}(s),\vec{v}} = 0$ for all $s$, then $c = 0$, and orthonormal expansion yields $\vec{v} = \vec{0}$, contradicting the definition of helix. Hence $\subg{\tau}{\alpha}(s)/\subg{\kappa}{\alpha}(s)$ is a constant. 

 Conversely, assume that $\subg{\tau}{\alpha}(s)=c\subg{\kappa}{\alpha}(s)$, for some $c\in\R$. If $c=0$ then $\vec{\alpha}$ is a plane curve and then
  $\sub{B}{\alpha}(s)=\vec{B}$ defines the helical axis for $\vec{\alpha}$. If $c\neq0$, we seek a constant vector  \[\vec{v}=v_1(s)\sub{T}{\alpha}(s)+v_2(s)\sub{N}{\alpha}(s)+v_3(s)\sub{B}{\alpha}(s)\]
  such that $\pair{\sub{T}{\alpha}(s),\vec{v}}$ is also constant. This condition, in turn, is equivalent to $v_1(s)=v_1$ being constant. Differentiating the expression for $\vec{v}$ gives us that
  \begin{align*}
    \vec{0}&=-\subg{\epsilon}{\alpha}\subg{\eta}{\alpha}\subg{\kappa}{\alpha}(s)v_2(s)\sub{T}{\alpha}(s)\\
             &\qquad+\left(v_1\subg{\kappa}{\alpha}(s)+v_2'(s)+(-1)^{\nu+1}\subg{\epsilon}{\alpha}c\subg{\kappa}{\alpha}(s)v_3(s)\right)\sub{N}{\alpha}(s)\\
           &\qquad
             +\left(c\subg{\kappa}{\alpha}(s)v_2(s)+v_3'(s)\right)\sub{B}{\alpha}(s).
  \end{align*}
  Now, linear independence implies that \[
    \begin{cases}
      0=-\subg{\epsilon}{\alpha}\subg{\eta}{\alpha}\subg{\kappa}{\alpha}(s)v_2(s)\\
      0=v_1\subg{\kappa}{\alpha}(s)+v_2'(s)+(-1)^{\nu+1}\subg{\epsilon}{\alpha}c\subg{\kappa}{\alpha}(s)v_3(s),\\
      0=c\subg{\kappa}{\alpha}(s)v_2(s)+v_3'(s).
    \end{cases}\] and hence
  \[v_2(s)=0\quad\mbox{e}\quad
    v_3(s)=\frac{(-1)^\nu}{c}\subg{\epsilon}{\alpha}v_1.\] Effectively, we have parametrized the helical axis for $\vec{\alpha}$, using $v_1$ as the real parameter. For example, setting $v_1=1$ we may see that
  \[\vec{v}\doteq\sub{T}{\alpha}(s)+\frac{(-1)^\nu}{c}\subg{\epsilon}{\alpha}\sub{B}{\alpha}(s)\]
  defines the helical axis for $\vec{\alpha}$.
\end{dem}

\begin{obs}
  \mbox{}
  \begin{itemize}
    \item In particular, this proof ensures the existence of precisely one helical axis for a given helix.
    \item If $\vec{\alpha}$ is a parabolic helix, then
    $\subg{\tau}{\alpha}(s)=\pm\subg{\kappa}{\alpha}(s)$. The converse holds provided that $\subg{\eta}{\alpha}=1$.
  \end{itemize}
\end{obs}

\begin{cor}\label{cor:class_helice_cte}
  A unit speed admissible helix $\vec{\alpha}\colon I \to \R^3_\nu$ with both constant curvature and constant torsion is congruent, for a certain choice of $a,b \in \R$, to a piece of precisely one of the following \emph{standard helices}:
  \begin{itemize}
    \item $\vec{\beta}_1(s) = \big(a \cos (s/c), a \sin (s/c), bs/c\big)$;
    \item $\vec{\beta}_2(s)=\big(a \cos(s/c), a\sin(s/c),bs/c\big)$;
      \item $\vec{\beta}_3(s) = \big(bs/c,a\cosh(s/c),
      a\sinh(s/c)\big)$;
      \item
      $\vec{\beta}_4(s) = \big(bs/c, a\sinh(s/c), a\cosh(s/c) \big)$;
      \item $\vec{\beta}_5(s) = \big(a s^2/2, a^2s^3/6,s+a^2s^3/6\big)$;
      \item
      $\vec{\beta}_6(s)=\big(as^2/2, s-a^2s^3/6,-a^2s^3/6\big)$,
    \end{itemize}
    where $\vec{\beta}_1$ is seen in $\R^3$, the remaining ones in $\LM^3$, $c \doteq \sqrt{a^2+b^2}$ for $\vec{\beta}_1$ and
    $\vec{\beta}_4$, and \linebreak[4]$c \doteq \sqrt{|a^2-b^2|}$ for $\vec{\beta}_2$
    and $\vec{\beta}_3$;
\end{cor}

\begin{dem}
  Let's denote the curvature and torsion of $\vec{\alpha}$, respectively, by $\kappa$ and $\tau$. If $\vec{\alpha}$ is seen in $\R^3$, then it is congruent to
  $\vec{\beta}_1$. We focus then on what happens in $\LM^3$. One vector spanning the helical axis is
  \[\vec{v}=\sub{T}{\alpha}(s)-\frac{\subg{\epsilon}{\alpha}\kappa}{\tau}\sub{B}{\alpha}(s),\]
  whence  $\pair{\vec{v},\vec{v}}_L=\subg{\epsilon}{\alpha}\big(1-\subg{\eta}{\alpha}\kappa^2/\tau^2\big)$. In general, the causal type of all the curves given in the statement of the result is determined by the constants $a$ and $b$. Thus, a timelike helix is:
  \begin{itemize}
    \item hyperbolic if $\kappa>|\tau|$, and hence congruent to $\vec{\beta}_3$;
    \item elliptic if $\kappa<|\tau|$, and hence congruent to $\vec{\beta}_2$, and;
    \item parabolic if $\kappa=|\tau|$, and hence congruent to $\vec{\beta}_5$.
  \end{itemize}  
  Similarly, a spacelike helic with timelike normal is necessarily hyperbolic, and so it is congruent to $\vec{\beta}_4$. Lastly, a spacelike helix with timelike binormal is:
  \begin{itemize}
    \item hyperbolic if $\kappa<|\tau|$, and hence congruent to
    $\vec{\beta}_3$;
    \item elliptic if $\kappa>|\tau|$, and hence congruent to
    $\vec{\beta}_2$, and;
    \item parabolic if $\kappa=|\tau|$, and hence congruent to $\vec{\beta}_6$.
  \end{itemize}
\end{dem}

\begin{obs}
  In each case above, it is possible to find out what $a$ and $b$ should be in terms of $\kappa$ and $\tau$. Have fun (or not).
\end{obs}

Now, we move on to non-admissible curves. Since every semi-lightlike curve is plane, it is automatically a helix. For lightlike curves the situation becomes interesting again, and we have the:

\begin{teo}[Lancret, lightlike version]\label{teo:lancret_luz}
  Let $\vec{\alpha}\colon I \to \LM^3$ be a lightlike curve with arc-photon parameter. Then $\vec{\alpha}$ is a helix if and only if its pseudo-torsion $\subg{\ptau}{\alpha}$ is constant.
\end{teo}

\begin{dem}
  Assume that $\vec{\alpha}$ is a helix and let $\vec{v} \in \LM^3$ define the helical axis. Namely, $\vec{v}$ is such that $\pair{\subb{T}{\alpha}(\phi), \vec{v}}_L = c \in \R$ is
  constant. Differentiating that relation twice we directly obtain
  \[\pair{\subb{N}{\alpha}(\phi),\vec{v}}_L= \subg{\ptau}{\alpha}(\phi)c +
    \pair{\subb{B}{\alpha}(\phi),\vec{v}}_L = 0 \]for all $\phi \in
  I$. We claim that $c \neq 0$ and that $\pair{\subb{B}{\alpha}(\phi),\vec{v}}_L$ is constant, whence $\subg{\ptau}{\alpha}(\phi)$ is also constant. To wit, if $c=0$ then Lemma
  \ref{lem:expansao_luz} (p. \pageref{lem:expansao_luz}) says that $\vec{v}=\vec{0}$, contradicting the definition of helix. Moreover, we have
  \[ \frac{{\dd}}{{\dd}\phi} \pair{\subb{B}{\alpha}(\phi),\vec{v}}_L
    = \subg{\ptau}{\alpha}(\phi) \pair{\subb{N}{\alpha}(\phi),\vec{v}}_L =
    0.  \]Conversely, assume that $\subg{\ptau}{\alpha}(\phi) = \ptau$ is a constant. If $\ptau = 0$,
  then $\vec{v} = \subb{B}{\alpha}(\phi)$ defines the helical axis for $\vec{\alpha}$. If $\ptau \neq 0$, define
  \[\vec{v} \doteq \subb{T}{\alpha}(\phi) - \frac{1}{\ptau}
    \subb{B}{\alpha}(\phi).\]Indeed, we have that
  \[ \frac{{\dd}\vec{v}}{{\dd}\phi} = \subb{N}{\alpha}(\phi) -
    \frac{1}{\ptau} \ptau \subb{N}{\alpha}(\phi) = \vec{0} \]so that $\vec{v}$
  is constant. It follows that
  $\pair{\subb{T}{\alpha}(\phi), \vec{v}}_L = 1/\ptau$ is constant, as wanted.
\end{dem}

\begin{cor}\label{cor:class_helice_luz}
  A lightlike helix $\vec{\alpha}\colon I \to \LM^3$ is congruent, for a certain choice of $r>0$, to a piece of precisely one of the following \emph{standard helices}:
    \begin{itemize}
      \item $\vec{\gamma}_1(\phi) = \big(\sqrt{r}\phi,
      r\cosh(\phi/\sqrt{r}), r\sinh(\phi/\sqrt{r})\big)$;
      \item
      $\vec{\gamma}_2(\phi)=\bigl(r \cos(\phi/\sqrt{r}), r \sin
      (\phi/\sqrt{r}), \sqrt{r}\phi\bigr)$;
      \item
      $\vec{\gamma}_3(\phi) =
      \left(-\dfrac{\phi^3}{4}+\dfrac{\phi}{3},\dfrac{\phi^2}{2},
      -\dfrac{\phi^3}{4}-\dfrac{\phi}{3}\right)$.
    \end{itemize}
\end{cor}

\begin{dem}
  Let's denote the constant pseudo-torsion of $\vec{\alpha}$ simply by $\ptau$. We know from the previous proof that a vector defining the helical axis of $\vec{\alpha}$ if $\ptau \neq 0$ is
  \[\vec{v}=\subb{T}{\alpha}(\phi)-\frac{1}{\ptau}\subb{B}{\alpha}(\phi),\]
  whence $\pair{\vec{v},\vec{v}}_L=2/\ptau$, while we may take
  $\vec{v} = \subb{B}{\alpha}(\phi)$ if $\ptau = 0$ (and hence $\pair{\vec{v},\vec{v}}_L=0$). Thus, we have that $\vec{\alpha}$ is
  \begin{itemize}
    \item hyperbolic if $\ptau > 0$, and hence congruent to $\vec{\gamma}_1$;
    \item elliptic if $\ptau < 0$, and hence congruent to $\vec{\gamma}_2$, and;
    \item parabolic if $\ptau = 0$, and hence congruent to $\vec{\gamma}_3$.
  \end{itemize}  
\end{dem}

\newpage
\section*{Problems}\addcontentsline{toc}{subsection}{Problems for the second section}

\begin{problem}
  Let $\vec{\alpha}\colon I \to \LM^n$ be a timelike future-directed curve (i.e., each $\vec{\alpha}'(t)$ is future-directed), and $a,b \in I$ with $a<b$. Show that:
  \begin{enumerate}[(a)]
  \item the difference $\vec{\alpha}(b)-\vec{\alpha}(a)$ is timelike and future-directed.
  \item $\displaystyle{\int_a^b \|\vec{\alpha}'(u)\|_L\,\d{u} \leq \|\vec{\alpha}(b)-\vec{\alpha}(a)\|_L}$, and equality holds if and only if the image of the restriction $\vec{\alpha}|_{]a,b[}$ is the line segment joining $\vec{\alpha}(a)$ and $\vec{\alpha}(b)$. What does this mean physically?
  \end{enumerate}
  \begin{dica}
    In (a), write $\vec{\alpha} = (\vec{\beta}, x_n)$, where $\vec{\beta}\colon I \to \R^{n-1}$, and estimate $\|\vec{\beta}(b)-\vec{\beta}(a)\|_E$. For (b), use the backwards Cauchy-Schwarz inequality (Proposition \ref{prop:backwards_cs}, p. \pageref{prop:backwards_cs}).
  \end{dica}
\end{problem}

\begin{problem}
  Let $\vec{\alpha}\colon I \to \LM^3$ be a unit speed admissible curve. Show that $\vec{\alpha}$ is a plane curve if and only if $\tau_{\vec{\alpha}}=0$.
\end{problem}

\begin{problem}
  Let $\vec{\alpha}\colon I \to \LM^n$ be a lightlike curve, and suppose that $\widetilde{\vec{\alpha}}_1\colon J_1 \to \LM^n$ and $\widetilde{\vec{\alpha}}_2\colon J_2 \to \LM^n$ are two arc-photon reparametrizations of $\vec{\alpha}$, so that $\widetilde{\vec{\alpha}}_1(\phi_1(t)) = \widetilde{\vec{\alpha}}_2(\phi_2(t))$. Show that $\phi_1(t) = \phi_2(t)+a$ for some $a \in \R$. What is the meaning of the constant $a$?
\end{problem}

\begin{problem}
  Check the remaining case $\epsilon_{\vec{\alpha}} = 1$ and $\eta_{\vec{\alpha}} = 0$ mentioned in the proof of Proposition \ref{prop:tripla_pos_luz} (p. \pageref{prop:tripla_pos_luz}).
\end{problem}

\begin{problem}
  Find the Cartan Trihedron and the pseudo-torsion of $\vec{\alpha}\colon \R \to \LM^3$ given by \[  \vec{\alpha}(\phi) = \left(\sqrt{r}\phi, r\cosh\left(\frac{\phi}{\sqrt{r}}\right), r\sinh \left(\frac{\phi}{\sqrt{r}}\right)\right), \]where $r>0$ is fixed.
\end{problem}

\begin{problem}
  Work through the proof of case (ii) in Theorem \ref{teo:fund_luz} (p. \pageref{teo:fund_luz}).
\end{problem}

\begin{problem}
  Show Corollary \ref{cor:uniqueness} (p. \pageref{cor:uniqueness}).
\end{problem}

\begin{problem}
  Show that every semi-lightlike curve $\vec{\alpha}\colon I \to \LM^3$ with non-zero constant pseudo-torsion $\sub{\ptau}{\alpha}(s) = \ptau \neq 0$, contained in the plane $\Pi\colon y=z$, is of the form \[ \vec{\alpha}(s) = \left(\pm s + a, \frac{b}{\ptau^2}e^{\ptau s}+cs+d,\frac{b}{\ptau^2}e^{\ptau s}+cs+d\right),  \]for some constants $a,b,c,d \in \R$ (perhaps up to reparametrization).
\end{problem}

\hrulefill

\newpage

\section{Surface theory}

\subsection{Causal characters (once more) and curvatures}

The usual definition of a regular surface in $\R^3$ (embedded, with no self-intersections) does not depend whatsoever of the product $\pair{\cdot,\cdot}_E$, and so it still makes perfect sense in $\LM^3$. This way, all the theory regarding the topology and calculus on surfaces is still valid and applicable here. In particular, we assume known:
\begin{itemize}
\item that inverse images of regular values of real-valued smooth functions on $\R^3$ are regular surfaces;
\item what is the tangent plane to a surface at any given point;
\item what is the differential of a smooth function defined in a surface, as well as what are its partial derivatives computed with respect to a given coordinate chart.
\end{itemize}
For example, since $1$ and $-1$ are both regular values of the scalar square function $F\colon \LM^3 \to \R$ given by $F(\vec{p}) = \pair{\vec{p},\vec{p}}_L$, we conclude that the \emph{de Sitter space}\index{de Sitter space} $\esf^2_1 = F^{-1}(1)$ and the \emph{hyperbolic plane}\index{Hyperbolic!plane} $\H^2$ (the upper connected component of $F^{-1}(-1)$) are regular surfaces:

  \begin{figure}[H]
    \centering \subfloat[$\esf^2_1$]
    {\includegraphics[width =
      .3\textwidth]{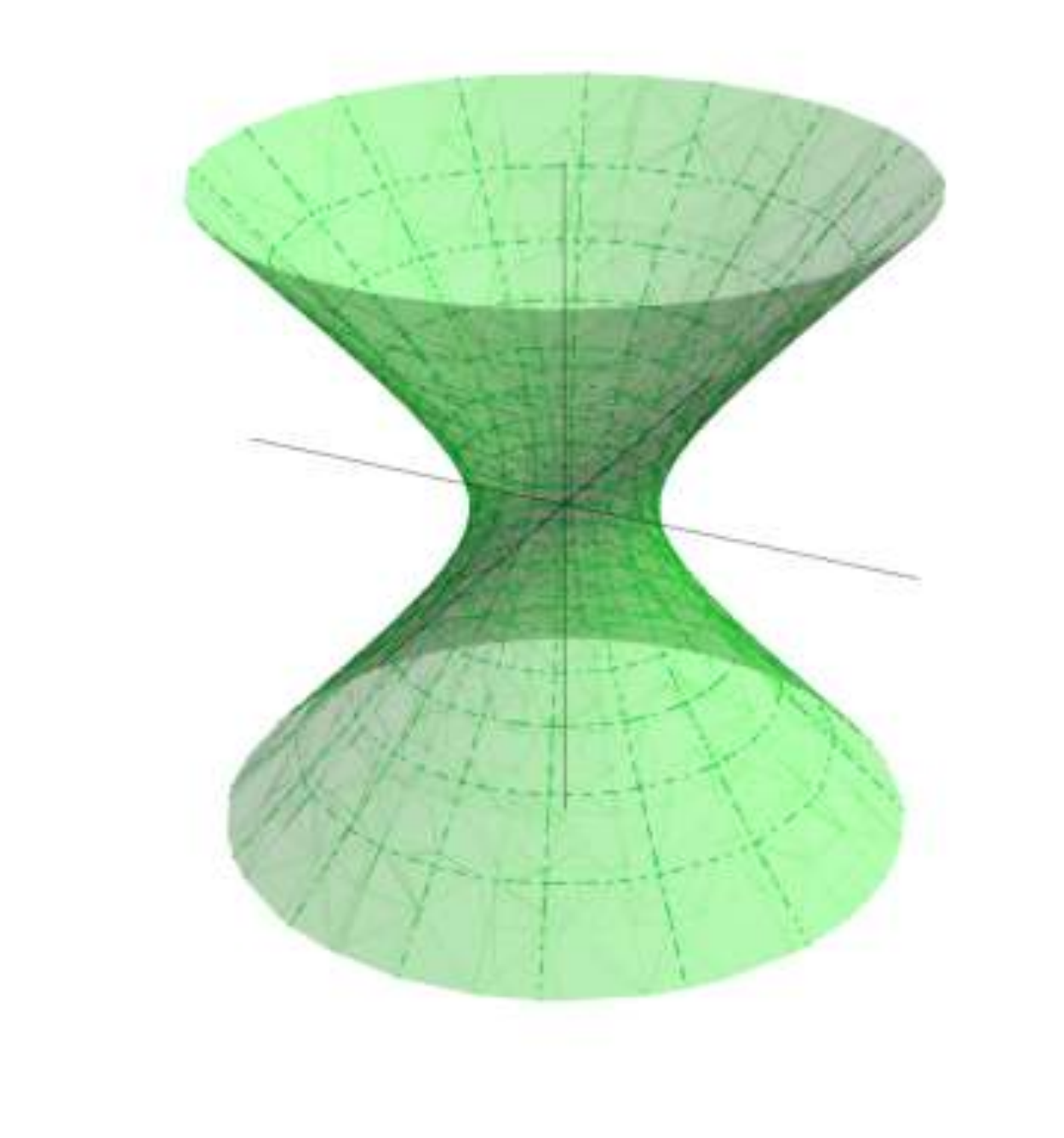}}\qquad
    \subfloat[$\H^2\cup\H^2_-$]
    {\includegraphics[width = .28\textwidth]
      {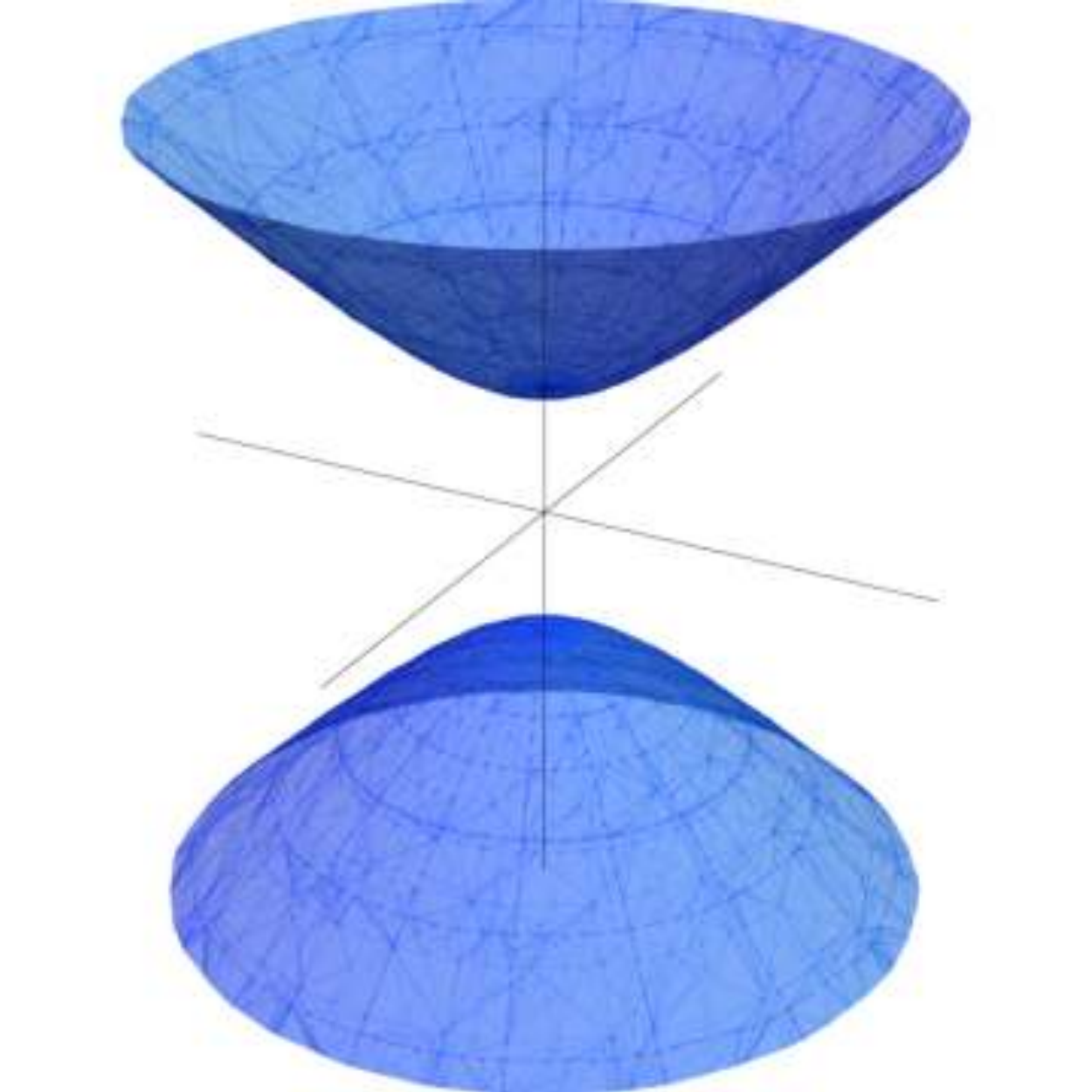}}
    \caption{The ``spheres'' in $\LM^3$.}
  \end{figure}

The product $\pair{\cdot,\cdot}_L$ comes into play when we want to generalize the notion of causal character to surfaces:

\begin{defn}
Let $M\subseteq \LM^3$ be a regular surface. We'll say that $M$ is:
\begin{enumerate}[(i)]
\item \emph{spacelike} if, for all $\vec{p} \in M$, $T_{\vec{p}}M$ is a spacelike plane;
\item \emph{timelike} if, for all $\vec{p} \in M$, $T_{\vec{p}}M$ is a timelike plane;
\item \emph{lightlike} if, for all $\vec{p} \in M$, $T_{\vec{p}}M$ is a lightlike plane.
\end{enumerate}In particular, we'll say that $M$ is \emph{non-degenerate} if no tangent plane $T_{\vec{p}}M$ is lightlike (and \emph{degenerate} otherwise). In this case, the \emph{indicator}\index{Indicator!of a surface} $\epsilon_M$ of $M$ will be $-1$ or $1$ according to whether $M$ is spacelike or timelike.
\end{defn}

\begin{Ex}\label{Ex:tipos_causais}\mbox{}
  \begin{enumerate}[(1)]
  \item Let $U\subseteq \R^2$ be open, and $f\colon U \to \R$ be a smooth function. The graph ${\rm gr}(f)$ is:
    \begin{itemize}
    \item spacelike, if $\|\nabla F\| < 1$;
    \item timelike, if $\|\nabla F\| > 1$;
    \item lightlike, if $\|\nabla F\|=1$.
    \end{itemize}
 Here's a picture:
     \begin{figure}[H]
      \centering
     \begin{picture}(0,0)%
\includegraphics{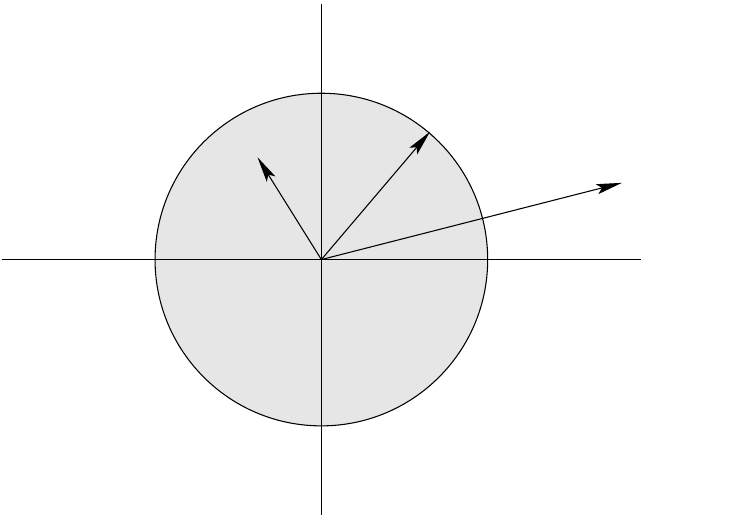}%
\end{picture}%
\setlength{\unitlength}{2693sp}%
\begingroup\makeatletter\ifx\SetFigFont\undefined%
\gdef\SetFigFont#1#2#3#4#5{%
  \reset@font\fontsize{#1}{#2pt}%
  \fontfamily{#3}\fontseries{#4}\fontshape{#5}%
  \selectfont}%
\fi\endgroup%
\begin{picture}(5175,3624)(2689,-5923)
\put(6211,-4426){\makebox(0,0)[lb]{\smash{{\SetFigFont{12}{14.4}{\rmdefault}{\mddefault}{\updefault}{\color[rgb]{0,0,0}$1$}%
}}}}
\put(7111,-3616){\makebox(0,0)[lb]{\smash{{\SetFigFont{12}{14.4}{\rmdefault}{\mddefault}{\updefault}{\color[rgb]{0,0,0}timelike}%
}}}}
\put(5806,-3166){\makebox(0,0)[lb]{\smash{{\SetFigFont{12}{14.4}{\rmdefault}{\mddefault}{\updefault}{\color[rgb]{0,0,0}lightlike}%
}}}}
\put(4996,-3481){\makebox(0,0)[lb]{\smash{{\SetFigFont{12}{14.4}{\rmdefault}{\mddefault}{\updefault}{\color[rgb]{0,0,0}$\nabla f$}%
}}}}
\put(3106,-3076){\makebox(0,0)[lb]{\smash{{\SetFigFont{12}{14.4}{\rmdefault}{\mddefault}{\updefault}{\color[rgb]{0,0,0}spacelike}%
}}}}
\end{picture}
      \caption{Finding the causal character of graphs over the plane $z=0$.}\label{fig:radar_grafico}
    \end{figure}
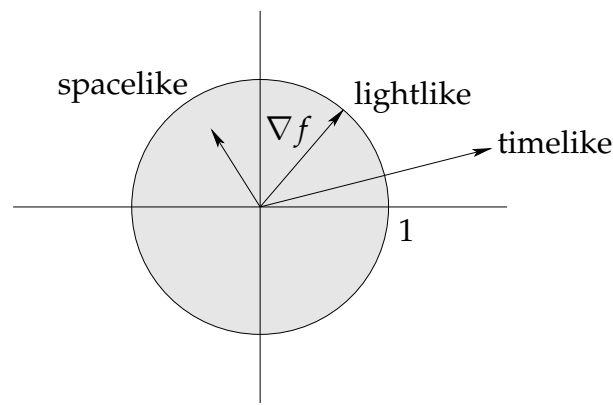

  \item Let $F\colon \LM^3 \to \R$ be a smooth function, $a \in \R$ a regular value for $F$, and $M = F^{-1}(a)$ a level surface. Then it follows from Theorem \ref{teo:causal_perp} (p. \pageref{teo:causal_perp}) that $M$ is spacelike (resp. timelike, lightlike) if and only if the usual gradient $\nabla F$ is always timelike (resp. spacelike, lightlike).
  \item If $\vec{\alpha}\colon I \to \LM^3$ is a smooth, regular and injective curve whose trace lies in the plane $y = 0$ but does not touch the $z$-axis, then we obtain a regular surface $M$ by rotating $\vec{\alpha}$ around the $z$-axis. The causal character of $M$ is the same one as $\vec{\alpha}$'s. One can understand this by noting that the parallels of revolution are always spacelike, so the only way of obtaining a lightlike or timelike direction comes from a possible contribution from $\vec{\alpha}$.
    \begin{figure}[H]
      \centering\includegraphics[width
      =.35\textwidth]{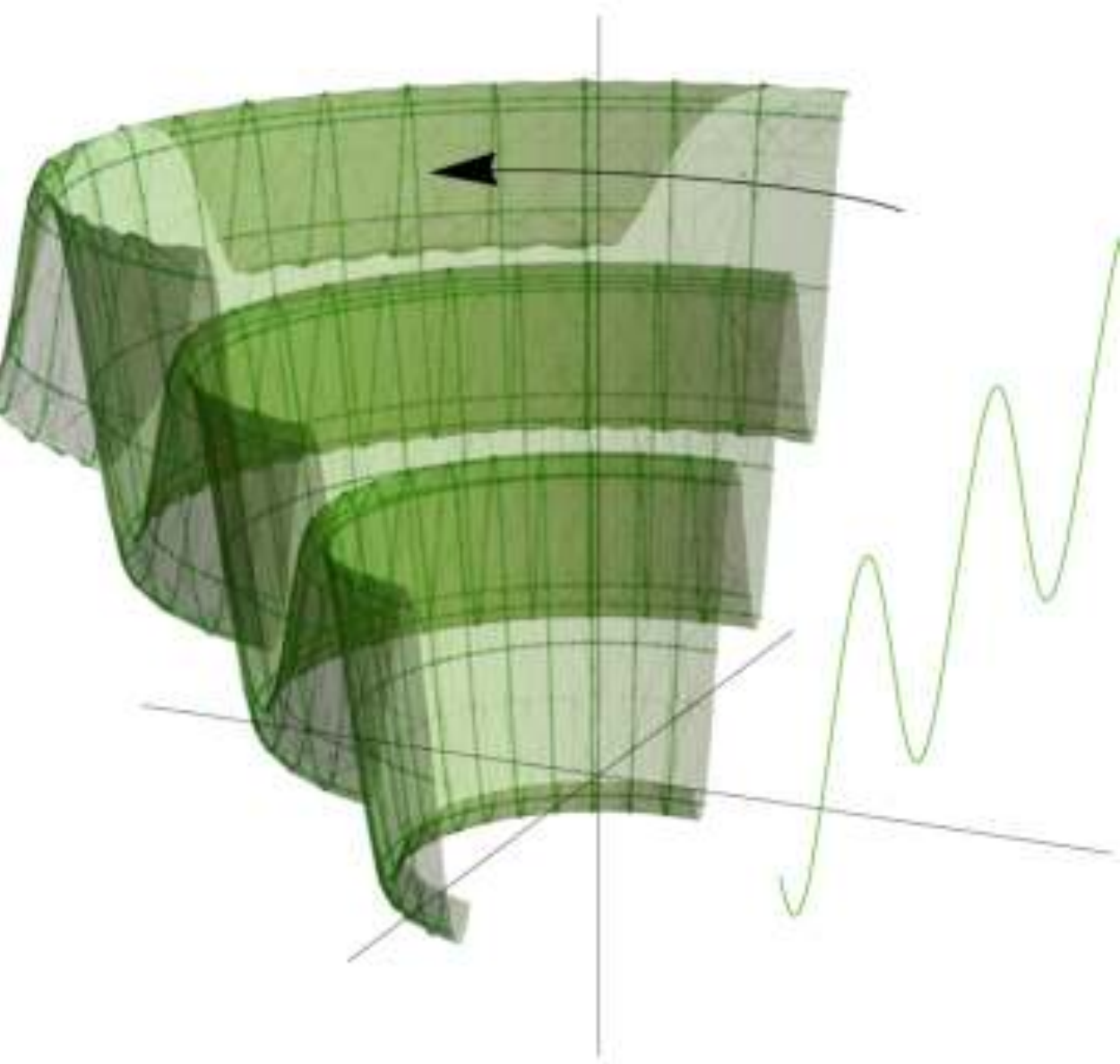}
      \caption{Spanning a surface of revolution in $\LM^3$.}
    \end{figure}
  \end{enumerate}
\end{Ex}

We have topological restrictions on the causal character of a surface:

\begin{prop}
  There is no compact regular surface with constant causal character in $\LM^3$.
\end{prop}

\begin{dem}
  Let $M \subseteq \LM^3$ be a compact regular surface. By compactness, both projections \[ M \ni (x,y,z) \mapsto x \in \R \quad\mbox{and}\quad M \ni (x,y,z) \mapsto z \in \R  \]admit critical points in $M$, say, $\vec{p}$ and $\vec{q}$. Then $T_{\vec{p}}M$ is timelike, while $T_{\vec{q}}M$ is spacelike.
\end{dem}

Just like for regular surfaces in $\R^3$, we will say that the restriction of $\pair{\cdot,\cdot}_L$ to the tangent planes of a regular surface $M\subseteq \LM^3$ is its \emph{First Fundamental Form}\index{First Fundamental Form}. If there's a first form, there should be at least a second one too. And as we might recall, for this we needed some orientability condition. So we'll say that a \emph{Gauss map}\index{Gauss map} for a \emph{non-degenerate} regular surface $M\subseteq \LM^3$ is a smooth choice of unit normal vectors along $M$, that is, a smooth map $\vec{N}\colon M \to \LM^3$ such that $\|\vec{N}(\vec{p})\|_L=1$ and $\vec{N}(\vec{p}) \perp T_{\vec{p}}M$, for all $\vec{p} \in M$. For surfaces in $\R^3$, the codomain of a Gauss map is automatically the sphere $\esf^2$, but in $\LM^3$ this depends on the causal character of $M$. Namely, the codomain of $\vec{N}$ is the de Sitter space $\esf^2_1$ if $M$ is timelike, while it is the hyperbolic plane $\H^2$ if $M$ is spacelike with $\vec{N}$ future-directed (or it's reflection through the plane $z=0$ if $\vec{N}$ is past-directed). Moreover, we see that if $M$ has a fixed causal character, then $\epsilon_M = \pair{\vec{N}(\vec{p}), \vec{N}(\vec{p})}_L$. This is useful for keeping track of the correct signs for some formulas we'll soon deduce.

To understand how a non-degenerate surface $M$ bends in space $\LM^3$ near a point $\vec{p} \in M$, we may focus on the ``linear approximation'' to $M$ at $\vec{p}$: the tangent plane $T_{\vec{p}}M$. Understanding how the tangent planes change near $\vec{p}$ is the same as understanding how their orthogonal complements $\vec{N}(\vec{p})$ change. The motto \[\mbox{``rate of change = derivative''}\] leads to the:

\begin{defn}
  Let $M\subseteq \LM^3$ be a non-degenerate regular surface, and $\vec{N}$ be a Gauss map for $M$. The \emph{Weingarten operator}\index{Weingarten operator} for $M$ at $\vec{p}$ is the differential $-\d{\vec{N}}_{\vec{p}}\colon T_{\vec{p}}M \to T_{\vec{p}}M$. The \emph{Second Fundamental Form}\index{Second Fundamental Form} of $M$ at $\vec{p}$ is the bilinear map $\II_{\vec{p}}\colon T_{\vec{p}}M\times T_{\vec{p}}M \to (T_{\vec{p}}M)^\perp$ characterized by the relation $\pair{\II_{\vec{p}}(\vec{v},\vec{w}), \vec{N}(\vec{p})}_L = \pair{-\d{\vec{N}}_{\vec{p}}(\vec{v}), \vec{w}}_L$, for all $\vec{v},\vec{w} \in T_{\vec{p}}M$. Its \emph{scalar} version $\widetilde{\II}_{\vec{p}}$ is just this common quantity, that is, $\widetilde{\II}_{\vec{p}}(\vec{v},\vec{w}) = \pair{\II_{\vec{p}}(\vec{v},\vec{w}), \vec{N}(\vec{p})}_L$.
\end{defn}

\begin{obs}
  Note that if $M$ is spacelike, then $T_{\vec{p}}M \cong T_{\vec{N}(\vec{p})}(\H^2)$, since both planes are the Lorentz-orthogonal complement of $\vec{N}(\vec{p})$. Similarly, is $M$ is timelike, for the same reason we have $T_{\vec{p}}M \cong T_{\vec{N}(\vec{p})}(\esf^2_1)$, and this is why we may regard $-\d{\vec{N}}_{\vec{p}}$ as a linear operator in $T_{\vec{p}}M$. The negative sign, by the way, is meant to reduce signs in further formulas, is \emph{not} related to the ambient $\LM^3$, and appears naturally in the context of submanifold theory in pseudo-Riemannian geometry, in general.
\end{obs}

One can prove, just like in $\R^3$, that $\d{\vec{N}}_{\vec{p}}$ is a self-adjoint operator with respect to $\pair{\cdot,\cdot}_L$, so that both $\II_{\vec{p}}$ and $\widetilde{\II}_{\vec{p}}$ are symmetric. We will conclude this preliminary discussion by giving precise definitions of ``curvature'' and formulas for expressing them in terms of a parametrization of the surface.

\begin{defn}
  Let $M\subseteq \LM^3$ be a non-degenerate regular surface. The \emph{mean curvature vector}\index{Mean curvature!vector} and the \emph{Gaussian curvature}\index{Gaussian curvature} of $M$ at a point $\vec{p}\in M$ are defined by
  \begin{align*}
    \vec{H}(\vec{p}) &\doteq  \frac{1}{2}\tr_{\pair{\cdot,\cdot}_L}(\II_{\vec{p}}) = \frac{1}{2}(\epsilon_{\vec{v}_1}\II_{\vec{p}}(\vec{v}_1,\vec{v}_1) + \epsilon_{\vec{v}_2}\II_{\vec{p}}(\vec{v}_2,\vec{v}_2)) \quad \mbox{and} \\ K(\vec{p}) &\doteq - {\rm det}_{\pair{\cdot,\cdot}_L}(\widetilde{\II}_{\vec{p}}) =  -\det\big((\II_{\vec{p}}(\vec{v}_i,\vec{v}_j))_{i,j=1}^2\big),
  \end{align*}where $(\vec{v}_1,\vec{v}_2)$ is any orthonormal basis for $T_{\vec{p}}M$.
\end{defn}
 
\begin{obs}\mbox{}
  \begin{itemize}
  \item  The negative sign in the definition of $K$ accounts for the loss of information we have when considering $\widetilde{\II}$ instead of $\II$ there.
  \item  If we write $\vec{H}(\vec{p}) = H(\vec{p})\vec{N}(\vec{p})$, $H(\vec{p})$ is called the \emph{mean curvature}\index{Mean curvature} of $M$ at $\vec{p}$. Choosing the opposite Gauss map changes the sign of $H$, but not of $\vec{H}$.
  \end{itemize}
  
\end{obs}

Still assuming this whole setup, we recall the classical notation for the coefficients of the fundamental forms. If $\vec{x}\colon U \to \vec{x}[U]\subseteq M$ is a parametrization, then we set \[ E \doteq \left\langle \frac{\partial\vec{x}}{\partial u},\frac{\partial \vec{x}}{\partial u}\right\rangle_L, \quad F \doteq \left\langle \frac{\partial\vec{x}}{\partial u},\frac{\partial \vec{x}}{\partial v}\right\rangle_L \quad\mbox{and}\quad G \doteq \left\langle \frac{\partial\vec{x}}{\partial v},\frac{\partial \vec{x}}{\partial v}\right\rangle_L,    \]as well as \[e \doteq \left\langle \frac{\partial^2\vec{x}}{\partial u^2},\vec{N}\circ \vec{x}\right\rangle_L, \quad f \doteq \left\langle \frac{\partial^2\vec{x}}{\partial u\partial v},\vec{N}\circ\vec{x}\right\rangle_L \quad\mbox{and}\quad g \doteq \left\langle \frac{\partial^2\vec{x}}{\partial v^2},\vec{N}\circ\vec{x}\right\rangle_L,\]so that (with a mild abuse of notation) we have \[ \II\left(\frac{\partial\vec{x}}{\partial u}\right) = \epsilon_M e \vec{N}, \quad \II\left(\frac{\partial\vec{x}}{\partial u},\frac{\partial \vec{x}}{\partial v}\right) = \epsilon_M f \vec{N}, \quad\mbox{and}\quad \II\left(\frac{\partial\vec{x}}{\partial v}\right) = \epsilon_M g \vec{N}.  \] 

To produce orthonormal bases for the tangent planes to $M$, needed for computing $\vec{H}$ and $K$ via the definitions, the Gram-Schmidt process comes to rescue. We obtain similar formulas for the ones in $\R^3$, which now take into account the causal character of $M$ itself:

\begin{prop}\label{prop:coords_KH}
  Let $M\subseteq \LM^3$ be a non-degenerate regular surface, and $\vec{x}\colon U \to \vec{x}[U]\subseteq M$ a parametrization for $M$. Then \[  H\circ\vec{x} = \frac{\epsilon_M}{2} \frac{Eg+eG-2Ff}{EG-F^2}\quad\mbox{and}\quad K\circ \vec{x} = \epsilon_M\frac{eg-f^2}{EG-F^2}.  \]
\end{prop}

Do note that setting $\epsilon_M=1$ if $M\subseteq \R^3$, the above gives also correct results for the mean and Gaussian curvatures of $M$. The details of these maybe-not-so-short calculations may be consulted, for example, in~\cite{TL}. They also follow from the more general theory developed in~\cite{ON2}. Here are some more examples:
\newpage
\begin{Ex}\mbox{}
  \begin{enumerate}[(1)]
  \item Planes admit a constant Gauss map, so the Weingarten operator vanishes. Hence we get $H=K=0$.
  \item The position map $\vec{N}(\vec{p})=\vec{p}$ is a Gauss map for both the de Sitter space $\esf^2_1$ and the hyperbolic plane $\H^2$. Taking the causal characters into account, we obtain $K=1$ and $H = -1$ for $\esf^2_1$, and $K=-1$ and $H=1$ for $\H^2$. If you studied anything about hyperbolic geometry before, this serves both as a quick sanity check (hyperbolic plane should have negative curvature) as well as another justification for the presence of the minus sign in the definition of $K$.
  \item If $f\colon U \subseteq \R^2 \to \R$ is a smooth function for which the graph ${\rm gr}(f)\subseteq \LM^3$ is non-degenerate, by applying the coordinate formulas given in Proposition \ref{prop:coords_KH} (p. \pageref{prop:coords_KH}), we obtain \[ K = \frac{f_{uv}^2 - f_{uu}f_{vv}}{(-1+f_u^2+f_v^2)^2} \quad\mbox{and}\quad H =\frac{f_{uu}(-1+f_v^2) - 2f_uf_vf_{uv}+f_{vv}(-1+f_u^2)}{|-1+f_u^2+f_v^2|^{3/2}}.  \]
  \end{enumerate}
\end{Ex}

\subsection{The Diagonalization Problem}

We know from linear algebra the \emph{Real Spectral Theorem}\index{Real Spectral Theorem}: that if $(V,\pair{\cdot,\cdot})$ is a finite-dimensional real vector space equipped with a positive-definite inner product, and $T\colon V \to V$ is a linear operator which is self-adjoint with respect to $\pair{\cdot,\cdot}$, then $V$ admits an orthonormal basis of eigenvectors of $T$. This result is no longer true if $\pair{\cdot,\cdot}$ is not positive-definite, and non-degeneracy alone is not strong enough to ensure any good conclusions. There is one adaptation, though: if $\dim V \geq 3$ and $\pair{T(\vec{v}),\vec{v}} \neq 0$ for all non-zero $\vec{v} \in V$ with $\pair{\vec{v},\vec{v}} = 0$, then $V$ admits an orthonormal basis of eigenvectors of $T$. A very surprising proof using integration and homotopy, due to Milnor, may be found in~\cite{Greub}.

We have seen that the Weingarten operator of any non-degenerate surface $M\subseteq \LM^3$ is still self-adjoint with respect to the First Fundamental Form of $M$. So we conclude that if $M$ is spacelike, then $-\d{\vec{N}}_{\vec{p}}$ is diagonalizable: the eigenvalues $\kappa_1(\vec{p})$ and $\kappa_2(\vec{p})$ are called the \emph{principal curvatures}\index{Principal!curvatures} of $M$ at $\vec{p}$, and the (orthogonal) eigenvectors are called the \emph{principal directions}\index{Principal!directions} of $M$ at $\vec{p}$. We cannot guarantee the existence of principal directions for timelike surfaces in $M$, even with the sharpened version of the Spectral Theorem mentioned above, since $\dim T_{\vec{p}}M = 2 < 3$. 

That being said, our goal here is to understand precisely when do we have principal directions for timelike surfaces in $\LM^3$.

\begin{prop}
  Let $M\subseteq \LM^3$ be a non-degenerate regular surface with diagonalizable Weingarten operators. Then
  \[ H(\vec{p}) = \epsilon_M
    \frac{\kappa_1(\vec{p})+\kappa_2(\vec{p})}{2} \quad\mbox{and}\quad
    K(\vec{p}) = \epsilon_M \kappa_1(\vec{p})\kappa_2(\vec{p}).  \]
\end{prop}

\begin{obs}
  Usually one defines $H$ and $K$ for surfaces in $\R^3$ by the above formulas (setting $\epsilon_M=1$, of course). The reason why we went through the hassle of using metric traces and determinants to define them in $\LM^3$ was just so we could have a unified approach that worked in all the cases simultaneously, even when we could not use principal curvatures. Also note that the expression for $H$ justifies the name ``mean'' curvature.
\end{obs}

We might as well start understanding a class of surfaces which, in general, have diagonalizable Weingarten operators.

\begin{defn}
  Let $M\subseteq \LM^3$ be a non-degenerate regular surface, and $\vec{p} \in M$. The point $\vec{p}$ is called \emph{umbilic}\index{Umbilic point} if there is $\lambda(\vec{p}) \in \R$ such that
  \[\widetilde{\vphantom{E}\II}_{\vec{p}}(\vec{v},\vec{w}) =
    \lambda(\vec{p})\pair{\vec{v},\vec{w}},\] for all $\vec{v},\vec{w}\in T_{\vec{p}}M$. We will also say that $M$ is
  \emph{totally umbilic}\index{Totally umbilic surface} if all its points are umbilic.
\end{defn}

Informally, a point is umbilic if there the two fundamental forms of $M$ are ``linearly dependent''. In umbilical points, we have $-\d{\vec{N}}_{\vec{p}} = \lambda(\vec{p}){\rm Id}_{T_{\vec{p}}M}$. Indeed, for all vectors $\vec{v},\vec{w} \in T_{\vec{p}}M$ we have that $\pair{\lambda(\vec{p})\vec{v},\vec{w}}_L = \widetilde{\II}_{\vec{p}}(\vec{v},\vec{w}) = \pair{-\d{\vec{N}}_{\vec{p}}(\vec{v}),\vec{w}}_L$, and the conclusion follows from non-degeneracy of $\pair{\cdot,\cdot}_L$ restricted to $T_{\vec{p}}M$.

You might remember from the classical theory in $\R^3$ that there, the only totally umbilic surfaces are spheres and planes. Since the de Sitter space $\esf^2_1$ and the hyperbolic plane $\H^2$ (together with its reflection $\H^2_-$) play the role of spheres in $\LM^3$, the following result (with the same proof as in $\R^3$) should not be a surprise:

\begin{teo}[Characterization of totally umbilic surfaces in $\LM^3$]\label{teo:umb}
  Let $M\subseteq \LM^3$ be a non-degenerate, regular, connected and totally umbilic surface. Then $M$ is contained in some plane, or there is a center $\vec{c} \in \R^3_\nu$ and a radius $r>0$ such that
  \begin{enumerate}[(i)]
    \item if $M$ is spacelike, then $M\subseteq
    \H^2(\vec{c},r)$ or $M \subseteq \H^2_-(\vec{c},r)$;
    \item if $M$ is timelike, then $M\subseteq \esf^2_1(\vec{c},r)$.
  \end{enumerate}
\end{teo}

\begin{obs}
Here, we mean $\esf^2_1(\vec{c},r) = \{\vec{p} \in \LM^3 \mid \pair{\vec{p}-\vec{c},\vec{p}-\vec{c}}_L=r^2\}$, etc.. Moreover, in the timelike case, what decides between $\H^2(\vec{c},r)$ or $\H^2_-(\vec{c},r)$ is the direction of the timelike vector $\vec{p}-\vec{c}$ for some (hence all) $\vec{p} \in M$ (due to connectedness).
\end{obs}

\begin{figure}[H]
  \centering
  \includegraphics[width =.4\textwidth]{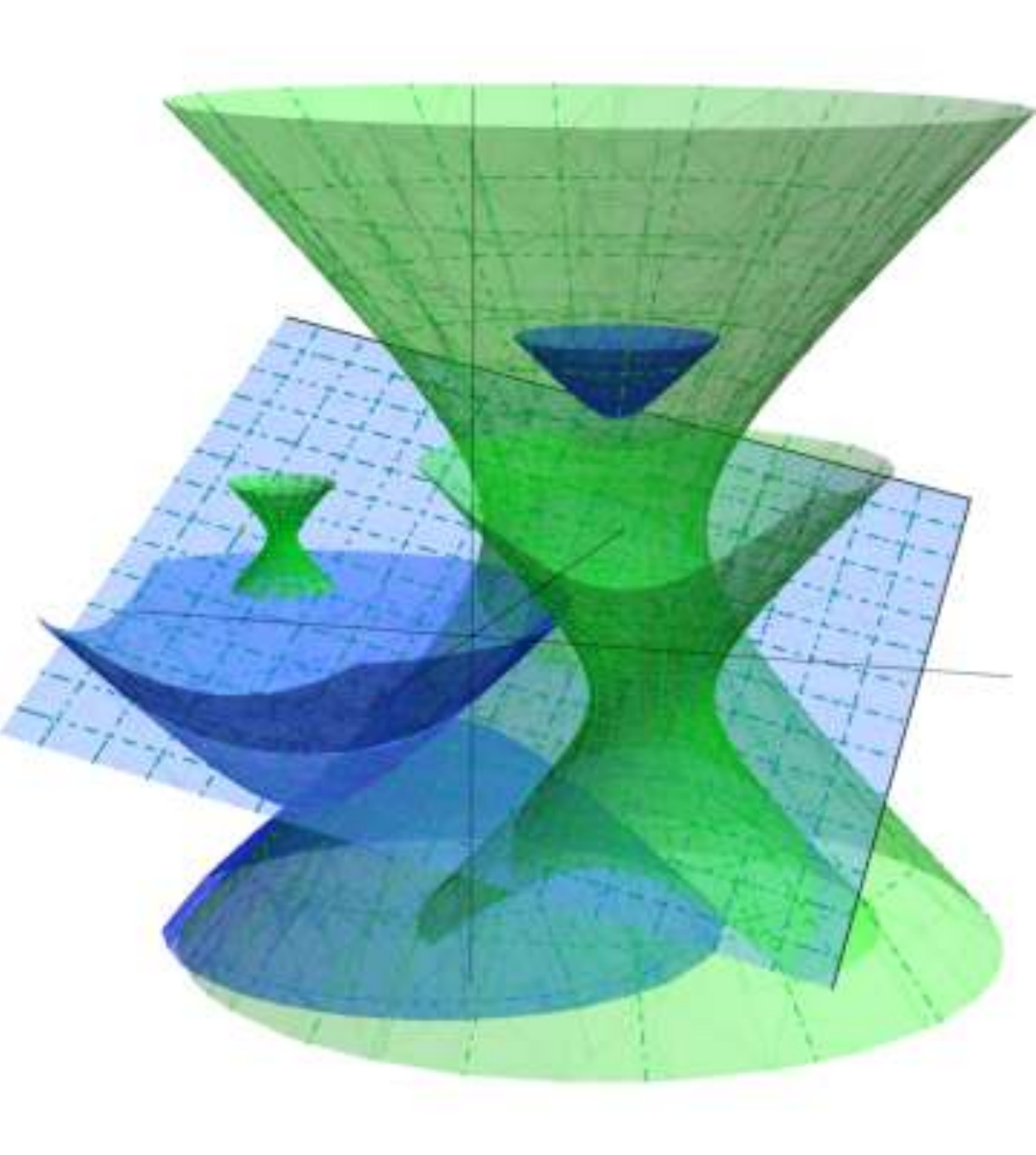}
  \caption{The totally umbilic surfaces in $\LM^3$.}
\end{figure}

Back to the diagonalization problem. Let's see necessary conditions for an affirmative answer to the problem.

\begin{prop}
  Let $M\subseteq \LM^3$ be a non-degenerate regular surface, and $\vec{p}\in M$ such that the Weingarten operator at $\vec{p}$ is diagonalizable. Then $H(\vec{p})^2 - \epsilon_M K(\vec{p}) \geq 0$, with equality holding if and only if $\vec{p}$ is umbilic.
\end{prop}

\begin{dem}
  Directly, we have: \begin{align*} 0 &\leq \left(\frac{\kappa_1(\vec{p})-\kappa_2(\vec{p})}{2}\right)^2 = \frac{\kappa_1(\vec{p})^2-2\kappa_1(\vec{p})\kappa_2(\vec{p}) + \kappa_2(\vec{p})^2}{4} \\ &= \frac{\kappa_1(\vec{p})^2+2\kappa_1(\vec{p})\kappa_2(\vec{p}) + \kappa_2(\vec{p})^2}{4} - \kappa_1(\vec{p})\kappa_2(\vec{p}) \\ &=\left(\frac{\kappa_1(\vec{p})+\kappa_2(\vec{p})}{2}\right)^2 - \kappa_1(\vec{p})\kappa_2(\vec{p}) \\ &= (\epsilon_M H(\vec{p}))^2 - \epsilon_M
    K(\vec{p}) = H(\vec{p})^2 - \epsilon_M K(\vec{p}).
  \end{align*}Equality holds if and only if $\kappa_1(\vec{p}) =
  \kappa_2(\vec{p})$, that is to say, if $\vec{p}$ is umbilic.
\end{dem}

So we have a necessary, but not sufficient condition for the diagonalizability of the Weingarten operators. What we can see, though, is that the quantity $H(\vec{p})^2 - \epsilon_M K(\vec{p})$ will play a big role in our analysis, which will be done in full detail in the proof of the desired:

\begin{teo}[Diagonalization in $\LM^3$]\label{teo:diagL3}
    Let $M\subseteq \LM^3$ be a non-degenerate regular surface, $\vec{N}$ a Gauss map for $M$, and $\vec{p}\in M$. Then:

  \begin{enumerate}[(i)]
    \item if $H(\vec{p})^2 - \epsilon_M K(\vec{p}) > 0$, $-{\dd}\vec{N}_{\vec{p}}$ is diagonalizable;
    \item if $H(\vec{p})^2 - \epsilon_M K(\vec{p}) < 0$, $-{\dd}\vec{N}_{\vec{p}}$ is not diagonalizable;
    \item if $H(\vec{p})^2 - \epsilon_M K(\vec{p}) = 0$ and $M$ is spacelike, then $\vec{p}$ is umbilic, and hence $-{\dd}\vec{N}_{\vec{p}}$ is diagonalizable.
  \end{enumerate}
\end{teo}

\begin{obs}
  If $H(\vec{p})^2 - \epsilon_M K(\vec{p}) = 0$ and $M$ is timelike, the criterion is inconclusive and the Weingarten operator may or may not be diagonalizable.
\end{obs}

\begin{dem}
  Consider the characteristic polynomial $c(t)$ of $-{\dd}\vec{N}_{\vec{p}}$, given by \[ c(t) = t^2 - {\rm tr}(-{\dd}\vec{N}_{\vec{p}})\, t + \det(-{\dd}\vec{N}_{\vec{p}}) = t^2 - 2\epsilon_M H({\vec{p}}) t + \epsilon_M
  K({\vec{p}}), \]whose discriminant is: \[ (-2\epsilon_M
  H({\vec{p}}))^2 - 4 (\epsilon_M K({\vec{p}})) = 4(H({\vec{p}})^2 -
  \epsilon_M K({\vec{p}})). \]
  \begin{itemize}
  \item If $H({\vec{p}})^2 - \epsilon_M K({\vec{p}}) > 0$, then $c(t)$ has two distinct roots, which are the eigenvalues of $-{\dd}\vec{N}_{\vec{p}}$, who then admits two linearly independent eigenvectors (hence diagonalizable).
  \item If $H({\vec{p}})^2 - \epsilon_M K({\vec{p}}) < 0$, $c(t)$ does not have any real roots. Thus $-{\dd}\vec{N}_{\vec{p}}$ has no real eigenvalues, and hence it is not diagonalizable.
  \item Now assume that $H({\vec{p}})^2 - \epsilon_M K({\vec{p}}) = 0$ and that
$M$ is spacelike, that is, that $K(\vec{p}) = -H({\vec{p}})^2$. From the expression given for the discriminant of $c(t)$, it follows that $-H({\vec{p}})$ is an eigenvalue of $-{\dd}\vec{N}_{\vec{p}}$. So, there is a unit (spacelike) vector $\vec{u}_1 \in T_{\vec{p}}M$ such that ${\dd}\vec{N}_{\vec{p}}(\vec{u}_1) =
H({\vec{p}})\vec{u}_1$. Consider then an orthogonal basis $\mathcal{B} \doteq (\vec{u}_1,\vec{u}_2)$ of $T_{\vec{p}}M$.
Then:\[ \left[{\dd}\vec{N}_{\vec{p}}\right]_{\mathcal{B}}
  = \begin{pmatrix} H({\vec{p}}) & a \\ 0 & b
  \end{pmatrix}, \quad \text{ where } {\dd}\vec{N}_{\vec{p}}(\vec{u}_2) =
  a\vec{u}_1 + b\vec{u}_2.  \]It suffices to check that $a = 0$ and $b = H(\vec{p})$ to conclude the proof. Applying $\pair{\cdot, \vec{u}_1}_L$, we have: \[ a = \pair{{\dd}\vec{N}_{\vec{p}}(\vec{u}_2),\vec{u}_1}_L = \pair{\vec{u}_2, {\dd}\vec{N}_{\vec{p}}(\vec{u}_1)}_L =
  \pair{\vec{u}_2, H({\vec{p}})\vec{u}_1}_L = H({\vec{p}})\pair{\vec{u}_2,\vec{u}_1}_L = 0.  \]On the other hand: \[ -H({\vec{p}})^2 =
  K({\vec{p}}) = -\det(-{\dd}\vec{N}_{\vec{p}}) = -\det({\dd}\vec{N}_{\vec{p}}) =
  -H({\vec{p}})b,\]so that $H({\vec{p}})b= H({\vec{p}})^2$. If $H({\vec{p}}) =
  0$, then ${\dd}\vec{N}_{\vec{p}}$ is the zero map (hence diagonalizable). If $H({\vec{p}}) \neq 0$, we obtain $b = H({\vec{p}})$, as wanted. Note that in this case $\vec{p}$ is umbilic.
  \end{itemize}
\end{dem}

Observe that in the above proof, we would not be able to control the causal type of the eigenvector $\vec{u}_1$ in the last case discussed if $M$ were timelike. If $\vec{u}_1$ were lightlike, we could not consider the basis
$\mathcal{B}$ to proceed with the argument. With this in mind, we obtain the following extension of the theorem:

\begin{cor}\label{cor:diag_tempo}
  Let $M\subseteq \LM^3$ be a timelike regular surface and $\vec{p}\in M$ be a point with $H({\vec{p}})^2 - K({\vec{p}}) = 0$. If
  $-{\dd}\vec{N}_{\vec{p}}$ has no lightlike eigenvectors, then it is diagonalizable and $\vec{p}$ is umbilic, with both principal curvatures equal to $-H({\vec{p}})$.
\end{cor}

Let's conclude the section exploring examples of timelike surfaces for which the equality $H({\vec{p}})^2 = K({\vec{p}})$ holds and anything can happen with the Weingarten operators.

\begin{Ex}\label{Ex:b_scroll}\mbox{}
  \begin{enumerate}[(1)]
    \item For the de Sitter space $\mathbb{S}_1^2$,
    we had $-{\dd}\vec{N}_{\vec{p}} =
    -{\rm Id}_{T_{\vec{p}}\left(\mathbb{S}^2_1\right)}$ (hence diagonalizable), with $K=1$ and $H = -1$, so that $H^2 - K = 0$.
  \item Consider a lightlike curve $\vec{\alpha}\colon I \to \LM^3$ with arc-photon parameter. Define the \emph{$\vec{B}$-scroll}\index{BScroll@$\vec{B}$-scroll (associated to a lightlike curve)} associated to $\vec{\alpha}$, $\vec{x}\colon I \times \R \to \LM^3$ given by $\vec{x}(\phi,t) \doteq \vec{\alpha}(\phi)+t \sub{B}{\alpha}(\phi)$. Restricting enough the domain of $\vec{x}$, we may assume that its image $M$ is a regular surface. Put, for each $\phi\in I$, $D(\phi) \doteq \det\big(\sub{T}{\alpha}(\phi), \sub{N}{\alpha}(\phi) ,\sub{B}{\alpha}(\phi)\big)>0$. Computing the derivatives \[ \vec{x}_\phi(\phi,t) = \sub{T}{\alpha}(\phi) + t \subg{\ptau}{\alpha}(\phi) \sub{N}{\alpha}(\phi) \quad\mbox{and}\quad \vec{x}_t(\phi,t) = \sub{B}{\alpha}(\phi), \]we immediately have that\[ (g_{ij}(\phi,t))_{1 \leq i,j \leq 2} =
\begin{pmatrix}
  t^2 \subg{\ptau}{\alpha}(\phi)^2 & -1  \\ -1 & 0
\end{pmatrix},
\]whence $M$ is timelike. Here, $g_{ij}$ is shorthand for the coefficients of the First Fundamental Form. Noting that
$|\det((g_{ij}(\phi,t))_{1 \leq i,j \leq 2})| = 1$, we directly obtain that
\[ \vec{N}(\vec{x}(\phi,t)) = \sub{T}{\alpha}(\phi)\times_L
  \sub{B}{\alpha}(\phi) +t\subg{\ptau}{\alpha}(\phi)
  \sub{N}{\alpha}(\phi)\times_L \sub{B}{\alpha}(\phi). \]Computing the second order derivatives
\begin{align*}
  \vec{x}_{\phi\phi}(\phi,t) &= t \subg{\ptau}{\alpha}(\phi)^2 \sub{T}{\alpha}(\phi) + (1+t\subg{\ptau}{\alpha}'(\phi)) \sub{N}{\alpha}(\phi) + t\subg{\ptau}{\alpha}(\phi)\sub{B}{\alpha}(\phi), \\ \vec{x}_{\phi t}(\phi,t) &= \subg{\ptau}{\alpha}(\phi) \sub{N}{\alpha}(\phi) \quad\mbox{and}\\ \vec{x}_{tt}(\phi,t) &= \vec{0},
\end{align*}we obtain the coefficients $h_{ij}$ of the Second Fundamental Form: \[ (h_{ij}(\phi,t))_{1 \leq i,j \leq 2} =
\begin{pmatrix}
  (-1-t\subg{\ptau}{\alpha}'(\phi)+t^2\subg{\ptau}{\alpha}(\phi)^3)D(\phi)
  & -\subg{\ptau}{\alpha}(\phi) D(\phi) \\ -\subg{\ptau}{\alpha}(\phi)
  D(\phi) & 0
\end{pmatrix}. \]It follows that
\[ K(\vec{x}(\phi,t)) = \subg{\ptau}{\alpha}(\phi)^2
  D(\phi)^2\quad\mbox{and}\quad H(\vec{x}(\phi,t)) =
  \subg{\ptau}{\alpha}(\phi) D(\phi).\]We then know that, in each point
$\vec{x}(\phi,t)$, $-{\dd}\vec{N}_{\vec{x}(\phi,t)}$ has only one eigenvalue (namely, $\subg{\ptau}{\alpha}(\phi) D(\phi)$). It suffices to check then that there are points in $M$ for which the associated eigenspace has dimension $1$ -- this shows that the Weingarten operators at those points are not diagonalizable. To wit, we have
\[ \left[-{\dd}\vec{N}_{\vec{x}(\phi,t)}\right]_{\mathcal{B}_{\vec{x}}} =
  D(\phi)\begin{pmatrix} \subg{\ptau}{\alpha}(\phi) & 0 \\
    1+t\subg{\ptau}{\alpha}(\phi) & \subg{\ptau}{\alpha}(\phi)
\end{pmatrix},
\]and the kernel of \[\begin{pmatrix}
   0 & 0 \\ 1+t\subg{\ptau}{\alpha}(\phi) & 0
\end{pmatrix}\] has always dimension $1$ when $1+t\subg{\ptau}{\alpha}(\phi) \neq 0$ (e.g., along $\vec{\alpha}$ itself, setting $t=0$).
  \end{enumerate}
\end{Ex}

\newpage
\section*{Problems}\addcontentsline{toc}{subsection}{Problems for the last section}

\begin{problem}
  Work through Example \ref{Ex:tipos_causais} (p. \pageref{Ex:tipos_causais}).
\end{problem}

\begin{problem}[Horocycles] Let $\vec{v} \in \LM^3$ be a future-directed lightlike vector, and $c<0$. The set $H_{\vec{v},c} \doteq \{\vec{x} \in \H^2 \mid \pair{\vec{x},\vec{v}}_L=c\}$ is called a \emph{horocycle\index{Horocycle} of $\H^2$ based on $\vec{v}$}. Let $\vec{\alpha}\colon I \to H_{\vec{v},c}$ be a unit speed curve (assume that $0 \in I$, reparametrizing if necessary).
  \begin{enumerate}[(a)]
  \item Show that \[  \vec{\alpha}(s) = -\frac{s^2}{2c}\vec{v} + s\vec{w}_1+\vec{w}_2, \]for some unit and orthogonal vectors $\vec{w}_1$ and $\vec{w}_2$, with $\vec{w}_1$ spacelike, $\vec{w}_2$ timelike, $\pair{\vec{w}_1,\vec{v}}_L=0$ and $\pair{\vec{w}_2,\vec{v}}_L=c$.
  \item Conclude that $\vec{\alpha}$ is a semi-lightlike curve whose pseudo-torsion identically vanishes.
  \end{enumerate}
\end{problem}

\begin{problem}
  Compute the Gaussian and mean curvatures for the surface of revolution spanned by a unit speed curve as in item (3) of Example \ref{Ex:tipos_causais} (p. \pageref{Ex:tipos_causais}).
  \begin{obs}
    One can also study surfaces of revolution in $\LM^3$ generated by hyperbolic rotations about the $x$-axis instead of the timelike $z$-axis. See~\cite{TL} for more about this.
  \end{obs}
\end{problem}

\begin{problem}
  Let $\vec{\alpha}\colon I \to \LM^3$ and $\vec{\beta}\colon J \to \LM^3$ be two smooth lightlike curves such that $\{\vec{\alpha}'(u),\vec{\beta}'(v)\}$ is linearly independent for all $(u,v) \in I\times J$. Then, reducing $I$ and $J$ if necessary, the image $M$ of the sum $\vec{x}\colon I\times J \to \LM^3$ given by $\vec{x}(u,v) = \vec{\alpha}(u)+\vec{\beta}(v)$ is a regular surface. Show that $M$ is timelike with $H=0$.
  \begin{obs}
    Actually, the ``converse'' holds: every timelike surface with $H=0$ admits parametrizations like this $\vec{x}$ above. See~\cite{Chen} for more details.
  \end{obs}
\end{problem}

\begin{problem}
  Prove Theorem \ref{teo:umb} (p. \pageref{teo:umb}).
\end{problem}

\begin{problem}
  Let $M\subseteq \LM^3$ be a non-degenerate regular surface, and $\vec{N}\colon M \to \LM^3$ a Gauss map for $M$. Show that the Weingarten operator $-\d{\vec{N}}_{\vec{p}}$ is self-adjoint with respect to $\pair{\cdot,\cdot}_L$, for all $\vec{p} \in M$. Namely, show that given $\vec{v},\vec{w} \in T_{\vec{p}}M$, we have \[  \pair{\d{\vec{N}}_{\vec{p}}(\vec{v}),\vec{w}}_L = \pair{\vec{v},\d{\vec{N}}_{\vec{p}}(\vec{w})}_L. \]
  \begin{dica}
    Use a parametrization of $M$ and do it locally.
  \end{dica}
\end{problem}

\begin{problem}
  Make sure you understand how to obtain Corollary \ref{cor:diag_tempo} (p. \pageref{cor:diag_tempo}) by adapting the proof of Theorem \ref{teo:diagL3} (p. \pageref{teo:diagL3}).
\end{problem}

\begin{problem}
  Consider the \emph{anti-de Sitter}\index{Anti-de Sitter space} space $\H^2_1 \doteq \{ \vec{p} \in \R^3_2 \mid \pair{\vec{p},\vec{p}}_2 = -1 \}$. Try to understand how to translate the results discussed for the ambient $\LM^3$ for the ambient $\R^3_2$ and show that $\H^2_1$ has constant Gaussian curvature $K=-1$. 
\end{problem}

\hrulefill

\newpage
\section*{Extra \#1: Riemann's classification of surfaces with constant $K$}\addcontentsline{toc}{section}{Extra \#1: Riemann's classification of surfaces with constant $K$}

Up to this moment, we know some surfaces with constant Gaussian curvature. Namely, we have met: 

\begin{itemize}
\item The planes $\R^2$ and $\LM^2$, with $K=0$;
\item The sphere $\esf^2$ and the de Sitter space $\esf^2_1$, with $K=1$;
\item The hyperbolic plane $\H^2$ and the anti-de Sitter $\H^2_1$, with $K=-1$.
\end{itemize}

Our goal here is to show that, locally, every surface with constant $K$ is one of those surfaces described above. More precisely, we want to prove the:

\begin{teo}[Riemann]\label{teo:class_riemann}\index{Riemann's!classification theorem}
  Let $(M,\pair{\cdot,\cdot})$ be a geometric surface with constant Gaussian curvature $K \in \{-1,0,1\}$. Then:

  \begin{enumerate}[(A)]
  \item if the metric is Riemannian, every point in $M$ has a neighborhood isometric to an open subset of
  \begin{enumerate}[(i)]
  \item $\R^2$, if $K=0$;
  \item $\esf^2$, if $K=1$;
  \item $\H^2$, if $K=-1$,
  \end{enumerate}
\item while if the metric is Lorentzian, to an open subset of
  \begin{enumerate}[(i)]
  \item $\LM^2$, if $K=0$;
  \item $\esf^2_1$, if $K=1$;
  \item $\H^2_1$, if $K=-1$.
  \end{enumerate}
  \end{enumerate}
\end{teo}

By \emph{geometric surface}\index{Geometric surface}, we mean an abstract surface ($2$-dimensional manifold) endowed with a metric tensor (called \emph{Riemannian} if positive-definite, or \emph{Lorentzian} if it has index $1$). The proof strategy consists in constructing parametrizations for which the metric assumes a simple form. To actually do this, we will use \emph{geodesics}\index{Geodesic}, which are know to be plentiful in any geometric surface. 

Recall here that given any regular parametrization $\vec{x}\colon U\subseteq \R^2 \to \vec{x}[U] \subseteq M$ of our surface, we set $g_{ij} = \pair{\vec{x}_u,\vec{x}_v}$, so that $(g^{ij})_{i,j=1}^2$ is the inverse matrix of $(g_{ij})_{i,j=1}^2$, and the \emph{Christoffel symbols}\index{Christoffel symbols} of $\vec{x}$ are defined by \[  \Gamma_{ij}^k = \sum_{r=1}^2 \frac{g^{kr}}{2} \left(\frac{\partial g_{ir}}{\partial u^j} + \frac{\partial g_{jr}}{\partial u^i} - \frac{\partial g_{ij}}{\partial u^r} \right), \]where we identify $u \leftrightarrow u^1$, $v \leftrightarrow u^2$, and $i,j,k \in \{1,2\}$. Geodesics are curves $\vec{\gamma}\colon I \to M$ with the property that given any parametrization $\vec{x}$ and writing $\vec{\gamma}(t) = \vec{x}(u(t),v(t))$, we have \[ \ddot{u}^k + \sum_{i,j=1}^k \Gamma_{ij}^k \dot{u}^i\dot{u}^j = 0, \qquad k=1,2.  \]
Further general facts about geodesics (which won't be necessary here) can be consulted in pretty much any book (we list here~\cite{dC1},~\cite{ON2},~\cite{Keti} or~\cite{TL}, for concreteness). To avoid singularities, we won't consider lightlike geodesics.

\emph{Thus, we fix once and for all a geometric surface $(M,\pair{\cdot,\cdot})$, with metric tensor of index $\nu \in \{0,1\}$, and a unit speed geodesic $\vec{\gamma}\colon I \to M$}. For each $v \in I$, consider a unit speed geodesic $\vec{\gamma}_v\colon J_v \to M$, which crosses $\vec{\gamma}$ orthogonally at the point $\vec{\gamma}_v(0) \doteq
\vec{\gamma}(v)$. Setting
\[ U \doteq \{ (u,v) \in \R^2 \mid v \in I \mbox{ e } u \in J_v
  \}, \]define $\vec{x}\colon U \to \vec{x}(U) \subseteq M$ by
$\vec{x}(u,v) = \vec{\gamma}_v(u)$.

   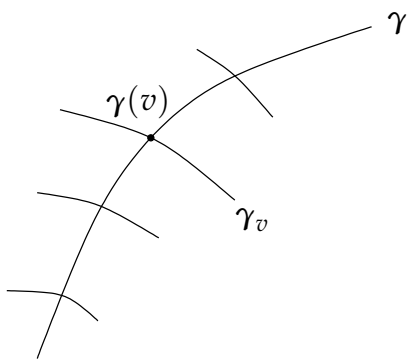
\begin{figure}[H]
      \centering
      \begin{picture}(0,0)%
\includegraphics{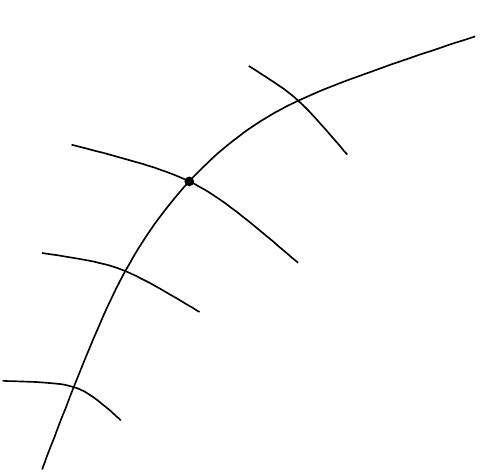}%
\end{picture}%
\setlength{\unitlength}{4144sp}%
\begingroup\makeatletter\ifx\SetFigFont\undefined%
\gdef\SetFigFont#1#2#3#4#5{%
  \reset@font\fontsize{#1}{#2pt}%
  \fontfamily{#3}\fontseries{#4}\fontshape{#5}%
  \selectfont}%
\fi\endgroup%
\begin{picture}(2277,2151)(349,-2053)
\put(1711,-1231){\makebox(0,0)[lb]{\smash{{\SetFigFont{12}{14.4}{\rmdefault}{\mddefault}{\updefault}{\color[rgb]{0,0,0}$\vec{\gamma}_v$}%
}}}}
\put(946,-556){\makebox(0,0)[lb]{\smash{{\SetFigFont{12}{14.4}{\rmdefault}{\mddefault}{\updefault}{\color[rgb]{0,0,0}$\vec{\gamma}(v)$}%
}}}}
\put(2611,-61){\makebox(0,0)[lb]{\smash{{\SetFigFont{12}{14.4}{\rmdefault}{\mddefault}{\updefault}{\color[rgb]{0,0,0}$\vec{\gamma}$}%
}}}}
\end{picture}
      \caption[Construction of a Fermi chart.]{Construction of a Fermi chart $\vec{x}$.}\label{fig:fermi_param}
    \end{figure}

\begin{defn}
  The chart $\vec{x}$ above defined is called a \emph{Fermi chart}\index{Fermi chart} for $M$, centered in $\vec{\gamma}$.
\end{defn}

We'll also fix until the end of the section this Fermi chart $\vec{x}\colon U \to \vec{x}(U)\subseteq M$ so constructed.

\begin{obs}\mbox{}
  \begin{itemize}
  \item  When $\pair{\cdot,\cdot}$ is Lorentzian, we'll have two types of Fermi charts, according to the causal character of $\vec{\gamma}$. Moreover, recalling that geodesics have automatically constant causal character (hence determined by a single velocity vector), it follows that if $\vec{\gamma}$ is spacelike (resp. timelike), then all the $\vec{\gamma}_v$ are timelike (resp. spacelike), since $\{ \vec{\gamma}'(v), \vec{\gamma}_v'(0) \}$ is a orthonormal basis of $T_{\vec{\gamma}(v)}M$, for all $v \in I$.
  \item When necessary, if $\vec{\gamma}$ is timelike, we might denote the coordinates by $(\tau, \vartheta)$ instead of $(u,v)$.    
  \end{itemize}
\end{obs}

\begin{prop}
  The Fermi chart $\vec{x}$ is indeed regular in a neighborhood of $\{0\}\times I$ (so that reducing $U$ if necessary, we may assume that $\vec{x}$ itself is regular).
\end{prop}

\begin{dem}
We'll show that for all $v \in I$, the vectors $\vec{x}_u(0,v)$ and $\vec{x}_v(0,v)$ are orthogonal. To wit, we have by construction that \[ \langle \vec{x}_u(0,v),\vec{x}_v(0,v)\rangle =
 \langle \vec{\gamma}_v'(0), \vec{\gamma}'(v)\rangle = 0. \]Since none of those vectors is lightlike, orthogonality implies linear independence. By continuity of $\vec{x}$, the vectors $\vec{x}_u(u,v)$ and $\vec{x}_v(u,v)$ remain linearly independent for small enough values of $u$.
\end{dem}

\begin{prop}The coordinate expression of $\pair{\cdot,\cdot}$ with respect to the Fermi chart $\vec{x}$ is \[ {\dd}s^2 = (-1)^\nu \epsilon_{\vec{\gamma}}\,{\dd}u^2 + G(u,v)\,{\dd}v^2. \]
\end{prop}

\begin{dem}
  All the $\vec{\gamma}_v$ are unit speed curves with the same indicator $\epsilon_{\vec{\gamma}_v}$. We have that \[ E(u,v) = \langle \vec{x}_u(u,v),\vec{x}_u(u,v)\rangle = \langle \vec{\gamma}'_v(u),\vec{\gamma}'_v(u)\rangle =
   \epsilon_{\gamma_v}.\] Now, $\epsilon_{\vec{\gamma}} \epsilon_{\vec{\gamma}_v} = (-1)^\nu$ for all $v \in I$, whence $E(u,v) = (-1)^\nu \epsilon_{\vec{\gamma}}$.

Proceeding, we see that by construction $F(0,v) = 0$ for all $v \in I$, so that it suffices to check that $F$ does not depend on the variable $u$. Fixed $v_0 \in I$, we have the expression $\vec{x}(u,v_0) = \vec{\gamma}_{v_0}(u)$, and so the second geodesic equation for $\vec{\gamma}_{v_0}$ yields $\Gamma_{11}^2(u,v_0) = 0$. From the arbitrariety of $v_0$ it follows that $\Gamma_{11}^2 = 0$. On the other hand, by definition of $\Gamma_{11}^2$ we have \[ \Gamma_{11}^2(u,v) = \frac{(-1)^\nu \epsilon_{\vec{\gamma}}}{(-1)^\nu \epsilon_{\vec{\gamma}}G(u,v) - F(u,v)^2} F_u(u,v),  \]so that $F_u(u,v) = 0$, and we conclude that $F(u,v) = 0$ for all $(u,v) \in U$, as desired.
\end{dem}

\begin{obs}
  Since $G(0,v) = \epsilon_{\vec{\gamma}} \neq 0$, the continuity of $G$ allows us to assume, by reducing $U$ again if necessary, that $G(u,v)$ has the same sign as $\epsilon_{\vec{\gamma}}$ for all $(u,v) \in U$.
\end{obs}

\begin{cor}\label{cor:curv_fermi}
  The Gaussian curvature of $(M,\pair{\cdot,\cdot})$ is expressed in terms of the Fermi chart $\vec{x}$ by \[ K \circ \vec{x} = (-1)^{\nu+1}\epsilon_{\vec{\gamma}}\frac{(\sqrt{|G|})_{uu}}{\sqrt{|G|}}. \]
\end{cor}

Before starting the proof of Theorem \ref{teo:class_riemann} (p. \pageref{teo:class_riemann}), we only need to get one more technical lemma out of the way:

\begin{lem}[Boundary conditions]
  The Fermi chart $\vec{x}$ satisfies $G_u(0,v) = 0$, for all $v \in I$.
\end{lem}

\begin{dem}
  As $\vec{\gamma}(v) = \vec{x}(0,v)$, the first geodesic equation for $\vec{\gamma}$ boils down to $\Gamma_{22}^1(0,v) = 0$, for all $v \in I$. Since $F(0,v) = 0$, it directly follows that \[ \Gamma_{22}^1(0,v) =  - \frac{G_u(0,v)}{2\epsilon_{\vec{\gamma}}}, \]whence $G_u(0,v) = 0$, as desired.
\end{dem}

Finally:

\begin{dem}[of Theorem \ref{teo:class_riemann}]
  In all possible cases, the coefficient $G$ must satisfy the following differential equation: \[  (\sqrt{|G|})_{uu} + (-1)^\nu \epsilon_{\vec{\gamma}}K\sqrt{|G|} = 0. \]Now, we solve this equation (in each case) for $\sqrt{|G|}$, and use the boundary conditions $G(0,v) = \epsilon_{\vec{\gamma}}$ and $G_u(0,v) = 0$ to determine $G$ explicitly.
  \begin{enumerate}[(A)]
  \item Assume that $\pair{\cdot,\cdot}$ is Riemannian.
    \begin{enumerate}[(i)]
    \item For $K=0$, we have $(\sqrt{G})_{uu} = 0$, and so $\sqrt{G(u,v)} = A(v)u+B(v)$. The boundary conditions then give $A(v) = 0$ and $B(v) = 1$, so that $G(u,v) = 1$ for all $(u,v) \in U$, and ${\dd}s^2 = {\dd}u^2 + {\dd}v^2$.
    \item When $K=1$, we have $(\sqrt{G})_{uu} + \sqrt{|G|} = 0$, whose solutions are of the form $\sqrt{G(u,v)} = A(v)\cos u + B(v) \sen u$. Now, the boundary conditions give $A(v) = 1$ and $B(v) = 0$, and so $G(u,v) = \cos^2u$, and it follows that ${\dd}s^2 = {\dd}u^2 + \cos^2u\,{\dd}v^2$: the metric in $\esf^2$.
    \item If $K=-1$, the equation to be solved is $(\sqrt{G})_{uu} -
    \sqrt{|G|} = 0$. We have that $\sqrt{G(u,v)} = A(v)e^u + B(v)e^{-u}$, and now the boundary conditions give $A(v) = B(v) = 1/2$, whence
    $G(u,v) = \cosh^2 u$ and we obtain the local expression ${\dd}s^2 = {\dd}u^2 +
    \cosh^2u\,{\dd}v^2$. To recognize this in an easier way as the metric in $\H^2$, we may let $x
    = e^v \tanh u$ and $y = e^v\sech u$, so that \[  {\dd}s^2 =
      \frac{{\dd}x^2 + {\dd}y^2}{y^2}, \]as desired.
    \end{enumerate}
  \item Assume now that $\pair{\cdot,\cdot}$ is Lorentzian.
 \begin{enumerate}[(i)]
    \item For $K=0$, just like above, we have ${\dd}s^2 = -{\dd}u^2 + {\dd}v^2 = {\dd}\tau^2 - {\dd}\vartheta^2$.
    \item If $K=1$, we now have two cases to discuss. If $\vec{\gamma}$ is spacelike, we again obtain $(\sqrt{G})_{uu} - \sqrt{G} = 0$, from where it follows that $G(u,v) = \cosh^2u$ and we get the $\esf^2_1$ metric (expressed in the usual revolution parametrization): ${\dd}s^2 = -{\dd}u^2 + \cosh^2u\,{\dd}v^2$.

If $\vec{\gamma}$ is timelike instead, we have $(\sqrt{-G})_{\tau \tau}+\sqrt{-G} = 0$, whose solution is $G(\tau,\vartheta) = -\cos^2\tau$, and so ${\dd}s^2 = {\dd}\tau^2 - \cos^2\tau\,{\dd}\vartheta^2$.
    \item If $K=-1$, the situation is dual to the previous one, switching ``spacelike'' and ``timelike'', and also the signs of the metric expressions. Omitting repeated calculations, we obtain \[ {\dd}s^2 = -{\dd}u^2 + \cos^2u\,{\dd}v^2= {\dd}\tau^2 - \cosh^2\tau\,{\dd}\vartheta^2,  \]which is the metric of $\H^2_1$ in suitable coordinates.
    \end{enumerate}
  \end{enumerate}
\end{dem}

We will conclude the section by presenting surfaces in the ambients $\LM^3$ and $\R^3_2$ whose metric's coordinate expressions are the ones discovered in the proof above. For $K=0$ the situation is completely uninteresting. But for $K \neq 0$ we have the following:
\newpage
\begin{Ex}\mbox{}
  \begin{enumerate}[(1)]
  \item $K=1$:
    \begin{itemize}
    \item The metric ${\dd}s^2 = -{\dd}u^2 + \cos^2u\,{\dd}v^2$ may be realized by the usual revolution parametrization $\vec{x}\colon \R^2 \to \esf^2_1 \subseteq \LM^3$ given by \[\vec{x}(u,v) = (\cosh u \cos v, \cosh u \sin v, \sinh u),\]and also by $\vec{y}\colon \cosh^{-1}\big(]1,\sqrt{2}[\big)\times \R \to \R^3_2$ given by
      \[\vec{y}(u,v) = \left(\cosh u \cosh v, \cosh u \sinh v, \int_0^u\sqrt{2-\cosh^2t}\,{\dd}t\right).\]
    \item For ${\dd}s^2 = {\dd}\tau^2 - \cos^2\tau\,{\dd}\vartheta^2$, consider $\vec{x}\colon\left]0,2\pi\right[\times\R \to \esf^2_1\subseteq\LM^3$ given by \[\vec{x}(\tau,\vartheta) = (\sin\tau, \cos\tau \cosh \vartheta,
     \cos\tau \sinh \vartheta),\]and also by $\vec{y}\colon \left]-\pi/2,\pi/2\right[\times\R\to\R^3_2$, given by \[\vec{y}(\tau,\vartheta) = \left(\int_0^\tau \sqrt{1+\sin^2t}\,{\dd}t, \cos\tau\cos\vartheta, \cos\tau\sin\vartheta\right).\]
   \begin{obs}
     The periodicity condition $\vec{y}(\tau,\vartheta) = \vec{y}(\tau+\pi,\vartheta)$ in the last given parametrization along with the fact that translations are isometries in $\R^3_2$ allow us to restrict everything to the given domains, which is maximal for non-degenerability.
    \end{obs}
    \end{itemize}
 To summarize, when $K=1$ we have the following visualizations:
  \begin{figure}[H]
    \centering
    \subfloat[In $\LM^3$]{
      \includegraphics[height=4cm]{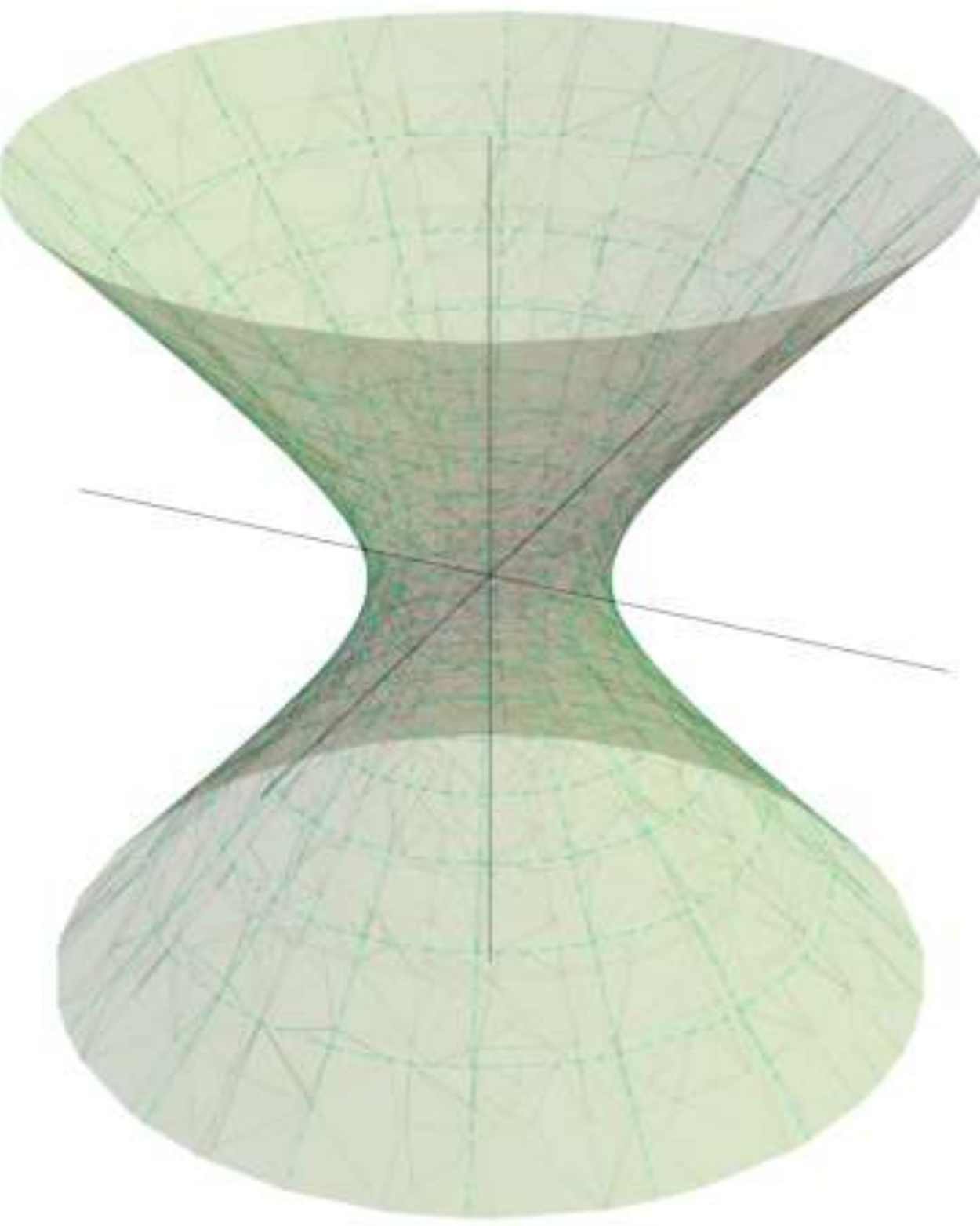}\label{fig:K1L3}}
    \ \subfloat[In $\R^3_2$]{
      \includegraphics[width=.3\textwidth]{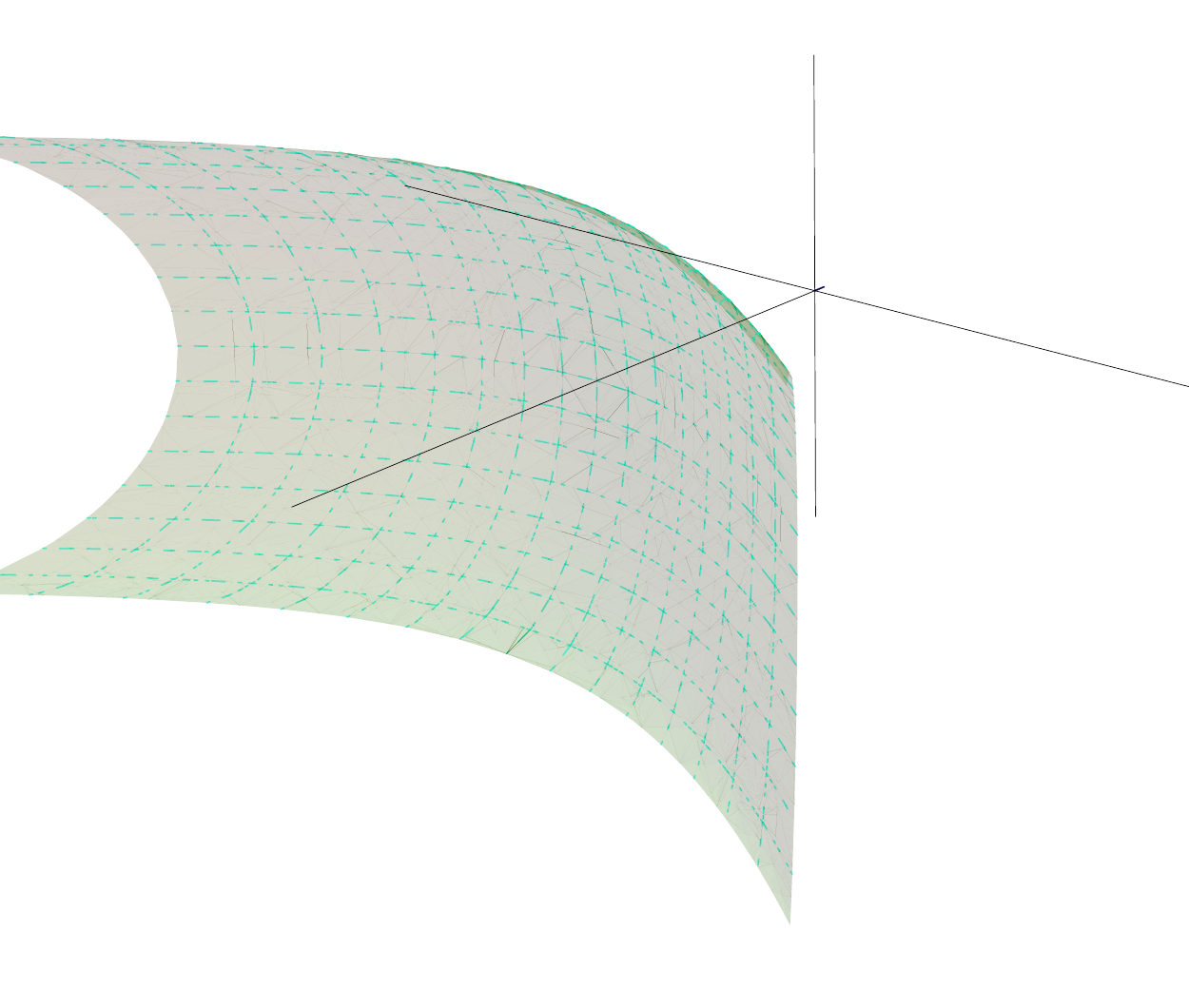} 
      \includegraphics[width=.3\textwidth]{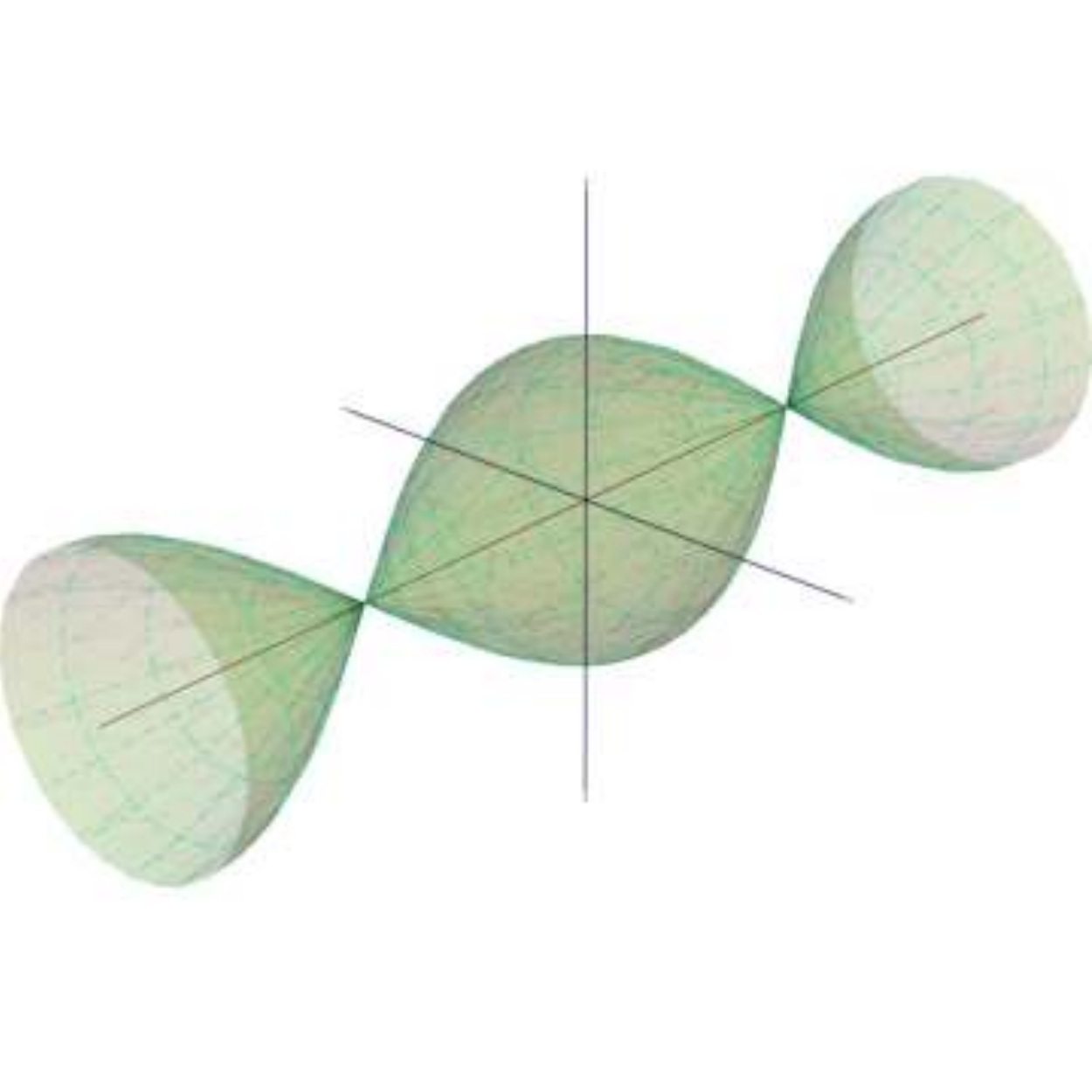}\label{fig:K1R32}}\ 
    \caption{Constant Gaussian curvature $K=1$.}
    \label{fig:K=1}
  \end{figure}

  \item $K=-1$:
    \begin{itemize}
    \item The metric ${\dd}s^2 = -{\dd}u^2 + \cos^2u\,{\dd}v^2$ may be realized by the parametrization $\vec{x}\colon\left]-\pi/2,\pi/2\right[\times\R \to
       \LM^3$, given by \[\vec{x}(u,v) = \left(\cos u\cos v, \cos u \sin v,
        \int_0^u\sqrt{1+\sin^2t}\,{\dd}t\right),\]and also by $\vec{y}\colon \left]0,2\pi\right[\times\R \to\H^2_1\subseteq\R^3_2$: \[\vec{y}(u,v) = (\cos u \sinh v, \cos u\cosh v, \sin u).\]In this case, the same remark made for $\vec{y}$ in the case $K=1$ holds for $\vec{x}$ here.
    \item The metric ${\dd}s^2 = {\dd}\tau^2 - \cosh^2\tau\,{\dd}\vartheta^2$ may be realized by the parametrization $\vec{x}\colon \cosh^{-1}\big(]1,\sqrt{2}[\big)\times \R\to\LM^3$ given by \[\vec{x}(\tau,\vartheta) = \left(\int_0^\tau \sqrt{2-\cosh^2t}\,{\dd}t,\cosh\tau\cosh\vartheta,\cosh\tau\sinh\vartheta\right)\]and by $\vec{y}\colon \R\times\left]0,2\pi\right[ \to \H^2_1 \subseteq \R^3_2$, \[\vec{y}(\tau,\vartheta)=(\sinh\tau,\cosh\tau\cos\vartheta,\cosh\tau\sin\vartheta).\]
    \end{itemize}
So in this case, we have:
  \begin{figure}[H]
    \centering
    \subfloat[In $\R^3_2$]{
      \includegraphics[width=.3\textwidth]{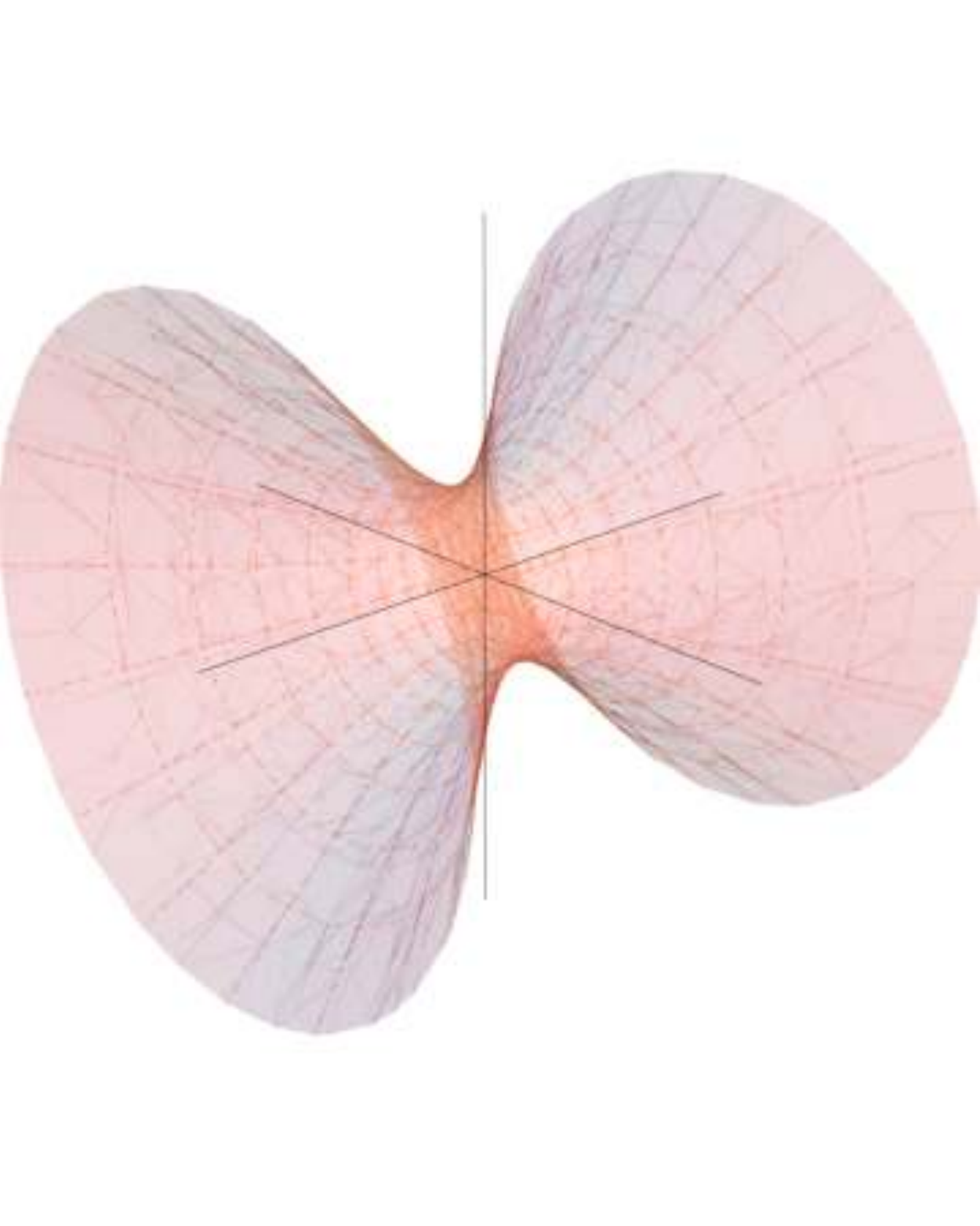}\label{fig:K-1R32}} \ 
    \subfloat[In $\LM^3$]{
      \includegraphics[width=.3\textwidth]{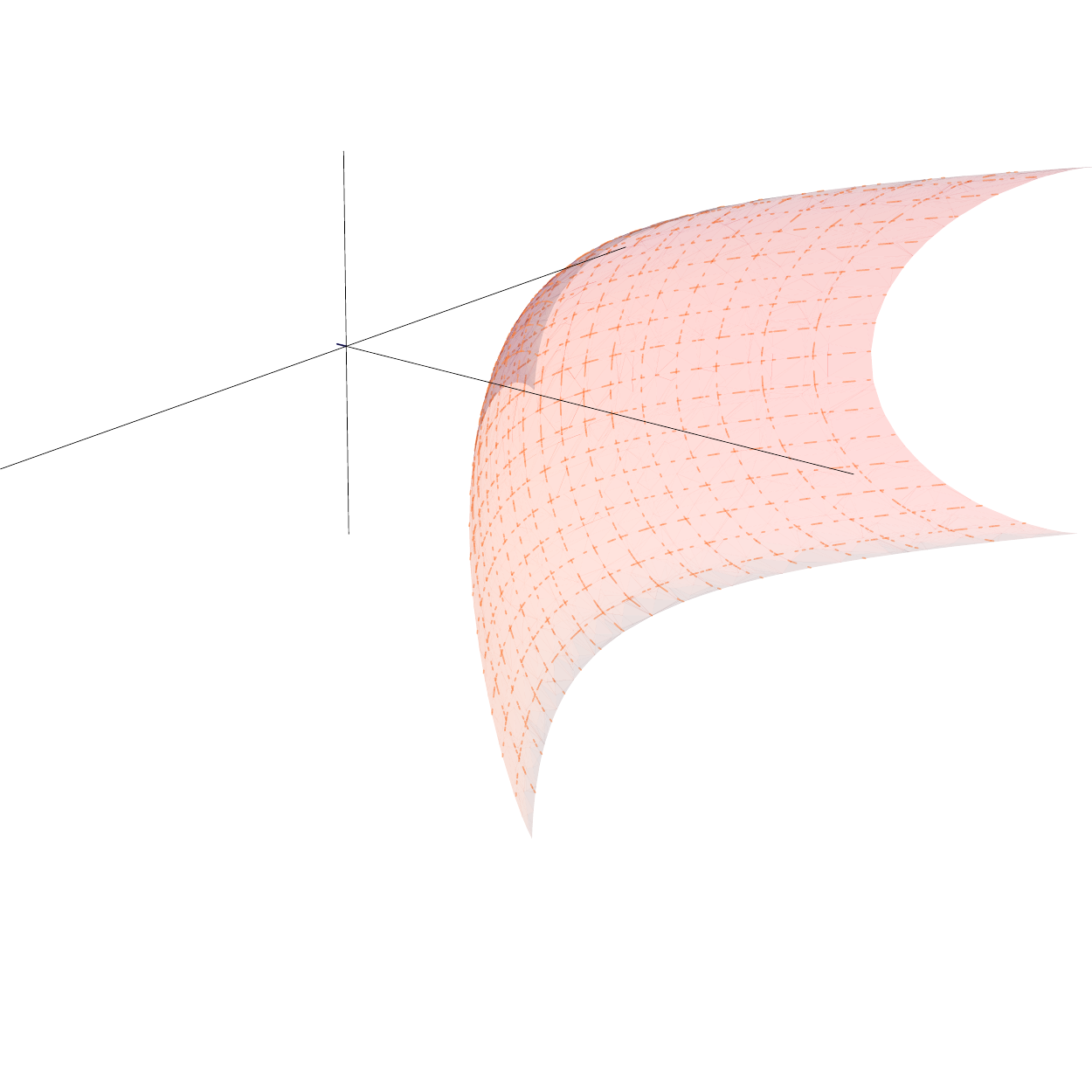} 
      \includegraphics[width=.3\textwidth]{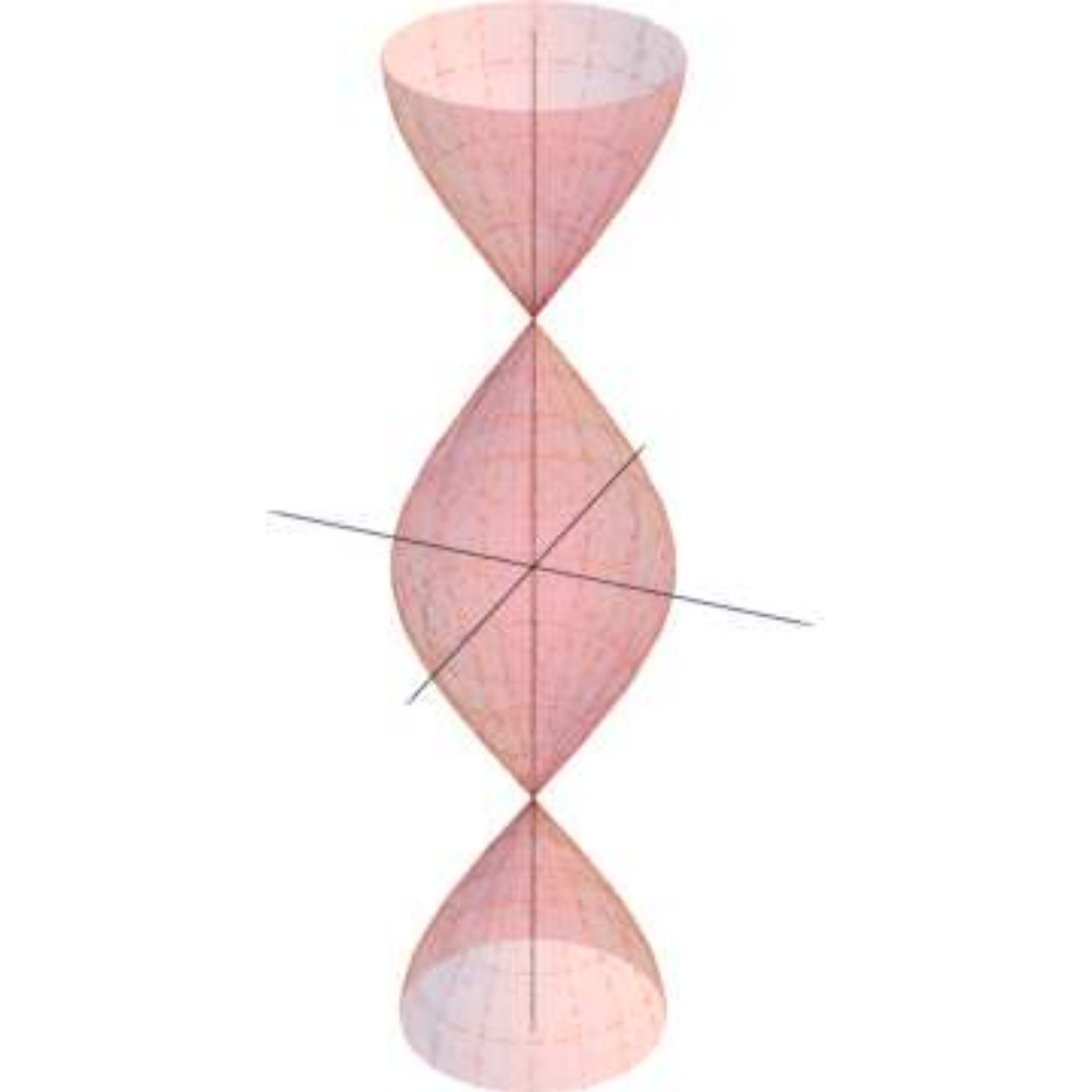}\label{fig:K-1L3}}\ 
    \caption{Constant Gaussian curvature $K=-1$.}
    \label{fig:K=-1}
  \end{figure}
  \end{enumerate}
  Lastly, we observe that the surfaces in the figures
  \ref{fig:K=1}\subref{fig:K1L3} and \ref{fig:K=-1}\subref{fig:K-1R32} are isometric when equipped by the metrics induced by $\R^3$, but on the pseudo-Riemannian ambients considered, they have rotational symmetry along axes of distinct causal characters. The same holds for the surfaces given in figures
  \ref{fig:K=1}\subref{fig:K1R32} and \ref{fig:K=-1}\subref{fig:K-1L3}. Furthermore, note that $\esf^2_1$ and $\H^2_1$ ``fit better'' in $\LM^3$ and $\R^3_2$, respectively -- switching the ambients require the use of parametrizations depending on certain elliptic integrals.
\end{Ex}

\newpage
\section*{Problems}\addcontentsline{toc}{subsection}{Problems}

\begin{problem}
  Prove Corollary \ref{cor:curv_fermi} (p. \pageref{cor:curv_fermi}).
\end{problem}

\begin{problem}
  Show that if $x = e^v \tanh u$ and $y = e^v \sech u$, then \[ \d{u}^2 + \cosh^2u\,\d{v} = \frac{\d{x}^2+\d{y}^2}{y^2}. \]
\end{problem}

   \begin{problem}[Riemann's Formula]\index{Riemann's!Formula}
    Let $(M,\pair{\cdot,\cdot})$ be a geometric surface equipped with a Riemannian metric tensor, and $\vec{x}\colon U \to \vec{x}(U) \subseteq M$ be a Fermi chart for $M$ (on which the metric is expressed by ${\dd}s^2 = {\dd}u^2 + G(u,v)\,{\dd}v^2$). In some adequate domain, consider the reparametrization $x = u \cos v$ and $y = u \sen v$. Show that \[ {\dd}s^2 = {\dd}x^2 + {\dd}y^2 + H(x,y) (x\,{\dd}y - y\,{\dd}x)^2,  \]where $H(x,y) = (G(u,v) - u^2)/u^4$.
    \begin{obs}
      The function $H$ measures, up to second order, how far is the metric from being Euclidean near the origin. The reason why is that one can show that if $H$ actually admits a continuous extension to the origin, then the Gaussian curvature at the point with coordinates $(x,y) = (0,0)$ is $-3H(0,0)$.
    \end{obs}
  \end{problem}

  \begin{problem}[Revolution surfaces with constant $K$]
    Let $\vec{\alpha}\colon I \to \R^3_\nu$ be smooth, regular, non-degenerate, injective and of the form $\vec{\alpha}(u) = (f(u),0,g(u))$, for certain functions $f$ and $g$ with $f(u) > 0$ for all $u \in I$, and let $M$ be the revolution surface spanned by $\vec{\alpha}$, around the $z$-axis. Assume that $\vec{\alpha}$ has unit speed, $M$ has constant Gaussian curvature $K$, and consider the parametrization $\vec{x}\colon I \times \left]0,2\pi\right[ \to I \to \vec{x}(U) \subseteq M$ given by \[ \vec{x}(u,v) = (f(u)\cos v, f(u) \sin v, g(u)). \]
    \begin{enumerate}[(a)]
    \item Show that, in general, $f$ and $g$ satisfy \[f''(u) + \subg{\epsilon}{\alpha} Kf(u) = 0\quad\mbox{and}\quad g(u) = \int \sqrt{(-1)^\nu(\subg{\epsilon}{\alpha} - f'(u)^2)}\,{\dd}u.\]
    \item Verify that  \[ f(u) = \begin{cases} A\cos(\sqrt{\subg{\epsilon}{\alpha}K}u)+B\sin(\sqrt{\subg{\epsilon}{\alpha}K}u), & \mbox{se }\subg{\epsilon}{\alpha}K > 0 \\ Au+B, & \mbox{se }K=0, \\ A\cosh(\sqrt{-\subg{\epsilon}{\alpha}K}u)+B\sinh(\sqrt{-\subg{\epsilon}{\alpha}K}u)& \mbox{se }\subg{\epsilon}{\alpha}K < 0,\end{cases}  \] where in the case $K = 0$ we necessarily have $|A| \leq 1$ if the ambient is $\R^3$, while $|A| \geq 1$ if the curve is spacelike in $\LM^3$ (for timelike curves there are no restrictions).
    \item Identify all the revolution surfaces with constant Gaussian curvature $K \in \{-1,0,1\}$.
    \end{enumerate}
  \end{problem}

\hrulefill

\newpage
\section*{Extra \#2: Weierstrass's representation of critical surfaces in $\LM^3$}\addcontentsline{toc}{section}{Extra \#2: Weierstrass's representation of critical surfaces in $\LM^3$}

\subsection*{An introduction to split-complex algebra}\addcontentsline{toc}{subsection}{An introduction to split-complex algebra}

We start recalling a possible construction of the complex numbers: define in $\R^2$ the operations \[ (a,b) + (c,d) \doteq (a+c,b+d) \quad\mbox{and}\quad (a,b)(c,d) \doteq (ac-bd, ad+bc).\]
Such operations turn $\R^2$ into a \emph{field}, which is then denoted by $\C$. Since we have the identities $(a,b) = (a,0) + (b,0)(0,1)$ and $(0,1)^2 = (-1,0)$, we may identify $\R$ with the set $\{ (a,0) \in \R^2 \mid a \in \R \}$ and put
$i \doteq (0,1)$, hence recovering the usual description\[ \C = \{ a+bi \mid a,b \in \R \mbox{ and } i^2=-1 \}.  \]Given $z = a+bi \in \C$, the projections ${\rm Re}(z) \doteq a$ and
${\rm Im}(z) \doteq b$ are called the \emph{real and imaginary parts of $z$}. The \emph{conjugate} of $z$ is defined as $\overline{z} \doteq a-bi$, and the \emph{absolute value} of $z$ as $|z| \doteq \sqrt{a^2+b^2} = \|(a,b)\|_E$. Moreover, if $z_1 = a_1+b_1i$ and $z_2 = a_2+b_2i$ are two complex numbers, we have \[ {\rm Re}(z_1\overline{z_2}) = \pair{(a_1,b_1),(a_2,b_2)}_E, \]which shows that $\C$ encodes the geometry of the usual inner product in $\R^2$. One then proceeds to develop Calculus in a complex variable.

 Our goal here is to define a Lorentzian version of $\C$ based on the above review, and briefly understand how calculus works in this new setting.

 \begin{defn}[Split-complex numbers]\label{def:para_complexos}
 The set $\C'$ of the \emph{split-complex numbers}\index{Split-complex!numbers} is the space $\LM^2$ equipped with the operations \[  (a,b) + (c,d) \doteq (a+c,b+d) \quad\mbox{and}\quad (a,b)(c,d) \doteq (ac+bd, ad+bc).\]
 \end{defn}
 \begin{obs}
 The split-complex numbers are also known as \emph{hyperbolic numbers}\index{Hyperbolic!numbers|seealso{split-complex numbers}}. To justify this terminology, work through Problem \ref{ex:complexos_generalizados} in the end of the section.
 \end{obs}

 It is easy to see that $\C'$ is a commutative ring with $1$. Since this time we have the identities $(a,b) = (a,0)+(b,0)(0,1)$ and $(0,1)^2 = (1,0)$, we may again identify $\R$ with $\{(a,0)\in \LM^2 \mid a \in \R \}$ and put $h \doteq (0,1)$ to obtain a similar description to the previous one given for $\C$: \[ \C' = \{ a+bh \mid a,b \in \R \mbox{ and } h^2=1 \}.  \]
 
 \begin{defn}
Let $w =a+bh \in \C'$. 
 \begin{enumerate}[(i)]
 \item The \emph{split-conjugate} of $w$ is defined by $\overline{w}\doteq a-bh$.
 \item The \emph{split-complex absolute value} of $w$ is given by $|w| \doteq \sqrt{|a^2-b^2|} = \|(a,b)\|_L$.
 \item The \emph{real part} of $w$ is given by ${\rm Re}(w) \doteq a$, and its \emph{imaginary part} is given by ${\rm Im}(w) \doteq b$.
 \end{enumerate}
 \end{defn}

Let's register some basic algebraic properties of $\C'$ in the following:

\begin{prop}\label{prop:basicos_para_comp}
  Let $w,w_1,w_2 \in \C'$.
  \begin{enumerate}[(i)]
  \item $\overline{w_1+w_2} = \overline{w_1} + \overline{w_2}$, $\overline{w_1w_2} = \overline{w_1}\,\overline{w_2}$, $\overline{\overline{w}} = w$, and $\overline{w} = w$ if and only if $w \in \R$. In fancier terms, conjugation in $\C'$ is still an involution preserving $\R$;
  \item if $1/w$ exists, then $\overline{1/w} = 1/\overline{w}$;
  \item $|w| = |\overline{w}|$, $|w\overline{w}| = |w|^2$;
  \item $|w_1w_2| = |w_1||w_2|$ and, if $1/w$ exists, $|1/w| = 1/|w|$. In particular, if $1/w$ exists, we necessarily have $|w| \neq 0$.
  \end{enumerate}
\end{prop}

 To justify that $\C'$ is indeed the Lorentzian version of $\C$ that we seek, note that if $w_1 = a_1+b_1h$ and $w_2 = a_2+b_2h$ are two split-complex numbers, then \[ {\rm Re}(w_1\overline{w_2}) = \pair{(a_1,b_1),(a_2,b_2) }_L,  \]which says that $\C'$ encodes the geometry of $\LM^2$ in the same way that $\C$ does it for $\R^2$. This also gives us a geometric interpretation for $\C'$ not being a field like $\C$: the zero divisors in $\C'$ correspond precisely to the lightlike directions in $\LM^2$.

We proceed with some calculus. We endow $\C'$ with the usual topology of the plane. That is to say, the open subsets of $\C'$ are the same ones as of $\C$, and the overall notion of continuity is the same. In particular, if $U\subseteq \C'$ is open and $f\colon U' \to \C'$ is written in the form \[ f(x+hy) = \phi(x,y) + h\psi(x,y)  \]for some real-valued functions $\phi$ and $\psi$, then $f$ is continuous if and only if both $\phi$ and $\psi$ are.

  To define holomorphicity in $\C'$, we will again mimic the definition used in $\C$, taking care to not divide by ``lightlike'' directions:
  
  \begin{defn}
  Let $U \subseteq \C'$ be an open set, $w_0 \in U$ and $f\colon U \to \C'$ a function. We'll say that $f$ is $\C'$-differentiable at $w_0$ if the limit \[  f'(w_0) \doteq \lim_{ {\Delta w \to 0}\atop{\Delta w \not\in C_L(\vec{0})}} \frac{f(w_0+\Delta w) - f(w_0)}{\Delta w}  \]exists. In this case, $f'(w_0)$ is called the \emph{derivative} of $f$ at $w_0$. And $f$ is called \emph{split-holomorphic}\index{Split-holomorphic function} in $w_0$ if it is $\C'$-differentiable in every point of some neighborhood of $w_0$.
  \end{defn}
  
  The usual rules hold:
  
  \begin{prop}
  Let $U\subseteq \C'$ be an open set, $w_0 \in U$ and $f,g\colon U \to \C'$ two \linebreak[4]$\C'$-differentiable functions at $w_0$. Then
  \begin{enumerate}[(i)]
  \item $f+g$ is $\C'$-differentiable at $w_0$ and $(f+g)'(w_0) = f'(w_0)+g'(w_0)$.
  \item $fg$ is $\C'$-differentiable at $w_0$ and $(fg)'(w_0) = g(w_0)f'(w_0) + f(w_0)g'(w_0)$.
  \item if $g$ does not assume any value in the lightlike directions of the plane, then $f/g$ is \linebreak[4]$\C'$-differentiable at $w_0$ and $(f/g)'(w_0) = (f'(w_0)g(w_0) - f(w_0)g'(w_0))/g(w_0)^2$.
  \end{enumerate}
  \end{prop}

  \begin{Ex}
    Constant functions and the identity $\C' \to \C'$ are clearly split-ho\-lo\-mor\-phic. It follows that all polynomials are split-ho\-lo\-mor\-phic, and its derivatives are given by the usual rules (e.g., the derivative of $f(w) = w^3+3w^2$ is $f'(w) = 3w^2+6w$). The same goes for rational functions, as long as the denominator does not take values in the lightlike directions of the plane.
  \end{Ex}

  \begin{prop}[Chain rule]\label{prop:cadeia_para_cplx}\index{Chain rule}
  Let $U_1,U_2 \subseteq \C'$ be open sets and $f\colon U_1 \to \C'$, $g\colon U_2 \to \C'$ be functions such that $f(U_1)\subseteq U_2$. If $f$ é $\C'$-differentiable at $w_0$ and $g$ is $\C'$-differentiable at $f(w_0)$, then $g \circ f$ is $\C'$-differentiable at $w_0$ and $(g\circ f)'(w_0) = g'(f(w_0))f'(w_0)$ holds.
  \end{prop}
  
In the usual complex calculus, we know that the real and imaginary parts of a holomorphic function must satisfy the \emph{Cauchy-Riemann equations}. In $\C'$, we should expect some sign change. Here's what we get (with almost the same proof):

 \begin{prop}[Revised Cauchy-Riemann]\index{Revised!Cauchy-Riemann equations}
  Let $U \subseteq \C'$ be an open set and fix $w_0 \in U$. If $f\colon U \to \C'$ is $\C'$-differentiable in $w_0$, and we write $f(x+hy) = \phi(x,y) + h\psi(x,y)$, then \[  \frac{\partial \phi}{\partial x}(w_0) =\frac{\partial \psi}{\partial y}(w_0) \quad\mbox{ and }\quad  \frac{\partial \phi}{\partial y}(w_0)=\frac{\partial \psi}{\partial x}(w_0)  .  \]
  \end{prop}
 
 \begin{obs}
   These revised equations may be expressed in a more concise way using split-complex versions of the so-called \emph{Wirtinger operators}\index{Wirtinger operators}:
   \[ \frac{\partial}{\partial w} \doteq
     \frac{1}{2}\left(\frac{\partial}{\partial
         x}+h\frac{\partial}{\partial y}\right) \quad \mbox{and}\quad
     \frac{\partial}{\partial \overline{w}} \doteq
     \frac{1}{2}\left(\frac{\partial}{\partial
         x}-h\frac{\partial}{\partial y}\right).  \]The revised Cauchy-Riemann equations become only  $\partial f/\partial \overline{w} =0$, in which case the formula  $f'(w)= (\partial f/\partial w)(w)$ holds.
 \end{obs}

\begin{Ex}
  Motivated by Euler's formula ${\rm e}^{x+iy} = {\rm e}^x(\cos y + i \sin y)$ in $\C$, we define $\exp_{\C'}\colon \C' \to \C'$ by $\exp_{\C'}(w) = {\rm e}^x(\cosh y +h \sinh y)$, where $w = x+hy$. We have that $\exp_{\C'}$ is split-holomorphic, with $(\exp_{\C'})' = \exp_{\C'}$. When there is no risk of confusion, one may simply write ${\rm e}^w$.
\end{Ex}

An important consequence of the revised Cauchy-Riemann equations is the analogue in $\C'$ of the well-known fact that the real and imaginary parts of a holomorphic function are harmonic. We have the:

 \begin{cor}
   Let $U \subseteq \C'$ be an open set and $f\colon U \to \C'$ a split-holomorphic function. If $f = \phi+h\psi$, then $\phi$ and $\psi$ are solutions of the \emph{wave equation}: $\square \phi = \square \psi =
   0$. Here $\square = \partial^2/\partial x^2 - \partial^2/\partial y^2$ is the \emph{wave operator}\index{d'Alembertian} $($\emph{d'Alembertian}$)$, and we say that $\phi$ and $\psi$ are
   \emph{Lorentz-harmonic}\index{Lorentz-harmonic functions}.
 \end{cor}

 \begin{obs}
 Note that $\square = 4  \dfrac{\partial}{\partial \overline{w}}\dfrac{\partial}{\partial w}$.
 \end{obs}

 Here we see another stark contrast between $\C$ and $\C'$. While it is difficult to solve explicitly the heat equation $\triangle \Phi = 0$ (an elliptic partial differential equation), there are explicit solutions for the wave equation $\square\Phi = 0$ (a hyperbolic partial differential equation). This can be used to completely classify all split-holomorphic functions with a convex domain. In particular, this gives us a good source of split-holomorphic functions. It also follows from this classification that while being holomorphic and complex-analytic are the same thing, split-holomorphic functions are not necessarily ``split-analytic'' (or even of class $\mathcal{C}^\infty$). We will not pursue this further here, but you can see the details in~\cite{TL}.

We'll conclude the discussion about differentiation stating the next two definitions, necessary for what will come later.

\begin{defn}
  Let $U\subseteq\C'$ be an open set, $w_0\in U$ and
  $f\colon U\setminus\{w_0\}\to\C'$. We'll say that $w_0$ is a \emph{pole of order $k\geq 1$}\index{Pole} of $f$ if $k$ is the \emph{least} integer for which $(w-w_0)^kf(w)$ is split-holomorphic.
\end{defn}

\begin{defn}
  Let $U\subseteq\C'$ be an open set and $P\subseteq U$ be discrete. We'll say that a split-holomorphic function
  $f\colon U \setminus P \to \C'$ is \emph{split-meromorphic in $U$}\index{Split-meromorphic function} if $P$ is precisely the set of poles of $f$.
\end{defn}

Let's also register the bare minimum we need about integration in $\C'$:

\begin{defn}
  Let $U\subseteq \C'$ be an open set, $f\colon U \to \C'$ be a continuous function and $\vec{\gamma}\colon I \to U$ a smooth curve. The \emph{integral of $f$ along $\vec{\gamma}$}\index{Integral along a curve} is defined as \[  \int_{\vec{\gamma}}f(w)\,{\dd}w \doteq \int_I  f(\vec{\gamma}(t))\vec{\gamma}'(t)\,{\dd}t. \]
\end{defn}

\begin{obs}
  This split-complex line integral can (obviously?) be expressed in terms of real line integrals. Moreover, this definition is naturally extended for piecewise smooth curves in $\C'$, and if $\vec{\gamma}$ is closed we'll just write $\oint_{\vec{\gamma}} f(w)\,\d{w}$ as usual.
\end{obs}

Probably the most important aspect of this integral is that we still have the:

 \begin{teo}[Fundamental Theorem of Calculus]\index{Fundamental Theorem!of Calculus for split-complex integrals}
 Let $U\subseteq \C'$ be an open set, \linebreak[4]$f\colon U \to \C'$ a continuous function, and $\vec{\gamma}\colon [a,b] \to U$ a piecewise smooth curve (actually $\mathcal{C}^1$ is enough). If $F\colon U \to \C'$ is a \emph{primitive} of $f$ (i.e., $F$ is split-holomorphic with $F' = f$), then \[  \int_{\vec{\gamma}} f(w)\,{\dd}w = F(\vec{\gamma}(b)) - F(\vec{\gamma}(a)). \]
 \end{teo}

In the next section, we will need some split-complex integrals to depend only on the endpoints of the curve we're integrating upon. In $\C$, we had the Cauchy-Goursat Theorem. In $\C'$ we still have the same theorem, with the same proof (e.g., using Green's Theorem):

\begin{teo}[Revised Cauchy-Goursat]\index{Revised!Cauchy-Goursat Theorem}
Let $U\subseteq \C'$ be a simply-connected open set, $f\colon U \to \C'$ be a split-holomorphic function with continuous derivative, and $\vec{\gamma}\colon [a,b] \to U$ a closed and piecewise smooth curve (again, $\mathcal{C}^1$ suffices), injective in $\left]a,b\right[$. Then \[ \oint_{\vec{\gamma}}f(w)\,{\dd}w = 0. \]
\end{teo}

\begin{obs}
  Note (again) that the assumption of $f'$ being continuous, which is automatically satisfied in $\C$, has to be explicitly stated here.
\end{obs}

The two main corollaries are:

\begin{cor}
  The integral of a split-holomorphic function (in the setting of the previous theorem) along a given curve depends only on the endpoints of the curve. In this case, we denote
  \[\int_{\vec{\gamma}} f(\omega)\,{\dd}\omega = \int_{w_0}^w
    f(\omega)\,{\dd}\omega,\] where $\vec{\gamma}$ joins $w_0$ and $w$.
\end{cor}

And also:

\begin{cor}
  Let $U \subseteq \C'$ be a simply-connected open set, $w_0 \in U$ and
  $f\colon U \to \C'$ be continuous. Then $F\colon U \to \C'$ given by
  \[ F(w) = \int_{w_0}^w f(\omega)\,{\dd}\omega \] is split-holomorphic and satisfies $F'=f$.
\end{cor}

We just need to get one more result out of our way:

\begin{defn}\label{defn:lor_conj}
  Let $U\subseteq \C'$ be an open set. Two functions $\phi,\psi\colon U \subseteq \R^2 \to\R$ are called
  \emph{Lorentz-con\-ju\-ga\-te} if $\phi_u =\psi_v$ and $\phi_v = \psi_u$. Such condition implies that both $\phi$ and $\psi$ are Lorentz-harmonic.
\end{defn}

\begin{teo}\label{teo:existe_conj}
  Let $U\subseteq \R^2\equiv \C'$ be a simply-connected set, and $\phi\colon U
  \to \R$ be a Lorentz-harmonic function. Then there exists a split-holomorphic function $f\colon U \to \C'$ which has $\phi$ as its real part. In particular, there is a Lorentz-conjugate function to $\phi$.
\end{teo}

\begin{dem}
  Define $g \doteq \phi_u+h\phi_v$. The condition $\square\phi = 0$ ensures that $g$ is split-holomorphic and, in particular, continuous, so that $U$ being simply-connected gives us the existence of a primitive $G =
  \psi+h\zeta$ for $g$. With this, $G'=g$ together with the revised Cauchy-Riemann equations for $G$ yield \[ \phi_u+h\phi_v =
    \psi_u+h\zeta_u = \psi_u+h\psi_v. \]Equating real and imaginary parts, we obtain $\psi = \phi + c$ for some $c \in \R$. So $f \doteq G - c$ is the desired function.
\end{dem}

The next step would be to look for a Cauchy-like formula in $\C'$. Unfortunately, there is not such a formula in this new setting. We are ready to move on and apply what we have seen here for surfaces in $\LM^3$ with $H=0$. For more about split-complex numbers, see for example,~\cite{An},~\cite{Ca} and~\cite{TL}.

\newpage
\subsection*{Weierstrass-Enneper representation formulas}\addcontentsline{toc}{subsection}{Weierstrass-Enneper representation formulas}

Given a spacelike surface $M\subseteq \R^3_\nu$ and a parametrization ${\vec{x}}\colon U \to \vec{x}[U]\subseteq M$, we may identify $\R^2$ with $\C$ in the usual way and use $z = u+iv$ as a parameter. We may then write
\[{\vec{x}}(z,\overline{z}) =
  (x^1(z,\overline{z}),x^2(z,\overline{z}),x^3(z,\overline{z})),\]noting the explicit dependence on the conjugate variable $\overline{z}$ is due to the fact that we don't know whether the components of $\vec{x}$ are holomorphic functions. For timelike surfaces in $\LM^3$, the parametrization domains will be identified with open subsets of $\C'$ instead. Recall from calculus in a complex variable that if
\[ \frac{\partial}{\partial z} =
  \frac{1}{2}\left(\frac{\partial}{\partial u} -
    i\frac{\partial}{\partial v}\right)\quad\mbox{and}\quad
  \frac{\partial}{\partial \overline{z}} =
  \frac{1}{2}\left(\frac{\partial}{\partial u}+i\frac{\partial}{\partial
      v}\right),  \]then the Laplacian operator can be expressed as
\[ \triangle = \frac{\partial^2}{\partial u^2}+\frac{\partial}{\partial
    v^2} = \left(\frac{\partial}{\partial u}+i\frac{\partial}{\partial
      v}\right)\left(\frac{\partial}{\partial
      u}-i\frac{\partial}{\partial v}\right) = 4\frac{\partial}{\partial
    \overline{z}}\frac{\partial}{\partial z}. \]

It is sometimes convenient to study such parametrizations as the real part of curves in $\C^3$. We then consider an extension of the product in $\R^3$ to $\C^3$, also to be denoted by
$\pair{\cdot,\cdot}_E$, defined by
\[ \pair{(z_1,z_2,z_3),(w_1,w_2,w_3)}_E =
z_1\overline{w}_1+z_2\overline{w}_2+z_3\overline{w}_3.\]

For timelike surfaces, we'll go instead to the \emph{complex Lorentzian space}\index{Complex!Lorentzian space} $\C^3_1$, with the extended product
$\pair{\cdot,\cdot}_L$ defined by
\[ \pair{(z_1,z_2,z_3),(w_1,w_2,w_3)}_L =
z_1\overline{w}_1+z_2\overline{w}_2-z_3\overline{w}_3.\]We'll maintain the usual causal character terminology used so far.

\begin{defn}
Let $U$ be a open subset of $\C$ or $\C'$, and ${\vec{x}}\colon U \to \R^3_\nu$
be a regular and non-degenerate parametrized surface.

\begin{enumerate}[(i)]
  \item 
The \emph{complex derivative}\index{Complex!derivative} of ${\vec{x}}$ is \[  \vec{\phi} \equiv
  \frac{\partial {\vec{x}}}{\partial z} \equiv {\vec{x}}_z\doteq
  \frac{1}{2}({\vec{x}}_u - i\,{\vec{x}}_v).  \]
\item  The \emph{split-complex derivative}\index{Split-complex!derivative} of ${\vec{x}}$ is \[  \vec{\psi} \equiv \frac{\partial {\vec{x}}}{\partial w} \equiv {\vec{x}}_w \doteq \frac{1}{2}({\vec{x}}_u + h\,{\vec{x}}_v).  \]
\end{enumerate}
\end{defn}

\begin{obs}
Note that $\pair{\vec{\phi},\vec{\phi}}_E = 0$ does not imply that $\vec{\phi} = \vec{0}$, since $\vec{\phi}(z,\overline{z})\in \C^3$
for all $z$, and not necessarily in $\R^3$. Same holds \emph{a fortiori} for $\vec{\psi}$.
\end{obs}

\begin{prop}\label{prop:isoterm_derluz}
  If $M\subseteq \R^3_\nu$ is a non-degenerate regular surface and
  ${\vec{x}}\colon U\to \vec{x}(U)\subseteq M$ be a parametrization of $M$, then ${\vec{x}}$ is \emph{isothermal}\index{Isothermal parametrization} (i.e., $|E|=|G|=\lambda^2$ for some smooth $\lambda$ and $F=0$) if and only if:
  \begin{enumerate}[(i)]
    \item $\pair{\vec{\phi},\vec{\phi}}_E=0$, when $M\subseteq\R^3$;
    \item $\pair{\vec{\phi},\vec{\phi}}_L=0$, when $M\subseteq\LM^3$ is spacelike;
    \item $\pair{\vec{\psi},\vec{\psi}}_L=0$, when $M\subseteq\LM^3$ is timelike.
  \end{enumerate}
\end{prop}

\begin{dem}
  Let's work through the proof when $M$ is spacelike. We have
  \[ (x^j_z)^2 = \left(\frac{1}{2}(x^j_u -i\,x^j_v)\right)^2 =
    \frac{1}{4}((x^j_u)^2 -(x^j_v)^2 - 2i x^j_ux^j_v), \]and summing over $j=1,2,3$ gives
  \[ \pair{\vec{\phi},\vec{\phi}}=\frac{1}{4}(E-G-2iF), \]so that the conclusion follows from the fact that a complex number vanishes if and only if its real and imaginary part vanish. For timelike $M$, one obtains \[ \pair{\vec{\psi},\vec{\psi}}_L = \frac{1}{4}(E+G+2hF) \]instead.
\end{dem}

\begin{lem}\label{lem:param_conf}
  Let $M\subseteq \R^3_\nu$ be a non-degenerate regular surface and
  ${\vec{x}}\colon U \to \vec{x}[U]\subseteq M$ be an isothermal parametrization of $M$. Then:
  \begin{enumerate}[(i)]
    \item $\pair{\vec{\phi},\overline{\vec{\phi}}}_E=\lambda^2/2 \neq
    0$, when $M\subseteq\R^3$;
    \item $\pair{\vec{\phi},\overline{\vec{\phi}}}_L=\lambda^2/2 \neq
    0$, when $M\subseteq\LM^3$ is spacelike;
    \item $\pair{\vec{\psi},\overline{\vec{\psi}}}_L=\epsilon_u\lambda^2/2 \neq
    0$, when $M\subseteq\LM^3$ is timelike.
  \end{enumerate}
\end{lem}

\begin{dem}
  Let's work through the proof again assuming that $M$ is spacelike:
  \[x^j_z\overline{x^j_z} = \frac{1}{4}(x^j_u -
    i\,x^j_v)(x^j_u+i\,x^j_v) = \frac{1}{4}((x^j_u)^2+(x^j_v)^2
    -2ix^j_ux^j_v), \]and summing over $j=1,2,3$ yields
  \[ \pair{\vec{\phi},\overline{\vec{\phi}}} =
    \frac{1}{4}(\lambda^2 + \lambda^2 - 2i \cdot 0) =
    \frac{\lambda^2}{2}.\]
\end{dem}

\begin{prop}
Let $M \subseteq \R^3_\nu$ be a non-degenerate regular surface and
${\vec{x}}\colon U\to \vec{x}[U]\subseteq M$ be an isothermal parametrization of $M$. Then ${\vec{x}}$ is \emph{critical}\index{Critical/minimal surface} (i.e., $H=0$) if and only if $\vec{\phi}$ is holomorphic or $\vec{\psi}$ is split-holomorphic, according to whether $M$ is spacelike or timelike, respectively.
\end{prop}

\begin{dem}
It is a straightforward consequence of the formulas \[ \frac{\partial \vec{\phi}}{\partial \overline{z}} = \frac{\partial^2{\vec{x}}}{\partial \overline{z}\partial z} = \frac{1}{4} \triangle{\vec{x}}\quad\mbox{and}\quad\frac{\partial \vec{\psi}}{\partial \overline{w}} = \frac{\partial^2{\vec{x}}}{\partial \overline{w}\partial w} = \frac{1}{4} \square{\vec{x}}.   \]
\end{dem}

In view of this last result, we may conclude that at least locally any non-degenerate critical surface may be represented by a triple:
\begin{itemize}
  \item $\vec{\phi} = (\phi^1,\phi^2,\phi^3)$ of holomorphic functions satisfying \[(\phi^1)^2 + (\phi^2)^2 + (\phi^3)^2 = 0,\] if
  $M\subseteq\R^3$;
  \item $\vec{\phi} = (\phi^1,\phi^2,\phi^3)$ of holomorphic functions satisfying \[(\phi^1)^2 + (\phi^2)^2 - (\phi^3)^2 = 0,\] if
  $M\subseteq\LM^3$ is spacelike;
  \item $\vec{\psi} = (\psi^1,\psi^2,\psi^3)$ of split-holomorphic functions satisfying \[(\psi^1)^2 + (\psi^2)^2 - (\psi^3)^2 = 0,\] if
  $M\subseteq\LM^3$ is timelike.
\end{itemize}

A natural question at this point is: given such a triple, how to recover the starting surface? The key to answering this lies in the next:

\begin{prop}
  Let $M\subseteq \R^3_\nu$ be a non-degenerate, regular and critical surface, $U$ be a simply-connected open set, and ${\vec{x}}\colon U\to\vec{x}[U]\subseteq M$ be an isothermal parametrization of $M$. Then the components of $\vec{x}$ satisfy:
    \begin{enumerate}[(i)]
      \item
      $\displaystyle{x^j(z,\overline{z}) = c_j + 2\,{\rm Re}\int_{z_0}^z
        \phi^j(\xi)\,{\dd}\xi}$, for some $z_0\in U$, if $M$ is spacelike, and
      \item
      $\displaystyle{x^j(w,\overline{w}) = c_j + 2\,{\rm Re}\int_{w_0}^w
        \psi^j(\omega)\,{\dd}\omega}$, for some $w_0\in U$, if $M$ is timelike,
    \end{enumerate}
    where $c_j\in\R$ are convenient constants.
\end{prop}

\begin{dem}
  Just for a change, let's work this time the proof when $M$ is timelike. First note that since $U$ is simply connected, ${\vec{x}}$ is isothermal and $M$ is critical, then $\vec{\psi}$ is split-holomorphic and the integrals in the statement of the proposition are path-independent. With differentials, we have:
\begin{align*}
\psi^j\,{\dd}w &= \frac{1}{2}(x^j_u +h\,x^j_v)({\dd}u+h\,{\dd}v) = \frac{1}{2}(x^j_u\,{\dd}u + x^j_v\,{\dd}v + h(x^j_v\,{\dd}u+x^j_u\,{\dd}v)) \\ \overline{\psi^j}\,{\dd}\overline{w} &= \frac{1}{2}(x^j_u -h\,x^j_v)({\dd}u-h\,{\dd}v) = \frac{1}{2}(x^j_u\,{\dd}u + x^j_v\,{\dd}v -h(x^j_v\,{\dd}u+x^j_u\,{\dd}v))
\end{align*}Adding both expressions, we obtain \[ \d{x}^j = x^j_u\,{\dd}u+x^j_v\,{\dd}v =
  \psi^j\,{\dd}w + \overline{\psi^j}\,{\dd}\overline{w} = 2\,{\rm Re}\,\psi^j\,{\dd}w, \]whence \[ x^j(w,\overline{w}) = c_j +
  2\,{\rm Re}\int_{w_0}^w \psi^j(\omega)\,{\dd}\omega,  \]for some $c_j \in \R$ and $w_0\in U$.
\end{dem}

\begin{teo}[Enneper-Weierstrass I]\label{teo:ew1esp}\index{Enneper-Weierstrass representation formulas}
  Let $U\subseteq\C$ be a simply-connected open set, $z_0 \in U$, and
  $f,g\colon U\to \C$ functions with $f$ holomorphic, $g$ meromorphic, and $fg^2$ holomorphic. Then the map ${\vec{x}}\colon U\to \R^3_\nu$ defined by
  ${\vec{x}}(z,\overline{z}) =
  (x^1(z,\overline{z}),x^2(z,\overline{z}),x^3(z,\overline{z}))$, where
  \begin{enumerate}[(i)]
    \item $\displaystyle{x^1(z,\overline{z})={\rm Re}\int_{z_0}^z
      f(\xi)(1-g(\xi)^2)\,{\dd}\xi}$,\\
    $\displaystyle{x^2(z,\overline{z})={\rm Re}\int_{z_0}^z
      if(\xi)(1+g(\xi)^2)\,{\dd}\xi}$ and,\\
    $\displaystyle{x^3(z,\overline{z})= 2\,{\rm Re}\int_{z_0}^z
      f(\xi)g(\xi)\,{\dd}\xi}$, for $\vec{x}$ in $\R^3$ or;
  \item  $\displaystyle{x^1(z,\overline{z})
    ={\rm Re}\int_{z_0}^z f(\xi)(1+g(\xi)^2)\,{\dd}\xi}$,\\
    $\displaystyle{x^2(z,\overline{z})
    = \,{\rm Re}\int_{z_0}^z if(\xi)(1-g(\xi)^2)\,{\dd}\xi}$ and,\\
    $\displaystyle{x^3(z,\overline{z})
    = 2\,{\rm Re}\int_{z_0}^z -f(\xi)g(\xi)\,{\dd}\xi}$, for $\vec{x}$ in $\LM^3$
  \end{enumerate}
  is a parametrized surface, regular in the points where the zeros of $f$ have exactly twice the order than the order of the poles of $g$, and $|g| \neq 1$ (this last condition necessary only in $\LM^3$). Furthermore, its image is a spacelike critical surface.
\end{teo}

\begin{dem}
  The conditions over $U$, $f$ and $g$ ensure that all integrals are path-indepdendent. Moreover, in $\R^3$, the complex derivative of $\vec{x}$ is precisely \[\vec{\phi}=\left(\frac{1}{2}f(1-g^2),
  \frac{i}{2}f(1+g^2),fg\right),\] which satisfies
  \[\pair{\vec{\phi},\vec{\phi}}_E=\left(\frac{1}{2}f(1-g^2)\right)^2 +
    \left(\frac{i}{2}f(1+g^2)\right)^2 + (fg)^2=0,\] so that the expression given in Proposition
  \ref{prop:isoterm_derluz} (p. \pageref{prop:isoterm_derluz}) for $\pair{\vec{\phi},\vec{\phi}}_E$ gives us that $E=G$ and $F=0$. A similar computations gives us the same conclusion in $\LM^3$. This way, $\vec{x}$ is regular precisely when $E = G \neq 0$, which is equivalent to the condition given on zeros of $f$ and poles of $f$. The necessity of $|g|\neq 1$ in $\LM^3$ follows from the fact that
  \[\vec{\phi} =\left(\frac{1}{2}f(1+g^2), \frac{i}{2}f(1-g^2), -fg
    \right),\] then
  \[\frac{\lambda^2}{2}= \pair{\vec{\phi},\overline{\vec{\phi}}}_L = \left|\frac{1}{2}f(1+g^2)\right|^2+\left|\frac{i}{2}f(1-g^2)\right|^2-\left|-fg\right|^2=
    \frac{|f|^2}{2}(1-|g|^2)^2.  \] This being the case, $\vec{x}$ is isothermal and spacelike. Furthermore, $\vec{\phi}$ is holomorphic and hence the image of $\vec{x}$ is critical.
\end{dem}

The pair $(f,g)$ is called the \emph{Weierstrass data} for the surface. The geometry of the surface can be described by this data. For more details in $\R^3$, see for example~\cite{Oss}.

When $g$ is holomorphic and invertible, we may use it as a parameter itself and obtain an alternative representation:

\begin{teo}[Enneper-Weierstrass II]\label{teo:ew2esp}
  Let $U\subseteq\C$ be a simply-connected open set, \linebreak[4]$z_0 \in U$, and
  $F\colon U\to\C$ a holomorphic function. Then the map ${\vec{x}}\colon U\to \R^3_\nu$ defined by
  ${\vec{x}}(z,\overline{z}) =
  (x^1(z,\overline{z}),x^2(z,\overline{z}),x^3(z,\overline{z}))$,
  where
  \begin{enumerate}[(i)]
    \item $\displaystyle{x^1(z,\overline{z})
    = {\rm Re}\int_{z_0}^z (1-\xi^2)F(\xi)\,{\dd}\xi}$,\\
    $\displaystyle{x^2(z,\overline{z})
    = {\rm Re}\int_{z_0}^z i(1+\xi^2)F(\xi)\,{\dd}\xi}$ and,\\
    $\displaystyle{x^3(z,\overline{z})
    = 2\,{\rm Re}\int_{z_0}^z \xi F(\xi)\,{\dd}\xi}$, for $\vec{x}$ in $\R^3$, or;
  \item $\displaystyle{ x^1(z,\overline{z})
    = {\rm Re}\int_{z_0}^z (1+\xi^2)F(\xi)\,{\dd}\xi}$,\\
  $\displaystyle{x^2(z,\overline{z})
    = \,{\rm Re}\int_{z_0}^z i (1-\xi^2)F(\xi)\,{\dd}\xi}$, and\\
  $\displaystyle{x^3(z,\overline{z}) = -2{\rm Re}\int_{z_0}^z\xi
    F(\xi)\,{\dd}\xi}$, for $\vec{x}$ and $\LM^3$,
  \end{enumerate}
  is a parametrized surface, regular in the points where $F(z)
  \neq 0$ (in $\R^3$) or $F(z) \neq 0$ and $|z|\neq 1$ (in $\LM^3$). Furthermore, its image is a spacelike critical surface.
\end{teo}

\begin{Ex}[Critical spacelike surfaces in $\R^3_\nu$]\mbox{}
  \begin{enumerate}[(1)]
    \item Enneper surface\index{Enneper surface in $\R^3$} in $\R^3$: consider the Weierstrass data $f(z)=1$ and $g(z)=z$. We obtain the parametrization
    $\vec{x}:\R^2\to\R^3$ given by
    \[\vec{x}(u,v)=\left(u-\frac {u^3}{3} + u v^2 , -v + \frac {v^3}
        {3} -u^2 v , u^2 - v^2 \right).\]
   \vspace*{-2cm} \begin{figure}[H]
      \centering
      \includegraphics[width=.4\textwidth]{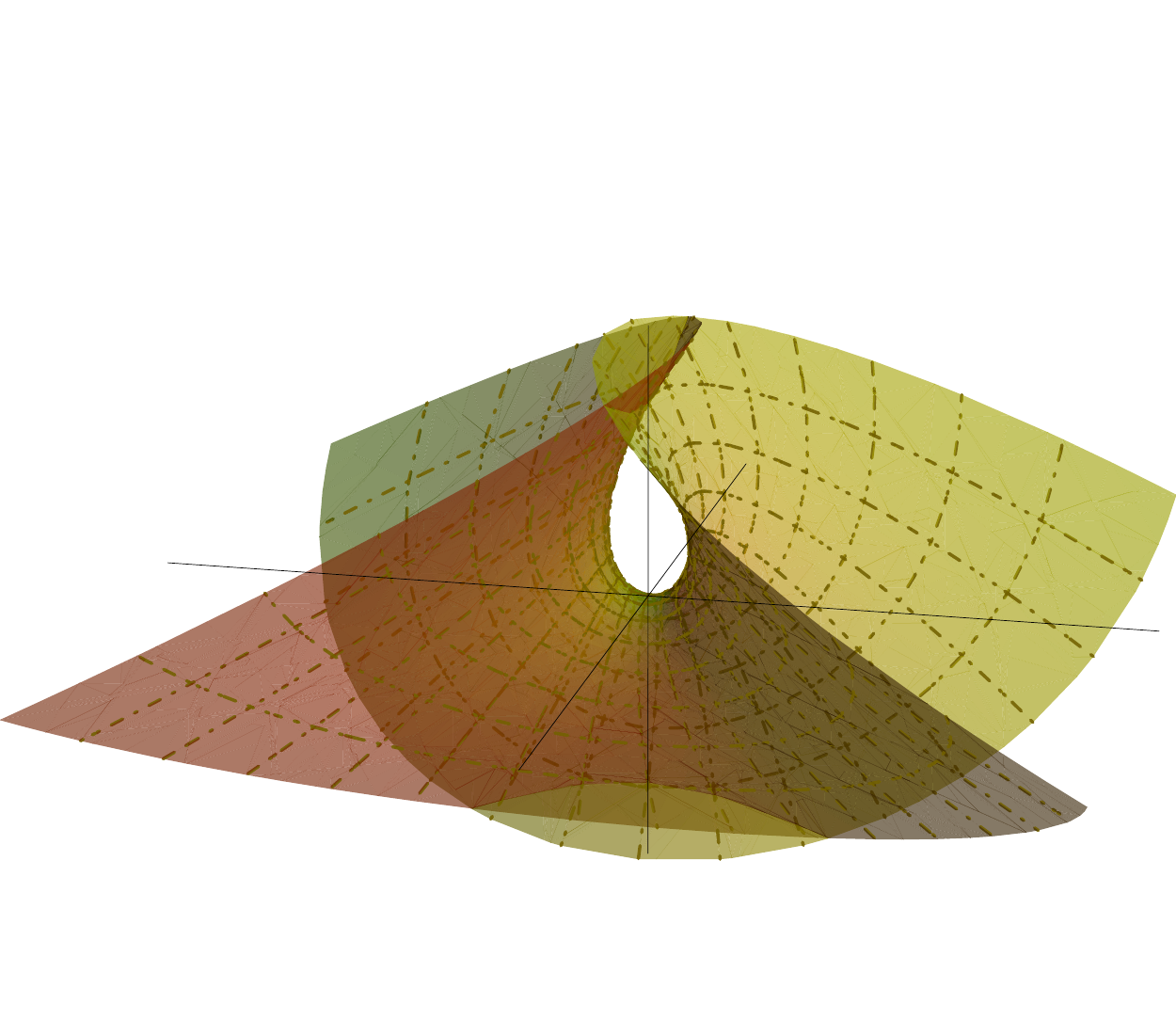}
      \caption{Enneper surface in $\R^3$.}
      \label{fig:enneperR3}
    \end{figure}
       The same surface could be obtained by the single type II data $F(z)=1$.
    \item Spacelike Enneper surface\index{Spacelike!Enneper surface in $\LM^3$}  in $\LM^3$: consider the same data as before, now in the Lorentzian setting. We obtain the parametrization $\vec{x}:\R^2\to\R^3$
    given by
    \[\vec{x}(u,v)=\left(u+\frac {u^3}{3} - u v^2 , -v - \frac {v^3}
        {3} +u^2 v , v^2 - u^2 \right).\]
    \begin{figure}[H]
      \centering
      \includegraphics[width=.4\textwidth]{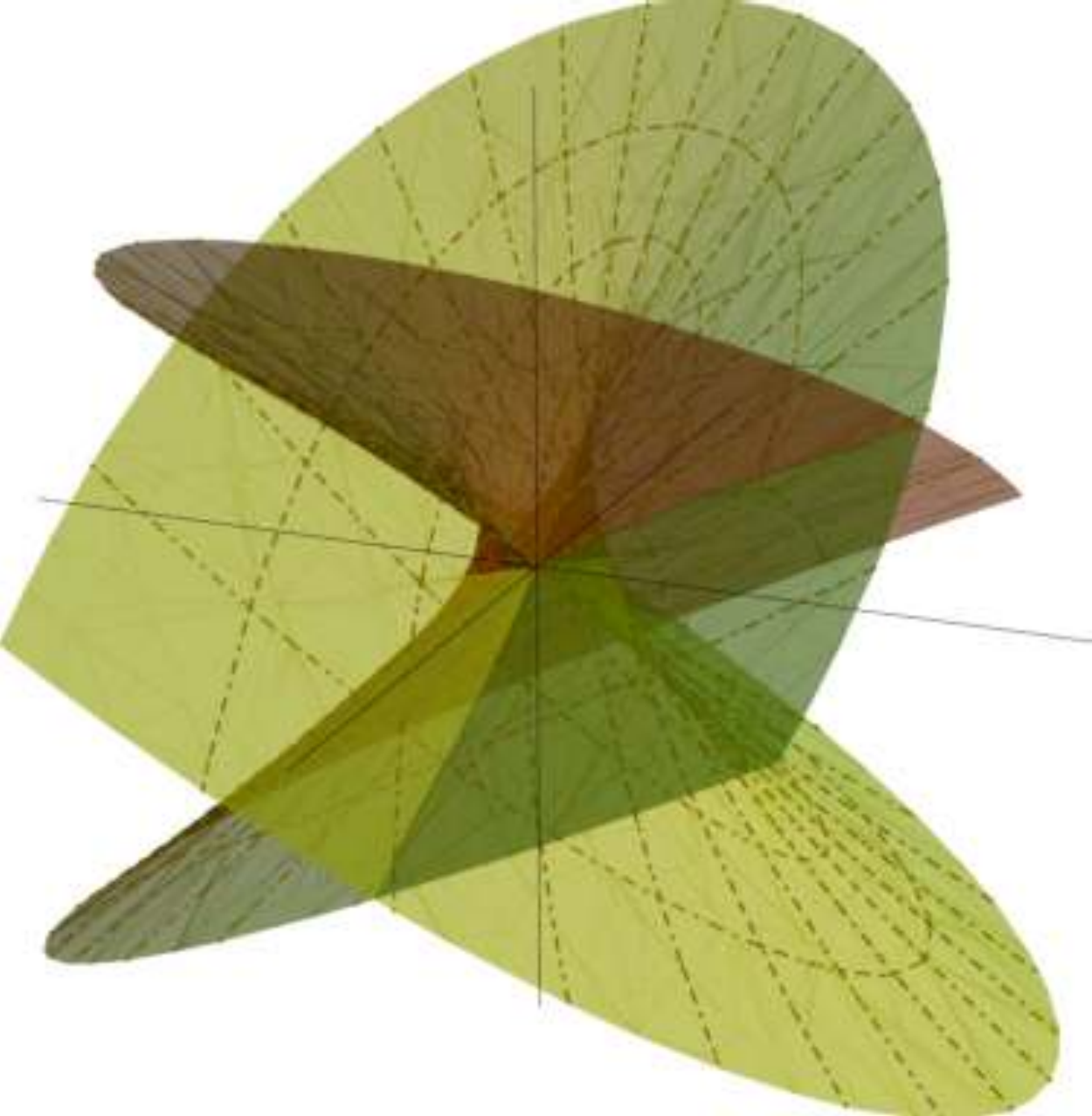}
      \caption{Spacelike Enneper surface in $\LM^3$.}
      \label{fig:enneperR3esp}
    \end{figure}
    \item Catalan surface\index{Catalan surface in $\R^3$} in $\R^3$: choose the type II data
    $F(z) = i\left(\frac{1}{z}-\frac{1}{z^3}\right)$, which yields
    \[{\vec{x}}(u,v) =\left(u - \sin u \cosh v,1-\cos u \cosh v,-4
        \sin\left(\frac{u}{2}\right)\sinh\left(\frac{v}{2}\right)\right)\]

    \begin{figure}[H]
      \centering
      \includegraphics[width=.4\textwidth]{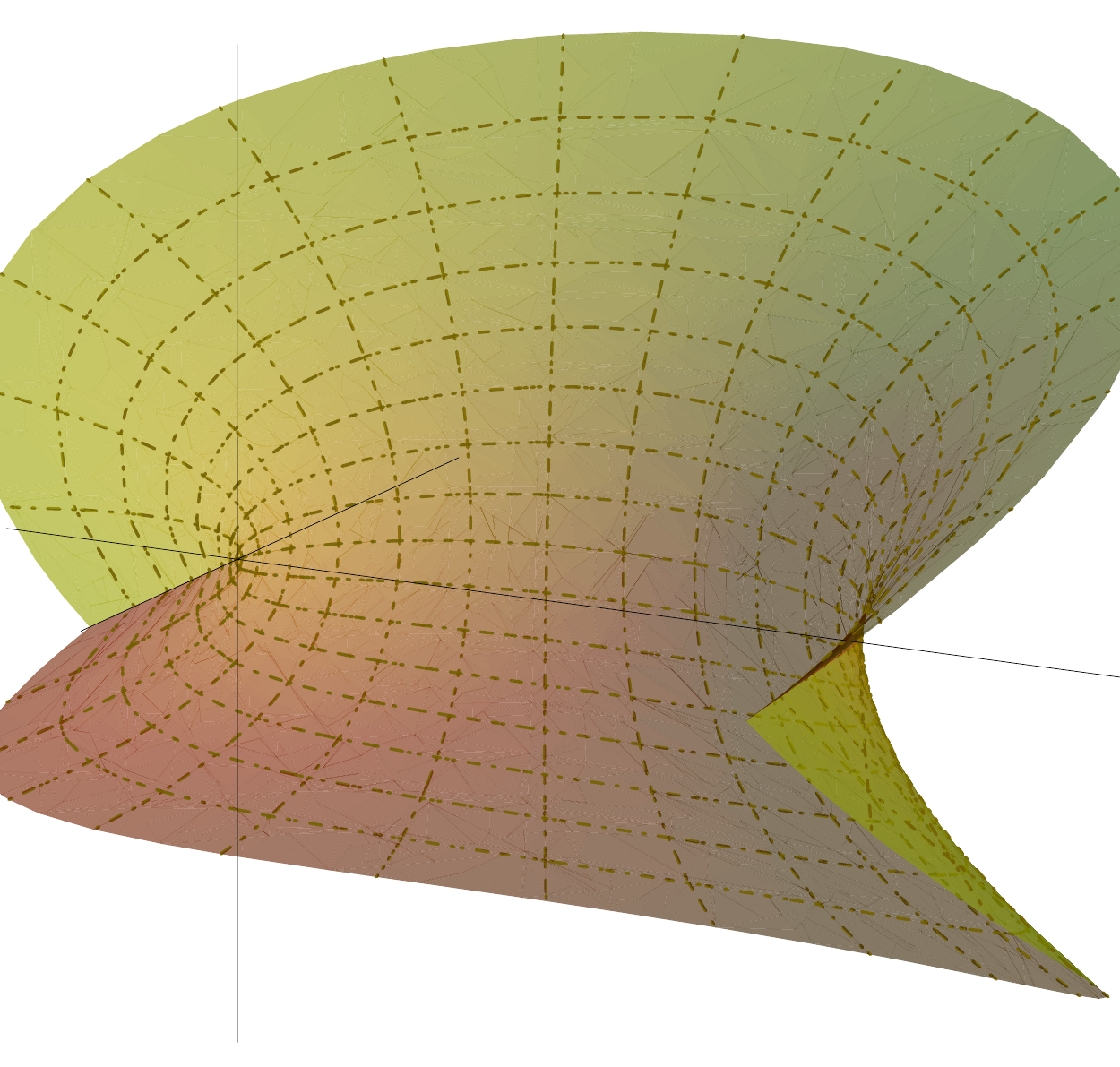}
      \caption{Catalan surface in $\R^3$.}
      \label{fig:catalanmin}
    \end{figure}
    \item Spacelike catenoid in $\LM^3$\index{Spacelike!catenoid in $\LM^3$}: choose the type II data $F(z)=1/z^2$, which gives the parametrization \[ \vec{x}(u,v) = \left(u -
        \frac{u}{u^2+v^2},v - \frac{v}{u^2+v^2}, -2\log (u^2+v^2)
      \right). \]
     \begin{figure}[H]
      \centering
      \includegraphics[width=.4\textwidth]{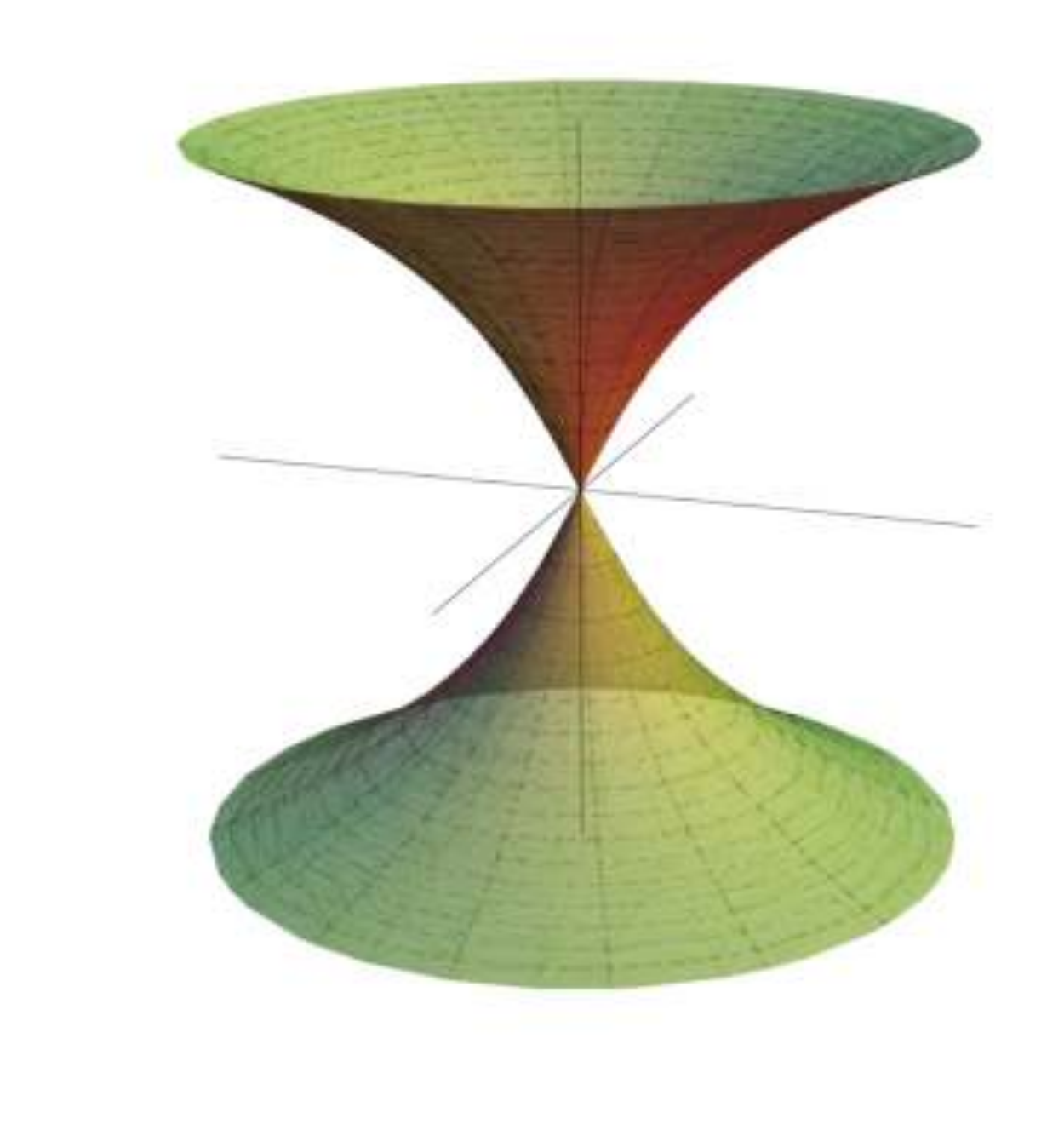}
      \caption{Spacelike Lorentzian catenoid.}
      \label{fig:catenoideespmax}
    \end{figure}
    \item Henneberg surface in $\R^3$\index{Henneberg surface in $\R^3$}: the type II data is $F(z) = 1-\frac{1}{z^4}$, producing the parametrization
    \begin{align*}
      {\vec{x}}(u&,v)=\left(2 \sinh u \cos v - (2/3) \sinh(3u) \cos(3v),\right.\\
                    &\phantom{=\left(\right.}\left.2 \sinh(u) \sin(v)+ (2/3)
                      \sinh(3u) \sin(3v), 2\cosh(2u)
                      \cos(2v)\right).
    \end{align*}
    \begin{figure}[H]	
      \centering
      \includegraphics[width=.4\textwidth]{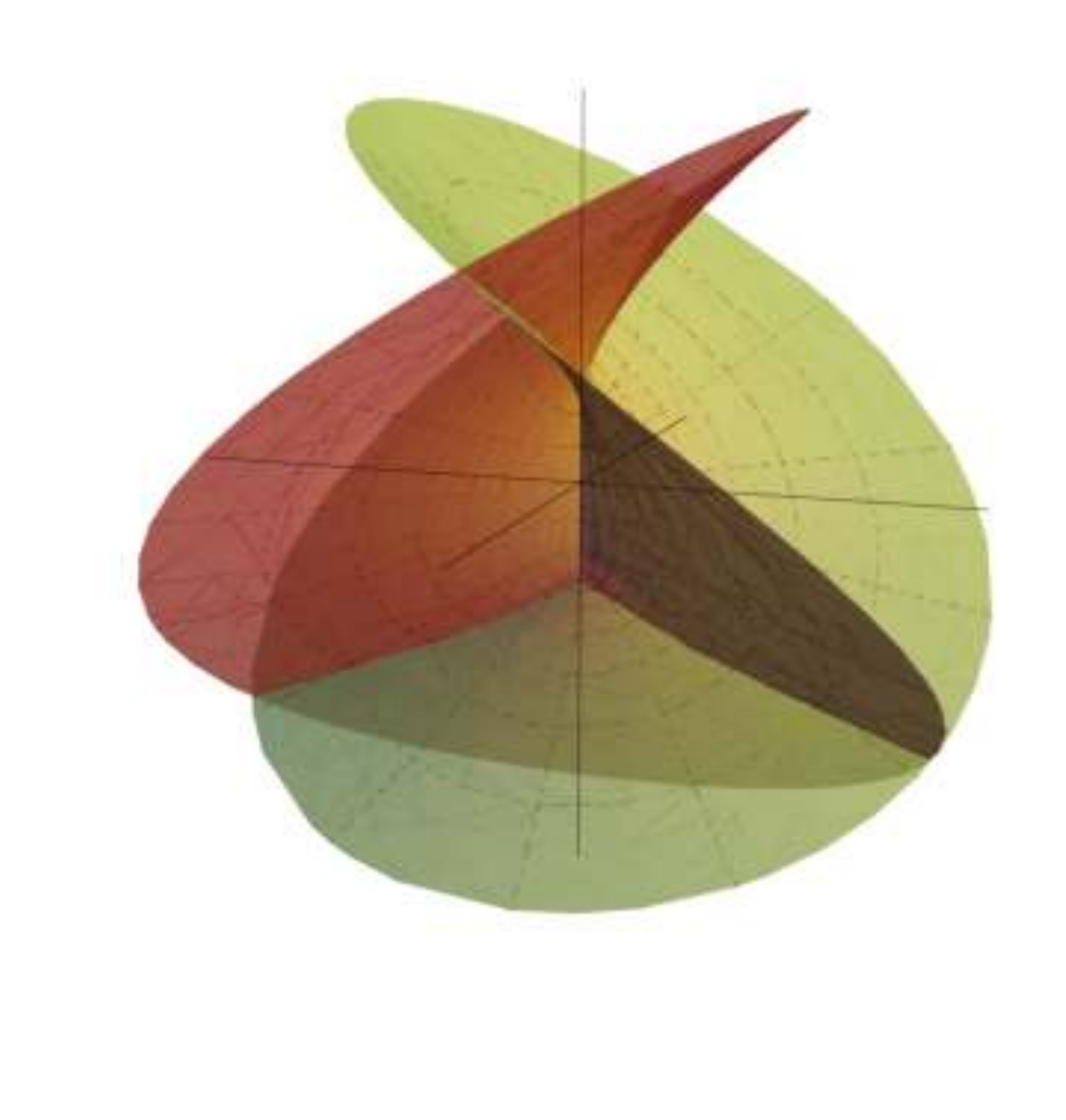}
      \caption{Henneberg surface in $\R^3$.}
      \label{fig:hennebergmin}
    \end{figure}
  \end{enumerate}
\end{Ex}

Now, let's repeat this for timelike surfaces in $\LM^3$:

\begin{teo}[Enneper-Weierstrass I -- timelike version]\label{teo:ew1temp}
  Let $U\subseteq\C'$ be a simply-connected open set, $w_0 \in U$, and
  $f,g\colon U\to \C$ functions with $f$ split-holomorphic, $g$ split-meromorphic, and $fg^2$ split-holomorphic. Then the map ${\vec{x}}\colon U\xrightarrow{\phantom{mmm}}{} \LM^3$ defined by \linebreak[4] ${\vec{x}}(w,\overline{w}) =
  (x^1(w,\overline{w}),x^2(w,\overline{w}),x^3(w,\overline{w}))$, where
  \begin{align*}
    x^1(w,\overline{w})
    &= {\rm Re}\int_{w_0}^w f(\omega)(1-g(\omega)^2)\,{\dd}\omega \\
    x^2(w,\overline{w})
    &= 2\,{\rm Re}\int_{w_0}^w  f(\omega)g(\omega)\,{\dd}\omega \\
    x^3(w,\overline{w})
    &= \,{\rm Re}\int_{w_0}^w f(\omega)(1+g(\omega)^2)\,{\dd}\omega,
  \end{align*}
  is a parametrized surface, regular in the points where the zeros of $f$ have exactly twice the order than the order of the poles of $g$, $f(w)$ is not a zero-divisor and $g(w)$ is not real. Furthermore, its image is a timelike critical surface.
\end{teo}

\begin{dem}
  The conditions over $U$, $f$ and $g$ again ensure that all the above integrals are path-independent. In this case, the split-complex derivative of $\vec{x}$ is  \[\vec{\psi}=\left(\frac{1}{2}f(1-g^2),fg,\frac{1}{2}f(1+g^2)\right).\] We also have that
  \begin{align*}
    \pair{\vec{\psi},\vec{\psi}}_L
    &=\left(\frac{1}{2}f(1-g^2)\right)^2 +
      (fg)^2-\left(\frac{1}{2}f(1+g^2)\right)^2=0\mbox{ and}\\
    \pair{\vec{\psi},\overline{\vec{\psi}}}_L
    &=\frac{f\overline{f}}{4}\left((1-g^2)(1-\overline{g}^2)+4g\overline{g}-(1+g^2)(1+\overline{g}^2)\right)=-\frac{f\overline{f}}{2}(g-\overline{g})^2,
  \end{align*}
  from where the conclusion follows.
\end{dem}

\begin{teo}[Enneper-Weierstrass II]\label{teo:ew2temp}\index{Enneper-Weierstrass representation formulas}
  Let $U\subseteq\C'$ be a simply-connected open set which does not touch the real axis, $w_0 \in U$, and $F\colon U\to\C'$ a split-holomorphic function. Then the map  ${\vec{x}}\colon U\to \LM^3$ defined by
  ${\vec{x}}(w,\overline{w}) =
  (x^1(w,\overline{w}),x^2(w,\overline{w}),x^3(w,\overline{w}))$, where
  \begin{align*}
    x^1(w,\overline{w})
    &= {\rm Re}\int_{w_0}^w (1-\omega^2)F(\omega)\,{\dd}\omega \\
    x^2(w,\overline{w})
    &=2\, {\rm Re}\int_{w_0}^w \omega F(\omega)\,{\dd}\omega \\
    x^3(w,\overline{w})
    &= {\rm Re}\int_{w_0}^w (1+\omega^2)F(\omega)\,{\dd}\omega,
  \end{align*}
  is a parametrized surface, regular in the points where $F(w)$ is not a zero divisor. Furthermore, its image is a timelike critical surface.
\end{teo}

\begin{Ex}[Timelike critical surfaces in $\LM^3$]\mbox{}
  \begin{enumerate}[(1)]
    \item Timelike Enneper surface\index{Timelike!Enneper surface}: for $F(w)=1$ we obtain
    \[\vec{x}(u,v)=\left(v-u^2v-\frac{v^3}{3},2uv,v+u^2v+\frac{v^3}{3}\right).\]
    \begin{figure}[H]
      \centering
      \includegraphics[width=.3\textwidth]{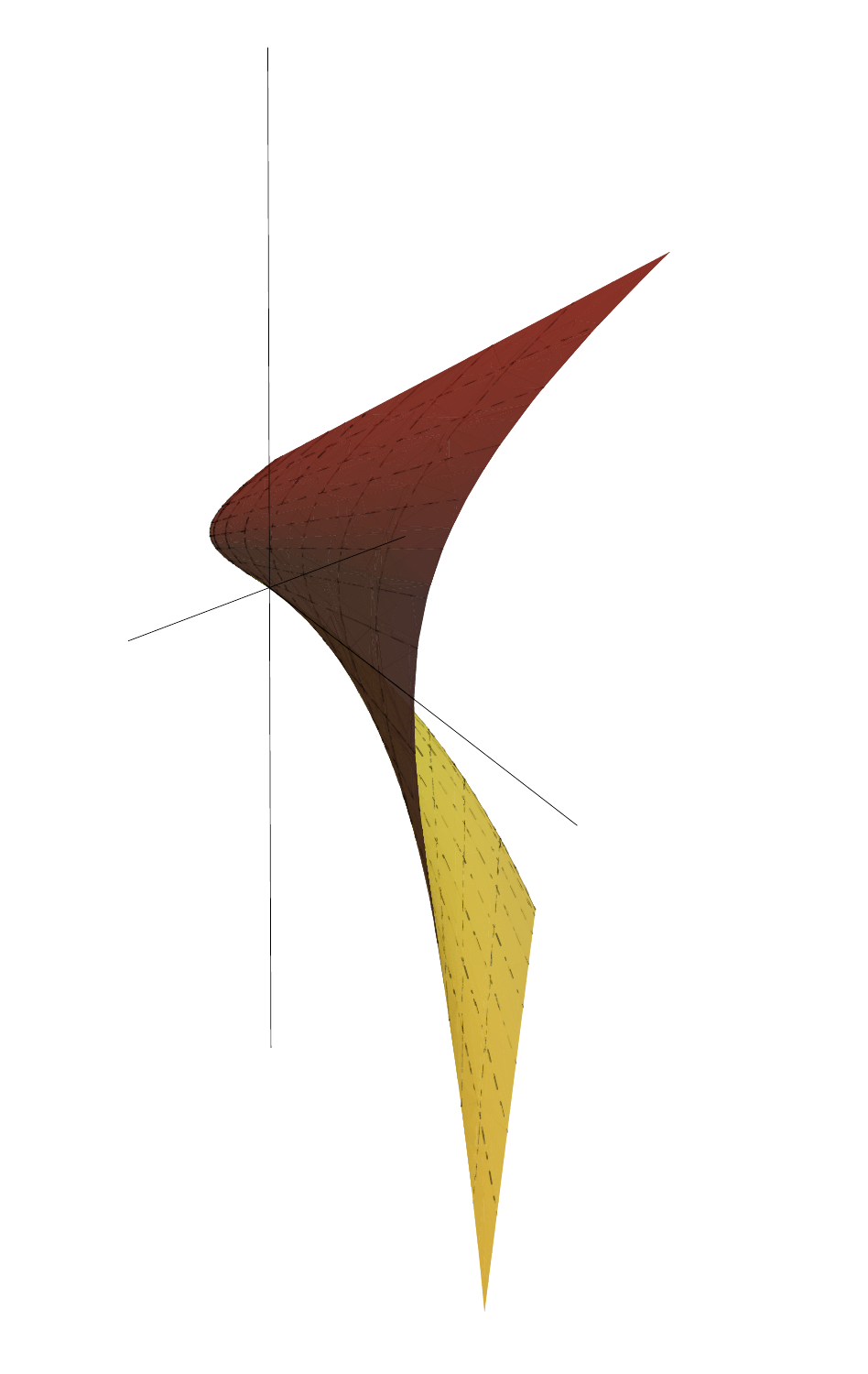}
      \caption{Timelike Enneper surface.}
      \label{fig:enneperL3tempo}
    \end{figure}
    \item Timelike catenoid\index{Timelike!catenoid}: for $F(w)=1/w^2$ we obtain
    \[\vec{x}(u,v)=\left(-\frac{u}{u^2-v^2}-u,\log\left((u^2-v^2)^2\right),-\frac{u}{u^2-v^2}+u\right),\]
    which is defined in all the plane $\LM^2$, except for the two null lines, and regular everywhere minus on the real axis ($v=0$).
    \begin{figure}[H]
      \centering
      \includegraphics[width=.4\textwidth]{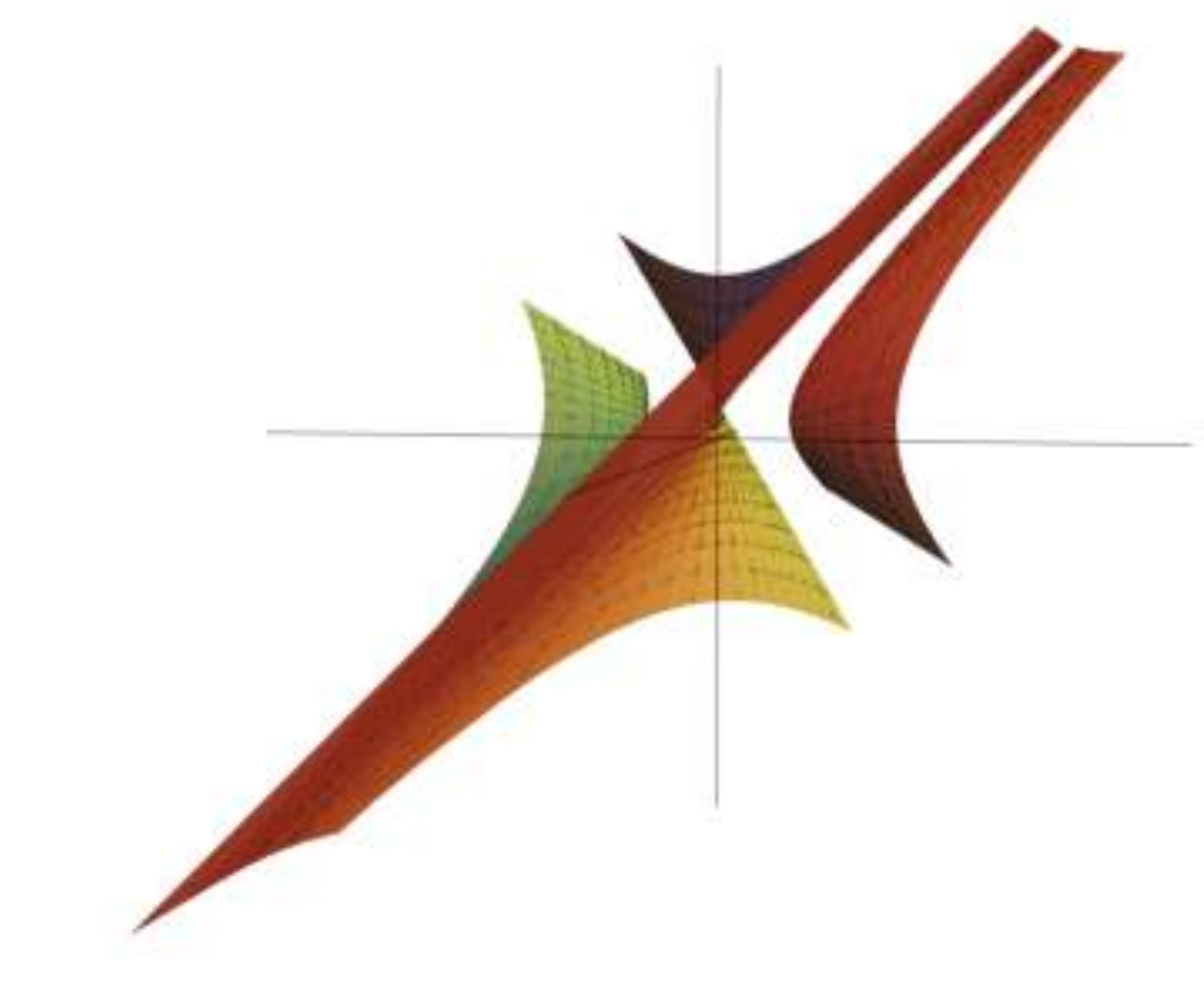}
      \caption{Timelike catenoid.}
      \label{fig:catenoideL3tempo}
    \end{figure}
  \end{enumerate}
  \item 
\end{Ex}

For more details about such representation formulas, you may consult~\cite{konderak} and~\cite{Magid}. Such techniques also have applications in the study of the so-called \emph{Bjorling problems} -- see, for example,~\cite{Alias},~\cite{Rosa},~\cite{D1} and~\cite{D2}.

\newpage
\section*{Problems}\addcontentsline{toc}{subsection}{Problems}

\begin{problem}[Generalized complex numbers]\label{ex:complexos_generalizados}
  Given \emph{structure constants}\index{Structure constants (for generalized complex numbers)} $\alpha,\beta \in \R$ and a \emph{symbol} $\mathfrak{u}$, define \[ \C_{\alpha,\beta} \doteq \{a+\mathfrak{u}b \mid a,b \in \R, \mathfrak{u}^2 = \alpha+\beta\mathfrak{u}\},  \]where the operations are defined in the obvious way.
  \begin{enumerate}[(a)]
  \item Show that $a+\mathfrak{u}b \in \C_{\alpha,\beta}$ is invertible if and only if $D\doteq a^2+\beta ab-\alpha b^2 \neq 0$.
  \item If $b \neq 0$, then $D/b^2 = 0$ may be regarded as a second degree equation in the variable $a/b$, whose discriminant is $\Delta = \beta^2+4\alpha$ (check). The position of the point $(\alpha,\beta)$ in the plane relative to the parabola $\Delta=0$ determines the possibility of realizing divisions in $\C_{\alpha,\beta}$. Show that:
    \begin{itemize}
    \item If $\Delta < 0$, $\C_{\alpha,\beta}$ is a field.
    \item If $\Delta \geq 0$, the zero divisors in $\C_{\alpha,\beta}$ are precisely the elements $a+\mathfrak{u}b$ such that $a+(\beta+\sqrt{\Delta})b/2=0$ or $a+(\beta-\sqrt{\Delta})b/2=0$, while all the other elements are invertible. 
    \end{itemize}
  \item One says that $\C_{\alpha,\beta}$ is an \emph{elliptic}, \emph{parabolic} or \emph{hyperbolic} system of numbers if $\Delta<0$, $\Delta = 0$ or $\Delta>0$, respectively. Justify this terminology by studying in terms of $\Delta$ the conic $x^2+\beta xy -\alpha y^2 = 0$ in the plane.
  \end{enumerate}
  \begin{obs}
    If one defines $\overline{a+\mathfrak{u}b} \doteq a+\beta b - \mathfrak{u}b$, $D$ is precisely the ``squared norm'' of the element $a+\mathfrak{u}b$. The map $D\colon \C_{\alpha,\beta} \to \R$ thus defined has its behavior controlled by $\Delta$. Namely, $D$ is positive-definite if $\Delta<0$, degenerate for $\Delta=0$ and indefinite for $\Delta > 0$. Polarizing $D$, we have that $\C_{\alpha,\beta}$ is an algebraic model for the geometry of the bilinear form \[  \pair{(a,b),(c,d)}_{\alpha,\beta} \doteq ac+\frac{\beta}{2}ad + \frac{\beta}{2}bc - \alpha bd  \]in $\R^2$.
  \end{obs}
\end{problem}

\begin{problem}
  Let $U\subseteq \C'$ be a connected open set, and $f\colon U \to \C'$ be a split-holomorphic function. Denote $\ell = (1+h)/2$. Show that given $s,t \in \R$, for all $w \in U$ such that $w+t\ell, w+s\overline{\ell} \in U$ we have $f(w) = \overline{\ell}f(w+t\ell)+\ell f(w+s\overline{\ell})$.
\end{problem}

\begin{problem}
  Show that $f\colon \C' \to \C'$ given by \[  f(x+hy) = \frac{1+h}{1+e^{-x}e^{-y}} \]is bounded and split-holomorphic (hence a counter-example for Liouville's Theorem in $\C'$).
\end{problem}

\begin{problem}
Let $\vec{x},\vec{y}\colon U \to \R^3_\nu$ be two regular and conjugate (or Lorentz-conjugate) parametrized surfaces. Show that if $\vec{x}$ is isothermal, then so is $\vec{y}$.
\end{problem}

\begin{problem}
  Let $\theta \in \R$. Show that the parametrized surface $\vec{x}\colon \left]0,2\pi\right[ \times \R \to \R^3$ given by \[ \vec{x}(u,v) = (u \cos \theta \pm \sin u \cosh v, v \pm \cos \theta \cos u \sinh v, \pm \sin \theta \cos u \cosh v)  \]is isothermal and minimal.
\end{problem}

\begin{problem}\label{ex:param_conf}
  Prove Lemma \ref{lem:param_conf} (p. \pageref{lem:param_conf})
  for timelike surfaces.
\end{problem}

\begin{problem}
  Let $M\subseteq \R^3_\nu$ be a critical spacelike surface and
  ${\vec{x}}\colon U\to M\subseteq \R^3_\nu$ be a type II Weierstrass parametrization defined by a holomorphic function $F$. Show that the Gaussian curvature is given by\[ K(\vec{x}(u,v)) = (-1)^{\nu+1}\frac{4}{|F(u+iv)|^2((-1)^\nu +u^2+v^2)^4}. \]
\end{problem}

 \begin{problem}
   Let $M\subseteq \R^3_\nu$ be a non-degenerate, regular and connected surface. Assume that $\vec{N}$ is a Gauss map for $M$, which is a locally conformal map. Show that $M$ is critical or is contained in a piece of a sphere, de Sitter space or hyperbolic plane.
 \end{problem}

\hrulefill

\printindex\addcontentsline{toc}{section}{Index}


\newpage

\begin{thebibliography}{AC}
\markboth{References}{References}
\addcontentsline{toc}{section}{References}




\bibitem{A} Alexandrov, D., \emph{A contribution to chronogeometry},
Canadian J. Math. 19, pp. 1119-1128, 1967.

\bibitem{Alias} Alias, L.; Chaves, R. M. B.; Mira, P.; \emph{Bjorling roblem for maximal surfaces in Lorentz-Minkowski space}, Math. Proc. Camb. Phil. Soc. (no. 134, pp. 289-316), 2003.

\bibitem{Rosa} Chaves, R. M. B.; Dussan, M. P.; Magid, M.; \emph{Bjorling Problem for timelike surfaces in the Lorentz-Minkowski space}, Journal of Mathematical Analysis and Applciations (377, no. 2, pp; 481-494), 2011. 

\bibitem{H} Anciaux, H., \emph{Minimal Submanifolds in
  Pseudo-Riemannian Geometry}, World Scientific, 2011.

\bibitem{An} Antonuccio, F., \emph{Semi--Complex Analysis }\&
\emph{Mathematical Physics (Corrected Version)}, eprint arXiv:gr-qc/9311032,
1993, \url{https://arxiv.org/pdf/gr-qc/9311032.pdf}.

\bibitem{BEE} Beem, J. K.; Ehrlich, P. E.; Easley, K. L., \emph{Global
  Lorentzian Geometry}, CRC Press, 1996.

\bibitem{Ca} Catoni et al., \emph{Geometry of Minkowski Spacetime}, Springer-Verlag (Springer Briefs in Physics), 2011.

\bibitem{Chen} Chen, B. Y.; \emph{Pseudo-Riemannian Geometry, $\delta$-invariants and Applications}, World Scientific, 2011.

\bibitem{dC1} do Carmo, M. P., \emph{Geometria Diferencial de Curvas e Superf\'{i}cies}, SBM (Universitary Texts Collection, volume 04), 2014.

\bibitem{D1} Dussan, M. P.; Franco Filho, A. P.; Magid, M.; \emph{The Bjorling Problem for timelike minimal surfaces in $\R^4_1$}, Annali di Matematica Pura ed Applicata (pp. 1-19), 2016.

\bibitem{D2} Dussan, M. P.; Magid, M.; \emph{Bjorling Problem for timelike surfaces in $\R^4_2$}, Journal of Geometry and Physics (no. 73, pp. 187-199), 2013.

\bibitem{Ig} Goldman, W. M.; Margulis, G. A., \emph{Flat Lorentz
  3-manifolds and cocompact Fuchsian groups. Crystallographic groups and
  their generalizations (Kortrijk, 1999)} Contemporary Mathematics 262, pp. 135–-145, 2000.

\bibitem{Greub} Greub, W. H., \emph{Linear Algebra}, Springer-Verlag, 1975.

\bibitem{Haw} Hawking, S., Ellis, G.; \emph{The Large Scale Structure of Spacetime}, Cambridge Monographs on Mathematical Physics, 1973.

\bibitem{Hi} Hitzer, E., \emph{Non-constant bounded holomorphic functions of hyperbolic
numbers – Candidates for hyperbolic activation functions}, Proccedings
of the First SICE Symposium on Computational Intelligence,
pp. 23--28, 2011.

\bibitem{KK} Kosheleva, O.; Kreinovich, V., \emph{Observable
  Causality Implies Lorentz Group: Alexandrov--Zeeman--Type Theorem for
  Space--Time Regions}, Mathematical Structures and Modeling 30,
pp. 4--14, 2014.

\bibitem{K} Kobayashi, O., \emph{Maximal Surfaces in the 3-Dimensional
  Minkowski Space $\LM^3$}, Tokyo J. of Math. Volume 06,
pp. 297--309, 1983.

\bibitem{konderak} Konderak, J.; \emph{A Weierstrass Representation Theorem for Lorentz Surfaces}, Complex Variables 50 (no. 5, pp.319-332), 2005.

\bibitem{RL} López, R., \emph{Differential Geometry of Curves and
  Surfaces in Lorentz--Minkowski space}, eprint arXiv:0810.3351, \url{https://arxiv.org/pdf/0810.3351}, 2008.

\bibitem{Magid} Magid, M.; \emph{Minimal Timelike Surfaces via the Split-Complex Numbers}, Proceedings of PADGE 2012, Shaker-Verlag, Aachen, 2013.

\bibitem{Na2} Naber, G. L., \emph{Spacetime and Singularities, An Introduction}, Cambridge University Press, 1988.

\bibitem{Na} Naber, G. L., \emph{The Geometry of Minkowski Spacetime: An Introduction to the Mathematics of the Special Theory of Relativity}, Springer--Verlag (Applied Mathematical Sciences 92), 1992.

\bibitem{ON2} O'Neill, B., \emph{Semi-Riemannian Geometry with Applications to Relativity}, Academic Press, 1983.

\bibitem{Oss} Osserman, R., \emph{A Survey of Minimal Surfaces}, Dover, 1986.

\bibitem{Pen} Penrose, R., \emph{Techniques of Differential Topology in Relativity}, AMS Coloquium Publications (SIAM, Philadelphia), 1972.

\bibitem{Keti} Tenenblat, K., \emph{Introduç\~{a}o \`{a} Geometria Diferencial}, Edgard Bl\"{u}cher, 2008.

\bibitem{TL} Terek, I.; Lymberopoulos, A.; \emph{Introduç\~{a}o \`{a} Geometria Lorentziana: Curvas e Superf\'{i}cies}, SBM (Universitary Texts Collection, volume 21), 2018.
\end{thebibliography}
\end{document}